\documentclass{amsart}
\usepackage{amssymb}

\newtheorem{theorem}{Theorem}

\newtheorem{conjecture}[theorem]{Conjecture}

\newtheorem{lemma}[theorem]{Lemma}

\newtheorem{proposition}[theorem]{Proposition}

\numberwithin{equation}{section}
\input{tcilatex}

\begin{document}
\title[Invariant area functionals]{Invariant surface area functionals and
singular Yamabe problem in 3-dimensional CR geometry}
\author{Jih-Hsin Cheng}
\address{Institute of Mathematics, Academia Sinica, Taipei and National
Center for Theoretical Sciences, Taipei Office, Taiwan, R.O.C.}
\email{cheng@math.sinica.edu.tw}
\thanks{2010 Mathematics Subject Classification: 32V05, 53C45, 53C17}
\author{Paul Yang}
\address{Department of Mathematics, Princeton University, Princeton, NJ
08544, U.S.A.}
\email{yang@Math.Princeton.EDU}
\thanks{}
\author{Yongbing Zhang}
\address{School of Mathematical Sciences, USTC, Hefei, Anhui, P.R.O.C.}
\email{ybzhang@amss.ac.cn, yongbing@princeton.edu}
\thanks{}
\keywords{CR invariant surface area functional, pseudohermitian geometry,
pseudohermitian torsion, Tanaka-Webster curvature, singular CR Yamabe
problem, volume renormalization}
\thanks{}

\begin{abstract}
We express two CR invariant surface area elements in terms of quantities in
pseudohermitian geometry. We deduce the Euler-Lagrange equations of the
associated energy functionals. Many solutions are given and discussed. In
relation to the singular CR Yamabe problem, we show that one of the energy
functionals appears as the coefficient (up to a constant multiple) of the
log term in the associated volume renormalization.
\end{abstract}

\maketitle

%\subjclass{32V05, 53C45, 53C17 }

%\address{\U{f8f8} }
%\email{\U{f8f8} }

\section{\textbf{Introduction and statement of the results}}

Motivated by the recent study (\cite{GrahamYamabe}) of the singular Yamabe problem and the associated
volume renormalization, we look into the analogous situation in CR
geometry. For a CR analogue of the Willmore energy in the surface case,
one of us found two CR invariant surface area elements $dA_{1},$ $%
dA_{2}$ in 1995 (see \cite{Ch}). Since there is a well developed local invariants for 
surfaces in the Heisenberg (see, for instance, 
\cite{CHMY}, \cite{CHMYCodazzi} and \cite{CCHY}), we can easily express $%
dA_{1},$ $dA_{2}$ in terms of quantities in pseudohermitian geometry. This
is done in Section 2.

To be more precise, let us review some basic notions for a nonsingular
surface $\Sigma $ in a pseudohermitian 3-manifold $(M,J,\theta ).$ We refer
the reader to \cite{CHMY} for more details. On $(M,J,\theta ),$ there is a
canonical connection $\nabla $, called Tanaka-Webster connection or
pseudohermitian connection. Associated to this connection, we have torsion $%
A_{11},$ (Tanaka-)Webster curvature $W.$ Associated to the contact form $%
\theta $, we have so called Reeb vector field $T.$ Associated to $\Sigma ,$
we have a special frame $e_{1},$ $e_{2}$ $:=$ $Je_{1}$ such that $e_{1}$ $%
\in $ $T\Sigma $ $\cap $ $\ker \theta $ and has unit length with respect to
the Levi metric $\frac{1}{2}d\theta (\cdot ,J\cdot ).$ We denote the coframe
dual to $e_{1},$ $e_{2}$ and $T$ as $e^{1},$ $e^{2}$ and $\theta .$ A
deviation function $\alpha $ on $\Sigma $ is defined so that $T+\alpha e_{2}$
$\in $ $T\Sigma .$ We defined mean curvature $H$ of $\Sigma $ in this
geometry, called $p$-mean curvature or horizontal mean curvature, so that $%
\nabla _{e_{1}}e_{1}$ $=$ $He_{2}.$

In Theorem \ref{cr-2-form} of Section 2, we obtain%
\begin{eqnarray*}
dA_{1} &=&|e_{1}(\alpha )+\frac{1}{2}\alpha ^{2}-\func{Im}A_{11}+\frac{1}{4}%
W+\frac{1}{6}H^{2}|^{3/2}\theta \wedge e^{1}, \\
dA_{2} &=&\{(T+\alpha e_{2})(\alpha ) \\
&&+\frac{2}{3}[e_{1}(\alpha )+\frac{1}{2}\alpha ^{2}-\func{Im}A_{11}+\frac{1%
}{4}W]H+\frac{2}{27}H^{3} \\
&&+\func{Im}[\frac{1}{6}W^{,1}+\frac{2i}{3}(A^{11})_{,1}]-\alpha (\func{Re}%
A_{\bar{1}}^{1})\}\theta \wedge e^{1}.
\end{eqnarray*}

\noindent (see (\ref{1-13})). We then have two energy functionals defined by%
\begin{equation*}
E_{1}(\Sigma ):=\int_{\Sigma }dA_{1}\text{ and }E_{2}(\Sigma ):=\int_{\Sigma
}dA_{2}.
\end{equation*}

\noindent The Euler-Lagrange equation for $E_{1}$  was derived in \cite%
{Ch} in terms of quantities in Cartan's theory of 3-dimensional CR geometry.
In Section 3 we express relevant quantities in terms of pseudohermitian
geometry. In Theorem \ref{Ee1} of Section 3, the Euler-Lagrange equation for 
$E_{1}$ reads%
\begin{eqnarray*}
0 &=&\frac{1}{2}e_{1}(|H_{cr}|^{1/2}\mathfrak{f)+}\frac{3}{2}%
|H_{cr}|^{1/2}\alpha \mathfrak{f} \\
&&+\frac{1}{2}sign(H_{cr})|H_{cr}|^{1/2}\{9h_{00}+6h_{11}h_{10}+\frac{2}{3}%
h_{11}^{3}\}
\end{eqnarray*}%
\noindent where $H_{cr},$ $\mathfrak{f}$ and $9h_{00}+6h_{11}h_{10}+\frac{2}{%
3}h_{11}^{3}$ are expressed in terms of pseudohermitian geometry in (\ref%
{e-1-1}), (\ref{e-10}) and (\ref{e-11}), resp. for the situation of free
torsion and constant Webster curvature. In Subsection 3.3 we provide an
alternative approach to deduce the first variation of $E_{1}.$ See (\ref%
{ScriptE1}) for the explicit formula. In Theorem \ref{Ee2} of Subsection
3.2, we state the Euler-Lagrange equation for $E_{2}$ in the case of free
torsion and constant Webster curvature. The Euler-Lagrange equation reads%
\begin{eqnarray*}
0&=&He_{1}e_{1}(H)+3e_{1}V(H)+e_{1}(H)^{2}+\frac{1}{3}H^{4} \\
&&+3e_{1}(\alpha )^{2}+12\alpha ^{2}e_{1}(\alpha )+12\alpha ^{4} \\
&&-\alpha He_{1}(H)+2H^{2}e_{1}(\alpha )+5\alpha ^{2}H^{2} \\
&&+\frac{3}{2}W(e_{1}(\alpha )+\frac{2}{3}H^{2}+5\alpha ^{2}+\frac{1}{2}W),
\end{eqnarray*}

\noindent (see (\ref{mathcalE})).

In Section 4 we provide many examples of critical points of $E_{1}$ and $%
E_{2}.$ Among others, we mention a couple of classes of closed
surfaces. In $\mathcal{H}_{1},$ shifted Heisenberg spheres defined by $%
(r^{2}+\frac{\sqrt{3}}{2}\rho _{0}^{2})^{2}+4t^{2}=\rho _{0}^{4}$ ($\rho
_{0} $ $>$ $0)$ satisfy the equation%
\begin{equation*}
H_{cr}:=e_{1}(\alpha )+\frac{1}{2}\alpha ^{2}+\frac{1}{6}H^{2}=0.
\end{equation*}

\noindent So they are minimizers for $E_{1}$ with zero energy. See Example 2
in Subsection 4.1. 

\vskip.1in\noindent
\underline{ Conjecture 1} The  shifted Heisenberg spheres are the
only closed minimizers for $E_{1}$ (with zero energy) in $\mathcal{H}_{1}.$
\vskip .1in

On the other hand, usual distance spheres (or Heisenberg spheres) defined by 
$r^{4}+4t^{2}$ $=$ $\rho _{0}^{4}$ ($\rho _{0}$ $>$ $0)$ are critical points
of higher energy level for $E_{1}.$ See the remark in the end of Example 3
of Subsection 4.1. Another interesting example is the Clifford torus in $%
S^{3}.${} It is a critical point of $E_{1}$ with positive energy. 

\vskip .1in\noindent
\underline{ Conjecture 2}
The Clifford torus is the unique minimizer among all surfaces
of torus type for $E_{1}$ up to CR automorphisms of $S^{3}.$ 
\vskip .1in
Critical points
of $E_{2}$ include vertical planes in $\mathcal{H}_{1},$ the surface defined
by $t$ $=$ $\frac{\sqrt{3}}{2}r^{2}$ in $\mathcal{H}_{1}$ and surfaces
foliated by a linear combination of $\mathring{e}_{1}$ and $\mathring{e}_{2}$
($\mathring{e}_{1}$ $:=$ $\partial _{x}+y\partial _{t},$ $\mathring{e}_{2}$ $%
:=$ $\partial _{y}-x\partial _{t}).$ We show that $E_{2}$ is unbounded
from below and above in general. See the remark in the end of Example 3 of
Subsection 4.2.

In Section 5 we study the expansion of a formal solution to the singular CR
Yamabe problem. Let $(M,J,\theta )$ be a $3$-dimensional pseudohermitian
manifold with boundary $\Sigma $ $=$ $\partial M$. Consider the conformal
change of $\theta $ $:$ $\tilde{\theta}$ $=$ $u^{-2}\theta .$ The singular
CR Yamabe problem is to find $u$ such that%
\begin{eqnarray}
\tilde{W} &=&-4\text{ on }M,  \label{SYP} \\
u &=&0\text{ on }\Sigma ,\text{ }u>0\text{ in the interior of }M  \notag
\end{eqnarray}

\noindent where $\tilde{W}$ is the Webster curvature of $(J,\tilde{\theta}).$
Consider a formal solution to (\ref{SYP}) of the following form%
\begin{equation*}
u(x,\rho )=c(x)\rho +v(x)\rho ^{2}+w(x)\rho ^{3}+z(x)\rho ^{4}+l(x)\rho
^{5}\log \rho+h(x)\rho ^{5} +O(\rho ^{6})
\end{equation*}

\noindent where $x$ is a regular (or nonsingular) point of $\Sigma $ and $%
\rho $ is a suitably chosen defining function for $\Sigma .$ We can
determine $c(x)$ $\equiv $ $1$ easily. In Section 5 we give explicit
expressions (\ref{v}), (\ref{w}) and (\ref{z}) for $v(x),$ $w(x)$ and $z(x),$
resp.. The coefficient $l(x)$ of the first log term is related to $\mathcal{E%
}_{2}$ and discussed in Section 6 about the volume renormalization.

In Section 6 we consider the volume renormalization for a formal solution to
(\ref{SYP}) as follows:%
\begin{equation*}
Vol(\{\rho >\epsilon \})=c_{0}\epsilon ^{-3}+c_{1}\epsilon
^{-2}+c_{2}\epsilon ^{-1}+L\log \frac{1}{\epsilon }+V_{0}+o(1)
\end{equation*}

\noindent (see (\ref{volumerenormalization})). We deduce explicit formulas
for the coefficients $c_{0},$ $c_{1},$ $c_{2}$ and $L.$ See (\ref{c0c1c2L}).
We show that 
\begin{equation*}
L(\Sigma )=2E_{2}(\Sigma )
\end{equation*}%
\noindent for $\Sigma $ being a closed, nonsingular surface (see (\ref{E2L}%
)). Finally we prove that 
\begin{equation*}
l(x)=\frac{1}{5}\mathcal{E}_{2}
\end{equation*}

\noindent (see (\ref{obstructionlandmathcalE})). This is the CR analogue of
the result for the singular Yamabe problem (\cite{GoverWaldron},\cite{GrahamYamabe}). In
particular, $\mathcal{E}_{2}\neq 0$ is also an obstruction to the smoothness
of solutions to the singular CR Yamabe problem.

\bigskip

\textbf{Acknowledgements}. J.-H. Cheng would like to thank the Ministry of
Science and Technology of Taiwan, R.O.C. for the support of the project:
MOST 106-2115-M-001-013- and the National Center for Theoretical Sciences
for the constant support. P. Yang would like to thank the NSF of the U.S.
for the support of the grant: DMS 1509505. Y. Zhang would like to thank 
Proessor Alice Chang for her invitation and the Department of Mathematics 
of Princeton University for its hospitality. He is supported by CSC scholarship
201606345025 and the Fundamental Research Funds for the Central Universities. 
P. Yang thanks Sean Curry for posing this question.

\textit{\bigskip }

\section{Two CR invariant surface area elements}

In \cite{Ch}, the first author constructed two $CR$ invariant area elements
on a nonsingular (noncharacteristic) surface $\Sigma $ in a strictly
pseudoconvex $CR$ $3$-manifold $(M,\xi ,J)$. Here $\xi $ denotes a contact
bundle and $J$ $:$ $\xi \rightarrow \xi $ is an endomorphism such that $%
J^{2} $ $=$ $-Id.$ We recall (\cite{CHMY}) that a point $p$ $\in $ $\Sigma $
is called singular if its tangent plane $T_{p}\Sigma $ coincides with the
contact plane $\xi _{p}$ at $p.$ We call a surface nonsingular if it
contains no singular points. When given a contact form $\theta $ on a $CR$
manifold, we can talk about pseudohermitian geometry (see \cite{We} or \cite%
{Lee}). In this section, we are going to express those two $CR$ invariant
area elements in terms of pseudohermitian geometric quantities.

First we recall in \cite{CHMY} that a moving frame associated to $\Sigma $
is chosen. Take $e_{1}$ $\in $ $T\Sigma \cap \xi $ of unit length with
respect to the Levi metric. Let $e_{2}$ $:=$ $Je_{1}.$ Let $T$ denote the
Reeb vector field uniquely determined by the condition: $\theta (T)$ $=$ $1,$
$T\lrcorner $ $d\theta $ $=$ $0.$ Let $(e^{1},$ $e^{2},$ $\theta )$ be the
coframe dual to $(e_{1},$ $e_{2},$ $T).$ It follows that%
\begin{eqnarray}
d\theta &=&2e^{1}\wedge e^{2}  \label{1-(-1)} \\
&=&i\theta ^{1}\wedge \theta ^{\bar{1}}  \notag
\end{eqnarray}

\noindent where $\theta ^{1}$ $:=$ $e^{1}+ie^{2}.$ So $(\theta ,\theta
^{1},\theta ^{\bar{1}};\phi =0)$ is an admissible $CR$ coframe satisfying
the following Levi equation%
\begin{equation}
d\theta =i\theta ^{1}\wedge \theta ^{\bar{1}}+\theta \wedge \phi  \label{1-0}
\end{equation}

\noindent (see \cite{CM} or \cite{We}). To apply the theory in \cite{Ch}, we
need to choose an admissible $CR$ coframe $(\theta ,\tilde{\theta}^{1},%
\tilde{\theta}^{\bar{1}};\tilde{\phi})$ such that $\tilde{\theta}^{1}$ $=$ $%
\omega ^{1}+i\omega ^{2}$ for $\omega ^{1},$ $\omega ^{2}$ real has the
property that%
\begin{equation}
\omega ^{2}=0\text{ on }\Sigma .  \label{1-1}
\end{equation}

\noindent Write%
\begin{equation}
\tilde{\theta}^{1}=\theta ^{1}u_{1}^{1}+\theta v^{1}  \label{1-2}
\end{equation}

\noindent for some complex valued functions $u_{1}^{1},$ $v^{1}$ $:=$ $%
v_{r}^{1}+iv_{c}^{1}$ where $v_{r}^{1}$ and $v_{c}^{1}$ denote real and
imaginary parts of $v^{1}$, resp.. Take $u_{1}^{1}$ $=$ $1.$ Comparing the
imaginary part of both sides in (\ref{1-2}), we obtain%
\begin{equation}
\omega ^{2}=e^{2}+\theta v_{c}^{1}  \label{1-3}
\end{equation}

\noindent Recall that the function $\alpha $ on $\Sigma $ is defined so that 
$T+\alpha e_{2}$ $\in $ $T\Sigma .$ So we compute%
\begin{eqnarray*}
0 &=&\omega ^{2}(T+\alpha e_{2})\text{ (by }(\ref{1-1})) \\
&=&v_{c}^{1}+\alpha \text{ (by (\ref{1-3})).}
\end{eqnarray*}

\noindent So we should take 
\begin{equation*}
v_{c}^{1}=-\alpha .
\end{equation*}

\noindent Now the real version of (\ref{1-2}) reads%
\begin{eqnarray}
\omega ^{1} &=&e^{1}+v_{r}^{1}\theta  \label{1-4} \\
\omega ^{2} &=&e^{2}-\alpha \theta  \notag
\end{eqnarray}

\noindent We may then take $v_{r}^{1}$ $=$ $0$ for simplicity. To make $%
(\theta ,\tilde{\theta}^{1},\tilde{\theta}^{\bar{1}};\tilde{\phi})$ satisfy
the Levi equation (\ref{1-0}), we need to take%
\begin{eqnarray}
\tilde{\phi} &=&-2v_{r}^{1}e^{2}-2\alpha e^{1}  \label{1-5} \\
&=&-2\alpha e^{1}  \notag
\end{eqnarray}

\noindent since $v_{r}^{1}$ $=$ $0.$ (\ref{1-4}) and (\ref{1-5}) give a
coframe transformation. We can then relate associated connection forms
according to a transformation formula 
\begin{equation}
\tilde{\Pi}=dh\cdot h^{-1}+h\text{ }\Pi \text{ }h^{-1}  \label{TF}
\end{equation}%
\noindent (see (1.3) in Part II of \cite{Ch} or \cite{CM}) where%
\begin{equation*}
\Pi =\left( 
\begin{array}{ccc}
\pi _{0}^{0} & \theta ^{1} & 2\theta \\ 
-i\phi ^{\bar{1}} & \phi _{1}^{1}+\pi _{0}^{0} & 2i\theta ^{\bar{1}} \\ 
-\frac{1}{4}\psi & \frac{1}{2}\phi ^{1} & -\overline{\pi _{0}^{0}}%
\end{array}%
\right)
\end{equation*}

\noindent with $\pi _{0}^{0}$ $=$ $(-1/3)(\phi _{1}^{1}+\phi )$ and 
\begin{equation*}
h=\left( 
\begin{array}{ccc}
t & 0 & 0 \\ 
t_{1} & t_{1}^{1} & 0 \\ 
\tau & \tau ^{1} & \bar{t}^{-1}%
\end{array}%
\right)
\end{equation*}%
\noindent subject to condition (1.4) in Part II of \cite{Ch}. The entries of 
$h$ are related to the change of admissible unitary coframes by (1.6) in
Part II of \cite{Ch}. For (1.5) in Part II of \cite{Ch} in our situation, we
have%
\begin{equation*}
u\equiv 1,u_{1}^{1}\equiv 1\text{ and }v^{1}=-i\alpha .
\end{equation*}

\noindent From (1.6) and (1.4) in part II of \cite{Ch}, we get $t^{3}$ $=$ $%
1.$ So $t$ is a constant $1$ or $e^{\frac{2\pi i}{3}}$ or $e^{\frac{4\pi i}{3%
}}.$ We may take $t\equiv 1.$ It follows that%
\begin{eqnarray}
t_{1}^{1} &=&t=1,  \label{1-6} \\
t_{1} &=&-\alpha ,\tau ^{1}=\frac{i\alpha }{2},\tau _{c}=-\frac{\alpha ^{2}}{%
4}  \notag
\end{eqnarray}

\noindent where $\tau _{c}$ is the imaginary part of $\tau .$ In general, we
will denote the real and imaginary parts of a complex-valued one-form (or
function) $\eta $ by $\eta _{r}$ and $\eta _{c},$ respectively. We now look
at the $(1,1)$ entry of the matrix transformation formula (\ref{TF}).
Comparing the corresponding imaginary parts of both sides, we obtain%
\begin{equation}
-\frac{1}{3}\tilde{\phi}_{1c}^{1}=-\frac{1}{3}\phi _{1c}^{1}+\alpha e^{2}-%
\frac{1}{2}\alpha ^{2}\theta .  \label{1-6-1}
\end{equation}

\noindent Here we have used (\ref{1-6}) several times. Recall (see \cite{We}
or page 227 in \cite{CL}) that 
\begin{equation}
\phi _{1}^{1}=\omega _{1}^{1}+\frac{i}{4}W\theta  \label{1-6-2}
\end{equation}

\noindent where $\omega _{1}^{1}$ and $W$ are pseudohermitian connection
form and Tanaka-Webster scalar curvature, resp.. Write $\omega
_{1}^{1}=i\omega $ where $\omega $ is a real one-form since $\omega _{1}^{1}$
is purely imaginary. So from (\ref{1-6-2}) we have%
\begin{eqnarray}
\phi _{1c}^{1} &=&\omega +\frac{1}{4}W\theta  \label{1-7} \\
&=&\omega (e_{1})e^{1}+\omega (e_{2})e^{2}+(\omega (T)+\frac{1}{4}W)\theta 
\notag \\
&=&He^{1}+(\omega (\alpha e_{2}+T)+\frac{1}{4}W)\theta  \notag
\end{eqnarray}

\noindent on $\Sigma ,$ where we have used 
\begin{equation}
\omega (e_{1})=H,  \label{1-7-1}
\end{equation}%
\noindent the $p$(or horizontal)-mean curvature (see page 136 in \cite{CHMY}%
) and 
\begin{equation}
e^{2}=\alpha \theta  \label{1-7-2}
\end{equation}%
\noindent on $\Sigma $ due to (\ref{1-1}) and (\ref{1-4}). Substituting (\ref%
{1-7}) into (\ref{1-6-1}) and using (\ref{1-7-2}), we obtain%
\begin{equation*}
\tilde{\phi}_{1c}^{1}=He^{1}+\{\omega (T+\alpha e_{2})+\frac{1}{4}W-\frac{3}{%
2}\alpha ^{2}\}\theta
\end{equation*}%
\noindent on $\Sigma .$ So the coefficients $h_{11}$ and $h_{10}$ of CR
second fundamental form (see (1.9.1) in \cite{Ch}) for $\tilde{\phi}%
_{1c}^{1} $ read%
\begin{eqnarray}
h_{11} &=&H  \label{1-8} \\
h_{10} &=&\omega (T+\alpha e_{2})+\frac{1}{4}W-\frac{3}{2}\alpha ^{2}  \notag
\\
&=&e_{1}(\alpha )+\frac{1}{2}\alpha ^{2}-\func{Im}A_{11}+\frac{1}{4}W. 
\notag
\end{eqnarray}

\noindent For the last equality in (\ref{1-8}), we have used an
integrability condition obtained from $e_{1},$ $T+\alpha e_{2}$ $\in $ $%
T\Sigma :$%
\begin{equation}
\omega (T+\alpha e_{2})=e_{1}(\alpha )+2\alpha ^{2}-\func{Im}A_{11}
\label{1-8-1}
\end{equation}

\noindent As a by-product of looking at $(1,1)$ entry of (\ref{TF}), we also
obtain%
\begin{equation}
\tau _{r}\text{ (the real part of }\tau )\text{ }=0  \label{1-9}
\end{equation}%
\noindent through comparing the corresponding real parts of both sides. In
fact, (\ref{1-9}) is determined by the condition (1.2) in Part II of \cite%
{Ch}), which, in our situation, reads%
\begin{eqnarray}
\phi _{1r}^{1} &=&\frac{1}{2}\phi =0,  \label{1-9-1} \\
\tilde{\phi}_{1r}^{1} &=&\frac{1}{2}\tilde{\phi}=-\alpha e^{1}  \notag
\end{eqnarray}

\noindent by (\ref{1-5}). We then look at the $(2,1)$ entry of (\ref{TF}).
Comparing the corresponding real parts of both sides gives%
\begin{equation}
\tilde{\phi}_{c}^{1}=d\alpha +\frac{1}{2}\alpha ^{2}e^{1}+\phi _{c}^{1}.
\label{1-10}
\end{equation}

\noindent Recall (see page 227 in \cite{CL}) that we have 
\begin{equation}
\phi ^{1}=A_{\bar{1}}^{1}\theta ^{\bar{1}}+\frac{i}{4}W\theta
^{1}+E^{1}\theta  \label{1-10-1}
\end{equation}

\noindent where $A_{\bar{1}}^{1}$ is the pseudohermitian torsion and 
\begin{equation*}
E^{1}=\frac{1}{6}W^{,1}+\frac{2i}{3}(A^{11})_{,1}.
\end{equation*}

\noindent It follows from (\ref{1-10-1}) that%
\begin{eqnarray}
\phi _{c}^{1} &=&[(\func{Im}A_{\bar{1}}^{1})+\frac{1}{4}W]e^{1}-(\func{Re}A_{%
\bar{1}}^{1})e^{2}  \label{1-11} \\
&&+\func{Im}[\frac{1}{6}W^{,1}+\frac{2i}{3}(A^{11})_{,1}]\theta .  \notag
\end{eqnarray}

\noindent So the coefficients $h_{00}$ of CR second fundamental form (see
(1.9.2) in \cite{Ch}) for $\tilde{\phi}_{c}^{1}$ reads%
\begin{eqnarray}
h_{00} &=&(T+\alpha e_{2})(\alpha )  \label{1-12} \\
&&+\func{Im}[\frac{1}{6}W^{,1}+\frac{2i}{3}(A^{11})_{,1}]-\alpha (\func{Re}%
A_{\bar{1}}^{1})  \notag
\end{eqnarray}

\noindent on $\Sigma $ by substituting (\ref{1-11}) into (\ref{1-10}). Here
we have used (\ref{1-7-2}) and%
\begin{equation}
e^{1}=\omega ^{1}  \label{1-12-0}
\end{equation}%
\noindent due to (\ref{1-4}) and $v_{r}^{1}$ $=$ $0$. Let%
\begin{eqnarray}
dA_{1} &:&=|h_{10}+\frac{1}{6}h_{11}^{2}|^{3/2}\theta \wedge \omega ^{1},
\label{1-12-1} \\
dA_{2} &:&=(h_{00}+\frac{2}{3}h_{10}h_{11}+\frac{2}{27}h_{11}^{3})\theta
\wedge \omega ^{1}  \notag
\end{eqnarray}

\noindent be those two CR invariant 2-forms in (1.13) of Part II in \cite{Ch}%
. By (\ref{1-8}) and (\ref{1-12}) we can express $dA_{1}$ and $dA_{2}$ in
terms of quantities in pseudohermitian geometry.

\bigskip

\begin{theorem}
\label{cr-2-form} For a nonsingular surface $\Sigma $ in a pseudohermitian
3-manifold $M$ with contact form $\theta $, we have%
\begin{eqnarray}
dA_{1} &=&|e_{1}(\alpha )+\frac{1}{2}\alpha ^{2}-\func{Im}A_{11}+\frac{1}{4}%
W+\frac{1}{6}H^{2}|^{3/2}\theta \wedge e^{1},  \label{1-13} \\
dA_{2} &=&\{(T+\alpha e_{2})(\alpha )  \notag \\
&&+\frac{2}{3}[e_{1}(\alpha )+\frac{1}{2}\alpha ^{2}-\func{Im}A_{11}+\frac{1%
}{4}W]H+\frac{2}{27}H^{3}  \notag \\
&&+\func{Im}[\frac{1}{6}W^{,1}+\frac{2i}{3}(A^{11})_{,1}]-\alpha (\func{Re}%
A_{\bar{1}}^{1})\}\theta \wedge e^{1}.  \notag
\end{eqnarray}

In particular, for $M$ being the Heisenberg group $\mathcal{H}_{1}$ of
dimension 3 or the standard pseudohermitian 3-sphere $\mathcal{S}^{3}$ (with 
$W\equiv 2),$ we have%
\begin{eqnarray}
dA_{1} &=&|e_{1}(\alpha )+\frac{1}{2}\alpha ^{2}+\frac{1}{4}W+\frac{1}{6}%
H^{2}|^{3/2}\theta \wedge e^{1},  \label{1-14} \\
dA_{2} &=&\{(T+\alpha e_{2})(\alpha )  \notag \\
&&+\frac{2}{3}[e_{1}(\alpha )+\frac{1}{2}\alpha ^{2}+\frac{1}{4}W]H+\frac{2}{%
27}H^{3}\}\theta \wedge e^{1}  \notag
\end{eqnarray}

\noindent where $W$ vanishes for $\mathcal{H}_{1}$ and equals constant $2$
for $\mathcal{S}^{3}.$
\end{theorem}

\bigskip

We then have two energy (or area) functionals defined by%
\begin{equation}
E_{1}(\Sigma ):=\int_{\Sigma }dA_{1}\text{ and }E_{2}(\Sigma ):=\int_{\Sigma
}dA_{2}.  \label{1-15}
\end{equation}

\noindent Note that the integral can only be taken over nonsingular region
of $\Sigma .$ So we assume $\Sigma $ has measure 0 singular set. The 2-forms
in (\ref{1-13}) or (\ref{1-14}) are CR invariant, i.e. invariant under the
contact form change. We write down the transformation laws of $e_{1},$ $%
e_{2},$ $T,$ and $\alpha ,$ $H$ under the change of contact form, $\tilde{%
\theta}$ $=$ $\lambda ^{2}\theta ,$ $\lambda $ $>$ $0:$%
\begin{eqnarray*}
\tilde{e}_{1} &=&\lambda ^{-1}e_{1}, \\
\tilde{e}_{2} &=&\lambda ^{-1}e_{2}, \\
\tilde{T} &=&\lambda ^{-2}T+\lambda ^{-3}(e_{2}(\lambda )e_{1}-e_{1}(\lambda
)e_{2}), \\
\tilde{\alpha} &=&\lambda ^{-1}\alpha +\lambda ^{-2}e_{1}(\lambda ), \\
\tilde{H} &=&\lambda ^{-1}H-3\lambda ^{-2}e_{2}(\lambda ), \\
Im\widetilde{A}_{11}&=&\lambda^{-2}ImA_{11}+\frac{1}{2}[\lambda^{-4}e_2(%
\lambda)^2-\lambda^{-4}e_1(\lambda)^2
+\lambda^{-1}(\lambda^{-1})_{22}-\lambda^{-1}(\lambda^{-1})_{11}], \\
\widetilde{W}&=&2\lambda^{-1}[(\lambda^{-1})_{11}+(\lambda^{-1})_{22}]
-4[\lambda^{-4}e_2(\lambda)^2+\lambda^{-4}e_1(\lambda)^2]+\lambda^{-2}W.
\end{eqnarray*}

\noindent for the reader's reference. The conformal invariance of $dA_1$ can
be verified directly by using the transformation laws.

\bigskip

\section{Euler-Lagrange equations for $E_{1}$ and $E_{2}$}

\subsection{Euler-Lagrange equation for $E_{1}$}

Let us make a brief review about the Euler-Lagrange equation for $E_{1}$
(see (\ref{1-12-1}), (\ref{1-15})), deduced in Section 2 of Part II in \cite%
{Ch}, by the Cartan-Chern approach of using pure differential forms. Let 
\begin{equation}
H_{cr}:=h_{10}+\frac{1}{6}h_{11}^{2}  \label{e-1}
\end{equation}

\noindent (note that the notation $H$ has been reserved for the $p$-(or
horizontal) mean curvature). Recall ((\ref{1-15}), (\ref{1-12-1}) and (\ref%
{1-14})) that in the torsion free case, the first energy functional reads%
\begin{eqnarray*}
E_{1}(\Sigma ) &:&=\int_{\Sigma }|H_{cr}|^{3/2}\theta \wedge e^{1} \\
&=&\int_{\Sigma }|e_{1}(\alpha )+\frac{1}{2}\alpha ^{2}+\frac{1}{4}W+\frac{1%
}{6}H^{2}|^{3/2}\theta \wedge e^{1}
\end{eqnarray*}

\noindent where 
\begin{equation}
H_{cr}=e_{1}(\alpha )+\frac{1}{2}\alpha ^{2}+\frac{1}{4}W+\frac{1}{6}H^{2}
\label{e-1-1}
\end{equation}%
\noindent On a nonsingular surface $\Sigma $, we define $h_{111},$ $h_{110},$
and $h_{100}$ by the following equations:%
\begin{eqnarray}
dh_{11}+3\tilde{\phi}_{r}^{1}-h_{11}\tilde{\phi}_{1r}^{1} &=&h_{111}\omega
^{1}+h_{110}\theta ,  \label{e-2} \\
dh_{10}-h_{11}\tilde{\phi}_{r}^{1}-h_{10}\tilde{\phi} &=&h_{101}\omega
^{1}+h_{100}\theta  \notag
\end{eqnarray}

\noindent (cf. (2.8.1), (2.8.2) in Section 2 of Part II in \cite{Ch}) with $%
h_{110}$ $=$ $h_{101},$ where $\tilde{\phi}_{r}^{1},$ $\tilde{\phi}%
_{1r}^{1}, $ and $\tilde{\phi}$ are as in the previous section$.$ Assume $%
H_{cr}$ $\neq $ $0$ (a solution surface for the Euler-Lagrange equation of $%
E_{1}$ with nonzero energy). Let%
\begin{eqnarray}
\mathfrak{f} &\mathfrak{:}&=|H_{cr}|^{-1}\{h_{10}h_{111}+\frac{1}{3}%
h_{11}^{2}h_{111}  \label{e-3} \\
&&+h_{11}h_{110}+\frac{3}{2}h_{100}\}.  \notag
\end{eqnarray}

\noindent From (2.16) in Section 2 of Part II in \cite{Ch}, the
Euler-Lagrange equation for $E_{1}$ reads%
\begin{eqnarray}
0 &=&\mathfrak{-}\frac{1}{2}e_{1}(|H_{cr}|^{1/2}\mathfrak{f)-}\frac{3}{2}%
|H_{cr}|^{1/2}\alpha \mathfrak{f}  \label{e-4} \\
&&-\frac{1}{2}sign(H_{cr})|H_{cr}|^{1/2}\{9h_{00}+6h_{11}h_{10}+\frac{2}{3}%
h_{11}^{3}\}.  \notag
\end{eqnarray}

\noindent Here we have used $\mathfrak{\tilde{\phi}(}e_{1})$ $=$ $-2\alpha $
since $\mathfrak{\tilde{\phi}}$ $=$ $-2\alpha e^{1}$ by the second formula
of (\ref{1-9-1}).

We can express all the quantities involved in (\ref{e-4}) in terms of
pseudohermitian geometry. First to find out the expression of $\tilde{\phi}%
_{r}^{1},$ we compare the imaginary part of $(2,1)$ entry in both sides of (%
\ref{TF}) to get%
\begin{equation}
\tilde{\phi}_{r}^{1}=\phi _{r}^{1}+\frac{3}{2}\alpha ^{2}e^{2}-\alpha \phi
_{1c}^{1}-\frac{1}{2}\alpha ^{3}\theta .  \label{e-5}
\end{equation}

\noindent From (\ref{1-10-1}) we have%
\begin{equation}
\phi _{r}^{1}=(\func{Re}A_{\bar{1}}^{1})e^{1}+(\func{Im}A_{\bar{1}}^{1}-%
\frac{1}{4}W)e^{2}+E_{r}^{1}\theta  \label{e-6}
\end{equation}

\noindent where%
\begin{equation*}
E_{r}^{1}=\func{Re}[\frac{1}{6}W^{,1}+\frac{2i}{3}(A^{11})_{,1}].
\end{equation*}

\noindent Substituting (\ref{e-6}) and (\ref{1-7}) into (\ref{e-5}), we
obtain%
\begin{eqnarray}
\tilde{\phi}_{r}^{1} &=&(\func{Re}A_{\bar{1}}^{1}-\alpha H)e^{1}+(\func{Im}%
A_{\bar{1}}^{1}-\frac{1}{4}W+\frac{3}{2}\alpha ^{2})e^{2}  \label{e-7} \\
&&+[E_{r}^{1}-\alpha (\omega (T+\alpha e_{2})+\frac{1}{4}W)-\frac{1}{2}%
\alpha ^{3}]\theta .  \notag
\end{eqnarray}

\noindent For the torsion free ($A_{\bar{1}}^{1}$ $=$ $0$) case, we can
reduce (\ref{e-7}) to%
\begin{eqnarray}
&&\tilde{\phi}_{r}^{1}  \label{e-8} \\
&=&-\alpha He^{1}+\{-\frac{1}{2}\alpha W+\alpha ^{3}+\frac{1}{6}\func{Re}%
W^{,1}-\alpha \omega (T+\alpha e_{2})\}\theta  \notag \\
&=&-\alpha He^{1}+\{-\frac{1}{2}\alpha W-\alpha ^{3}+\frac{1}{6}\func{Re}%
W^{,1}-\alpha e_{1}(\alpha )\}\theta  \notag
\end{eqnarray}

\noindent on $\Sigma $ where we have used (\ref{1-7-2}). In the last
equality of (\ref{e-8}), we have used (\ref{1-8-1}) with free torsion: $%
\omega (T+\alpha e_{2})$ $=$ $e_{1}(\alpha )+2\alpha ^{2}.$ Substituting (%
\ref{1-8}), (\ref{1-9-1}) and (\ref{e-8}) into (\ref{e-2}), we obtain $%
h_{111},$ $h_{110}$ ($=$ $h_{101})$ and $h_{100}$ in the torsion free case
as follows:%
\begin{eqnarray}
h_{111} &=&e_{1}(H)-2\alpha H,  \label{e-9} \\
h_{110} &=&(T+\alpha e_{2})(H)-3\alpha e_{1}(\alpha )-3\alpha ^{3}  \notag \\
&&-\frac{3}{2}\alpha W+\frac{1}{2}\func{Re}W^{,1},  \notag \\
h_{100} &=&(T+\alpha e_{2})(e_{1}(\alpha )+\frac{1}{2}\alpha ^{2}+\frac{1}{4}%
W)+\alpha He_{1}(\alpha )  \notag \\
&&+\alpha ^{3}H+\frac{1}{2}\alpha HW-\frac{1}{6}H\func{Re}W^{,1}.  \notag
\end{eqnarray}

\noindent We then compute the quantity $|H_{cr}|\mathfrak{f}$ in (\ref{e-3})
by (\ref{1-8}) and (\ref{e-9}) to obtain%
\begin{eqnarray}
|H_{cr}|\mathfrak{f} &\mathfrak{=}&h_{10}h_{111}+\frac{1}{3}%
h_{11}^{2}h_{111}+h_{11}h_{110}+\frac{3}{2}h_{100}  \label{e-10} \\
&=&e_{1}(H)(e_{1}(\alpha )+\frac{1}{2}\alpha ^{2}+\frac{1}{3}H^{2}+\frac{1}{4%
}W)  \notag \\
&&+H(T+\alpha e_{2})(H)+\frac{3}{2}(T+\alpha e_{2})(e_{1}(\alpha )+\frac{1}{2%
}\alpha ^{2})  \notag \\
&&-\frac{7}{2}\alpha He_{1}(\alpha )-\frac{5}{2}\alpha ^{3}H-\frac{2}{3}%
\alpha H^{3}-\frac{5}{4}\alpha HW  \notag
\end{eqnarray}

\noindent in the case of vanishing torsion ($A_{\bar{1}}^{1}$ $=$ $0$) and
constant Webster curvature ($W$ $=$ constant). Under the same assumption: $%
A_{\bar{1}}^{1}$ $=$ $0$ and $W$ $=$ constant, we compute$\neq $%
\begin{eqnarray}
&&9h_{00}+6h_{11}h_{10}+\frac{2}{3}h_{11}^{3}  \label{e-11} \\
&=&9(T+\alpha e_{2})(\alpha )+6H(e_{1}(\alpha )+\frac{1}{2}\alpha ^{2}+\frac{%
1}{4}W)+\frac{2}{3}H^{3}.  \notag
\end{eqnarray}

\noindent by (\ref{1-8}) and (\ref{1-12}). So in the situation of $A_{\bar{1}%
}^{1}$ $=$ $0$ and $W$ $=$ constant, in view of (\ref{e-1-1}), (\ref{e-10})
and (\ref{e-11}), we have an expression for (\ref{e-4}), the Euler-Lagrange
equation for $E_{1}$, in terms of $\alpha ,$ $H$ and their derivatives in
tangent directions $e_{1}$ and $T+\alpha e_{2}.$ We summarize what we have
obtained as a theorem.

\bigskip

\begin{theorem}
\label{Ee1} (\cite{Ch}) Let $\Sigma $ be a ($C^{\infty }$ smooth) surface in
a 3-dimensional pseudohermitian manifold $(M,J,\theta )$ with vanishing
torsion and constant Webster curvature$.$ Suppose $\Sigma $ is nonsingular
and $H_{cr}$ $\neq $ $0$ on $\Sigma $. Then $\Sigma $ satisfies the
Euler-Lagrange equation for the energy functional $E_{1}$ if and only if%
\begin{eqnarray}  \label{thesameEpsilon1}
0 &=&\frac{1}{2}e_{1}(|H_{cr}|^{1/2}\mathfrak{f)+}\frac{3}{2}%
|H_{cr}|^{1/2}\alpha \mathfrak{f} \\
&&+\frac{1}{2}sign(H_{cr})|H_{cr}|^{1/2}\{9h_{00}+6h_{11}h_{10}+\frac{2}{3}%
h_{11}^{3}\}  \notag
\end{eqnarray}

\noindent with $H_{cr},$ $\mathfrak{f}$ and $9h_{00}+6h_{11}h_{10}+\frac{2}{3%
}h_{11}^{3}$ expressed in (\ref{e-1-1}), (\ref{e-10}) and (\ref{e-11}),
resp..
\end{theorem}

\bigskip

\subsection{First variation of $E_2$}

Next we will deduce the Euler-Lagrange equation for $E_{2}$ in the case that 
$A_{\bar{1}}^{1}$ $=$ $0$ and $W=\text{constant}$. To deal with the term $%
(T+\alpha e_{2})(\alpha ),$ we observe that there holds%
\begin{eqnarray}
d(\alpha e^{1}) &=&d\alpha \wedge e^{1}+\alpha de^{1}  \label{e-12} \\
&=&[(T\alpha )\theta +e_{2}(\alpha )e^{2}]\wedge e^{1}-\alpha e^{2}\wedge
\omega  \notag \\
&=&[(T+\alpha e_{2})(\alpha )-\alpha ^{2}H]\theta \wedge e^{1}  \notag
\end{eqnarray}%
\noindent on $\Sigma .$ Here we have used (\ref{1-7-2}) and (\ref{1-7-1}).
Recall from (\ref{1-13}), we compute%
\begin{eqnarray}
&&dA_{2}  \label{e-12-0} \\
&=&\{(T+\alpha e_{2})(\alpha )+\frac{2}{3}[e_{1}(\alpha )+\frac{1}{2}\alpha
^{2}+\frac{1}{4}W]H+\frac{2}{27}H^{3}\}\theta \wedge e^{1}  \notag \\
&=&[\frac{2}{3}e_{1}(\alpha )+\frac{4}{3}\alpha ^{2}+\frac{1}{6}W+\frac{2}{27%
}H^{2}]H\theta \wedge e^{1}+d(\alpha e^{1}).  \notag
\end{eqnarray}

\noindent So we are reduced to computing the Euler-Lagrange equation of 
\begin{equation*}
\int_{\Sigma }[\frac{2}{3}e_{1}(\alpha )+\frac{4}{3}\alpha ^{2}+\frac{1}{6}W+%
\frac{2}{27}H^{2}]H\text{ }\theta \wedge e^{1}.
\end{equation*}

Let 
\begin{equation*}
\Sigma _{t}=F_{t}(\Sigma )
\end{equation*}%
be a family of immersions such that 
\begin{equation}
\frac{d}{dt}F_{t}=X=fe_{2}+gT.  \label{fg}
\end{equation}%
\noindent Here we let $e_{1}$ be the unit vector in $T\Sigma _{t}\cap \xi $
and $e_{2}=Je_{1}$. We assume $f$ and $g$ are supported in a domain of $%
\Sigma$ away from the singular set of $\Sigma$.

\bigskip

\begin{lemma}
Let $h:=f-\alpha g$ and $V:=$ $T+\alpha e_{2}.$ Then we have%
\begin{equation}
\omega (e_{2})=h^{-1}e_{1}(h)+2\alpha .  \label{omegae2}
\end{equation}%
Under the torsion free condition, we have%
\begin{equation}
\omega (T)=e_{1}(\alpha )-\alpha h^{-1}e_{1}(h),  \label{omegaT}
\end{equation}%
\begin{equation}
e_{2}(\alpha )=h^{-1}V(h).  \label{e2alpha}
\end{equation}%
\begin{equation}
e_{2}(H)=2W+4e_{1}(\alpha )+H^{2}+4\alpha ^{2}+h^{-1}e_{1}e_{1}(h)+2\alpha
h^{-1}e_{1}(h).  \label{e2H}
\end{equation}
\end{lemma}

\bigskip

%TCIMACRO{\TeXButton{Proof}{\proof} }%
%BeginExpansion
\proof
%EndExpansion
Suppose the surfaces $F_{t}(\Sigma )$ are the level sets of a defining
function $t$ such that 
\begin{equation*}
\frac{d}{dt}F_{t}=fe_{2}+gT,\quad t(F_{t}(\Sigma))=t.
\end{equation*}%
We observe that 
\begin{eqnarray*}
(fe_{2}+gT)t &=&1, \\
(T+\alpha e_{2})t &=&0.
\end{eqnarray*}
\noindent and hence we have%
\begin{eqnarray*}
&& T(t) =-\alpha e_{2}(t), \\
&& (f-\alpha g)e_{2}(t)=1.
\end{eqnarray*}%
\noindent Then for $f-\alpha g\neq 0,$ we have%
\begin{equation*}
e_{2}(t)=(f-\alpha g)^{-1}=h^{-1}.
\end{equation*}%
\noindent Note that $e_{1}(t)\equiv 0$ and $T=V-\alpha e_{2}$, so we have 
\begin{eqnarray*}
\lbrack e_{1},e_{2}](t)
&=&e_{1}e_{2}(t)-e_{2}e_{1}(t)=e_{1}e_{2}(t)=e_{1}(h^{-1})=-h^{-2}e_{1}(h) \\
&=&-\omega (e_{1})e_{1}(t)-\omega (e_{2})e_{2}(t)-2T(t) \\
&=&-\omega (e_{2})h^{-1}+2\alpha h^{-1},
\end{eqnarray*}%
\noindent which implies 
\begin{equation}
\omega (e_{2})=h^{-1}e_{1}(h)+2\alpha .  \label{e-14}
\end{equation}

Applying the first formula of (\ref{e-13}) to $t$, we find 
\begin{equation}
\omega (T)=e_{1}(\alpha )-\alpha h^{-1}e_{1}(h).  \label{e-15}
\end{equation}

By (\ref{e-13}) and%
\begin{equation*}
e_{2}(t)=(f-\alpha g)^{-1}=h^{-1},\quad e_{1}(t)=(T+\alpha e_{2})(t)=0,
\end{equation*}%
\noindent we have 
\begin{eqnarray*}
0 &=&e_{2}T(t)-Te_{2}(t) \\
&=&-e_{2}(\alpha h^{-1})-T(h^{-1}) \\
&=&h^{-2}[-he_{2}(\alpha )+\alpha e_{2}(h)+T(h)],
\end{eqnarray*}%
\noindent and hence 
\begin{equation*}
e_{2}(\alpha )=h^{-1}[\alpha e_{2}(h)+T(h)]=h^{-1}(T+\alpha
e_{2})(h)=h^{-1}V(h).
\end{equation*}

Substituting (\ref{1-7-1})$,$ (\ref{e-14}) and (\ref{e-15}) into%
\begin{equation}
-2W=e_{1}\omega (e_{2})-e_{2}\omega (e_{1})+\omega (e_{2})^{2}+\omega
(e_{1})^{2}+2\omega (T)  \label{e-15-1}
\end{equation}%
\noindent (obtained by applying the third formula of (\ref{e-13}) to $e_{1},$
$e_{2}$), we obtain%
\begin{equation*}
e_{2}(H)=2W+4e_{1}(\alpha )+H^{2}+4\alpha ^{2}+h^{-1}e_{1}e_{1}(h)+2\alpha
h^{-1}e_{1}(h).
\end{equation*}

%TCIMACRO{\TeXButton{End Proof}{\endproof}}%
%BeginExpansion
\endproof%
%EndExpansion

\bigskip

\begin{lemma}
Suppose the torsion vanishes, i.e., $A_{\bar{1}}^{1}$ $=$ $0$. Then on $%
\Sigma $ we have%
\begin{equation}
e_{1}e_{1}(\alpha )=-6\alpha e_{1}(\alpha )+V(H)-\alpha H^{2}-4\alpha
^{3}-2W\alpha .  \label{odealpha}
\end{equation}
\end{lemma}

\bigskip

%TCIMACRO{\TeXButton{Proof}{\proof} }%
%BeginExpansion
\proof
%EndExpansion
(\ref{odealpha}) had been proved in \cite{CHMYCodazzi}. For the convenience
of readers, we include a direct proof here. Note that (\ref{odealpha}) is
defined on $\Sigma$. To prove it, we extend the frame $e_{2},$ $e_{1}$ $=$ $%
-Je_{2}$ at points of $\Sigma $ to a neighborhood $U$ of $\Sigma $ (still
denoted by the same notation) so that $e_{2}$ $\in $ $\xi $ and 
\begin{equation*}
\nabla _{e_{2}}e_{2}=0,\quad e_{1}:=-Je_{2}
\end{equation*}%
\noindent in $U.$ Hence we have%
\begin{equation}
\omega (e_{2})=0.  \label{e-16}
\end{equation}%
\noindent Note that (\ref{e-16}) does not hold for $e_{2},$ $e_{1}$ defined
on $F_{t}(\Sigma )$ canonically. By (\ref{1-8-1}) we obtain 
\begin{eqnarray}
\omega (T) &=&e_{1}(\alpha )+2\alpha ^{2}\text{ and hence}  \label{e-17} \\
e_{1}e_{1}(\alpha ) &=&e_{1}(\omega (T))-4\alpha e_{1}(\alpha )  \notag
\end{eqnarray}%
\noindent on $\Sigma $. From (\ref{e-13}) and (\ref{e-16}), we obtain 
\begin{equation*}
0=d\omega (e_{1},T)=e_{1}(\omega (T))-T(\omega (e_{1})).
\end{equation*}%
\noindent It follows that 
\begin{eqnarray}
e_{1}(\omega (T)) &=&T(\omega (e_{1}))=(V-\alpha e_{2})(\omega (e_{1}))
\label{e-18} \\
&=&V(H)-\alpha (2W+H^{2}+2\omega (T)).  \notag
\end{eqnarray}

\noindent on $\Sigma .$ Here we have used%
\begin{equation*}
e_{2}(\omega (e_{1}))=2W+\omega (e_{1})^{2}+2\omega (T)
\end{equation*}%
\noindent by (\ref{e-15-1}) and (\ref{e-16}). Therefore from (\ref{e-17})
and (\ref{e-18}), we obtain (\ref{odealpha}).

%TCIMACRO{\TeXButton{End Proof}{\endproof}}%
%BeginExpansion
\endproof%
%EndExpansion

\bigskip

\begin{lemma}
\label{evolutionequations} Suppose the torsion vanishes, i.e., $A_{\bar{1}%
}^{1}$ $=$ $0$. Then we have%
\begin{equation}
\frac{d}{dt}[F_{t}^{\ast }(\theta \wedge e^{1})]=[-fH+V(g)]\theta \wedge
e^{1},  \label{e-18-1}
\end{equation}%
\begin{equation}
\frac{dH}{dt}=e_{1}e_{1}(h)+2\alpha e_{1}(h)+4h(e_{1}(\alpha )+\alpha ^{2}+%
\frac{1}{4}H^{2}+\frac{1}{2}W)+gV(H),  \label{e-18-2}
\end{equation}%
\begin{equation}
\frac{d\alpha }{dt}=V(h)+gV(\alpha ),  \label{e-18-3}
\end{equation}%
\begin{equation}
\frac{d}{dt}e_{1}(\alpha )=e_{1}V(h)+ge_{1}V(\alpha )+2fV(\alpha
)+fHe_{1}(\alpha ).  \label{e-18-4}
\end{equation}
\end{lemma}

\bigskip

%TCIMACRO{\TeXButton{Proof}{\proof} }%
%BeginExpansion
\proof
%EndExpansion
By (\ref{de}) in the Appendix we have%
\begin{eqnarray}
de^{1} &=&-e^{2}\wedge \omega  \label{e-19} \\
de^{2} &=&e^{1}\wedge \omega  \notag
\end{eqnarray}%
\noindent in the torsion free case. From (\ref{1-(-1)}), (\ref{e-19}) and (%
\ref{1-7-1})$,$ we have 
\begin{equation}
d(\theta \wedge e^{1})=-\theta \wedge (-e^{2}\wedge \omega )=-H\theta \wedge
e^{1}\wedge e^{2}.  \label{e-20}
\end{equation}%
\noindent We then compute 
\begin{eqnarray}
\frac{d}{dt}[F_{t}^{\ast }(\theta \wedge e^{1})] &=&(F_{t}^{\ast
})L_{fe_{2}+gT}(\theta \wedge e^{1})\text{ \ \ (pullback omitted)}
\label{e-20-1} \\
&=&i_{fe_{2}+gT}[-H\theta \wedge e^{1}\wedge e^{2}]+di_{fe_{2}+gT}(\theta
\wedge e^{1})  \notag \\
&=&-H(f\theta \wedge e^{1}+ge^{1}\wedge e^{2})+d(ge^{1})  \notag \\
&=&-(f-\alpha g)H\theta \wedge e^{1}+[e_{2}(g)\alpha +T(g)-\alpha gH]\theta
\wedge e^{1}  \notag \\
&=&[-fH+T(g)+\alpha e_{2}(g)]\theta \wedge e^{1}  \notag \\
&=&[-fH+V(g)]\theta \wedge e^{1}.  \notag
\end{eqnarray}

\noindent by (\ref{e-20}), (\ref{e-19}), (\ref{1-7-2}) and (\ref{1-7-1})$.$

We now compute the first variation of $H$ and $\alpha $. Observe that 
\begin{equation*}
\frac{d}{dt}(F_{t}^{\ast }\omega )=L_{fe_{2}+gT}\omega =d(f\omega
(e_{2})+g\omega (T))+2Wfe^{1}.
\end{equation*}%
\noindent It follows that%
\begin{eqnarray}
\frac{dH}{dt} &=&\frac{d}{dt}[\omega (e_{1})]  \label{e-21} \\
&=&(L_{fe_{2}+gT}\omega )(e_{1})+\omega ([fe_{2}+gT,e_{1}])  \notag \\
&=&(d(f\omega (e_{2})+g\omega (T)))(e_{1})+2Wf+\omega (-e_{1}(f)e_{2}  \notag
\\
&&+f[e_{2},e_{1}]-e_{1}(g)T+g[T,e_{1}])  \notag \\
&=&f[2W+e_{1}(\omega (e_{2}))+\omega (e_{1})^{2}+\omega (e_{2})^{2}+2\omega
(T)]  \notag \\
&&+g[e_{1}(\omega (T))+\omega (e_{2})\omega (T)].  \notag
\end{eqnarray}%
\noindent Substituting (\ref{omegae2}) and (\ref{omegaT}) into (\ref{e-21})
gives 
\begin{eqnarray*}
\frac{dH}{dt} &=&f[h^{-1}e_{1}e_{1}(h)+2\alpha h^{-1}e_{1}(h)+4e_{1}(\alpha
)+4\alpha ^{2}+H^{2}+2W] \\
&&+g[e_{1}e_{1}(\alpha )-\alpha h^{-1}e_{1}e_{1}(h)+2\alpha e_{1}(\alpha
)-2\alpha ^{2}h^{-1}e_{1}(h)] \\
&=&e_{1}e_{1}(h)+2\alpha e_{1}(h)+f[4e_{1}(\alpha )+4\alpha
^{2}+H^{2}+2W]+g[e_{1}e_{1}(\alpha )+2\alpha e_{1}(\alpha )].
\end{eqnarray*}%
\noindent It follows from (\ref{odealpha}) that 
\begin{equation*}
\frac{dH}{dt}=e_{1}e_{1}(h)+2\alpha e_{1}(h)+4h[e_{1}(\alpha )+\alpha ^{2}+%
\frac{1}{4}H^{2}+\frac{1}{2}W]+gV(H).
\end{equation*}

From (\ref{e-19}) and $e^{2}\wedge e^{1}$ $=$ $\alpha \theta \wedge e^{1}$
on $F_{t}(\Sigma )$ by (\ref{1-7-2}), we compute%
\begin{eqnarray*}
\frac{d}{dt}(F_{t}^{\ast }e^{1})
&=&L_{fe_{2}+gT}e^{1}=i_{fe_{2}+gT}(-e^{2}\wedge \omega )=-f\omega
(e_{1})e^{1}-f\omega (T)\theta +g\omega (T)e^{2}, \\
\frac{d}{dt}(F_{t}^{\ast }e^{2})
&=&L_{fe_{2}+gT}e^{2}=i_{fe_{2}+gT}(e^{1}\wedge \omega )+df=-(f\omega
(e_{2})+g\omega (T))e^{1}+df
\end{eqnarray*}%
\noindent and hence 
\begin{eqnarray*}
\frac{d}{dt}[F_{t}^{\ast }(e^{2}\wedge e^{1})] &=&(L_{fe_{2}+gT}e^{2})\wedge
e^{1}+e^{2}\wedge L_{fe_{2}+gT}e^{1} \\
&=&df\wedge e^{1}+e^{2}\wedge \lbrack -f\omega (e_{1})e^{1}] \\
&=&[\alpha e_{2}(f)+T(f)-\alpha Hf]\theta \wedge e^{1}.
\end{eqnarray*}%
\noindent Here we have used $e^{2}\wedge e^{1}$ $=$ $\alpha \theta \wedge
e^{1},$ $e^{2}\wedge \theta $ $=$ $0$ and $\omega (e_{1})$ $=$ $H$ on $%
F_{t}(\Sigma )$ by (\ref{1-7-2}) and (\ref{1-7-1})$.$ On the other hand, we
have 
\begin{eqnarray*}
\frac{d}{dt}[F_{t}^{\ast }(e^{2}\wedge e^{1})] &=&\frac{d}{dt}[F_{t}^{\ast
}(\alpha \theta \wedge e^{1})] \\
&=&\frac{d\alpha }{dt}\theta \wedge e^{1}+\alpha \lbrack -fH+V(g)]\theta
\wedge e^{1}
\end{eqnarray*}%
\noindent by (\ref{e-20-1}). So we get 
\begin{eqnarray*}
\frac{d\alpha }{dt} &=&\alpha e_{2}(f)+T(f)-\alpha V(g) \\
&=&V(f)-\alpha V(g)\text{ \ \ (}V:=T+\alpha e_{2}) \\
&=&V(h)+gV(\alpha )
\end{eqnarray*}

\noindent since $h$ $:=$ $f-\alpha g.$ We now compute%
\begin{eqnarray}
&&\frac{d}{dt}e_{1}(\alpha )  \label{e-22} \\
&=&\frac{d}{dt}[(d\alpha )(e_{1})]  \notag \\
&=&(L_{fe_{2}+gT}d\alpha )(e_{1})+d\alpha ([fe_{2}+gT,e_{1}])  \notag \\
&=&e_{1}(fe_{2}+gT)(\alpha )+[fe_{2}+gT,e_{1}](\alpha )  \notag \\
&=&fe_{1}e_{2}(\alpha )+ge_{1}T(\alpha )+f[e_{2},e_{1}](\alpha
)+g[T,e_{1}](\alpha )  \notag \\
&=&fe_{1}e_{2}(\alpha )+ge_{1}T(\alpha )  \notag \\
&&+f[He_{1}(\alpha )+\omega (e_{2})e_{2}(\alpha )+2T(\alpha )]+g\omega
(T)e_{2}(\alpha ).  \notag
\end{eqnarray}%
\noindent For the last equality of (\ref{e-22}), we have used (\ref{e-12-1})
and (\ref{e-13}). Substituting (\ref{omegae2}), (\ref{omegaT}) and (\ref%
{e2alpha}) into (\ref{e-22}) gives 
\begin{eqnarray*}
\frac{d}{dt}e_{1}(\alpha ) &=&fe_{1}e_{2}(\alpha )+ge_{1}V(\alpha
)-ge_{1}(\alpha )e_{2}(\alpha )-\alpha ge_{1}e_{2}(\alpha ) \\
&&+fHe_{1}(\alpha )+f(h^{-1}e_{1}(h)+2\alpha )e_{2}(\alpha )+2fV(\alpha
)-2\alpha fe_{2}(\alpha ) \\
&&+ge_{1}(\alpha )e_{2}(\alpha )-\alpha gh^{-1}e_{1}(h)e_{2}(\alpha ) \\
&=&fe_{1}e_{2}(\alpha )+ge_{1}V(\alpha )-\alpha ge_{1}e_{2}(\alpha ) \\
&&+fHe_{1}(\alpha )+fh^{-1}e_{1}(h)e_{2}(\alpha )+2fV(\alpha )-\alpha
gh^{-1}e_{1}(h)e_{2}(\alpha ) \\
&=&he_{1}e_{2}(\alpha )+ge_{1}V(\alpha )+fHe_{1}(\alpha
)+e_{1}(h)e_{2}(\alpha )+2fV(\alpha ) \\
&=&e_{1}V(h)+ge_{1}V(\alpha )+2fV(\alpha )+fHe_{1}(\alpha ).
\end{eqnarray*}

%TCIMACRO{\TeXButton{End Proof}{\endproof}}%
%BeginExpansion
\endproof%
%EndExpansion

\bigskip

\begin{lemma}
\label{integrationbypart} Suppose either $f_{1}$ or $f_{2}$ has compact
support in the nonsingular domain of $\Sigma .$ Then we have 
\begin{equation}
\int_{\Sigma }f_{1}e_{1}(f_{2})\theta \wedge e^{1}=\int_{\Sigma
}-[e_{1}(f_{1})+2\alpha f_{1}]f_{2}\theta \wedge e^{1},  \label{e-23}
\end{equation}%
\begin{equation}
\int_{\Sigma }f_{1}V(f_{2})\theta \wedge e^{1}=-\int_{\Sigma
}[V(f_{1})-\alpha Hf_{1}]f_{2}\theta \wedge e^{1}.  \label{e-24}
\end{equation}
\end{lemma}

\bigskip

%TCIMACRO{\TeXButton{Proof}{\proof} }%
%BeginExpansion
\proof
%EndExpansion
Note that $df=e_1(f)e^1+\frac{1}{1+\alpha^2}V(f)(\theta+\alpha e^2)$. We
have 
\begin{eqnarray*}
0 &=&\int_{\Sigma }d(f_{1}f_{2}\theta )=\int_{\Sigma
}e_{1}(f_{1}f_{2})e^{1}\wedge \theta +2f_{1}f_{2}e^{1}\wedge e^{2}
\end{eqnarray*}%
\noindent since $e^{2}\wedge \theta $ $=$ $0$ on $\Sigma $ and (\ref{1-(-1)}%
). By (\ref{1-7-2}) we have 
\begin{equation*}
\int_{\Sigma }f_{1}e_{1}(f_{2})\theta \wedge e^{1}=\int_{\Sigma
}-[e_{1}(f_{1})+2\alpha f_{1}]f_{2}\theta \wedge e^{1}.
\end{equation*}

On the other hand, we also have 
\begin{eqnarray*}
0 &=&\int_{\Sigma }d(f_{1}f_{2}e^{1})=\int_{\Sigma }[V(f_{1}f_{2})\theta
\wedge e^{1}+f_{1}f_{2}(-e^{2}\wedge \omega )] \\
&=&\int_{\Sigma }[V(f_{1}f_{2})-\alpha Hf_{1}f_{2}]\theta \wedge e^{1},
\end{eqnarray*}%
\noindent from (\ref{1-7-2}), (\ref{e-19}) and (\ref{1-7-1}). It follows
that 
\begin{equation*}
\int_{\Sigma }f_{1}V(f_{2})\theta \wedge e^{1}=-\int_{\Sigma
}[V(f_{1})-\alpha Hf_{1}]f_{2}\theta \wedge e^{1}.
\end{equation*}

%TCIMACRO{\TeXButton{End Proof}{\endproof}}%
%BeginExpansion
\endproof%
%EndExpansion

\bigskip

\begin{theorem}
\label{firstvariationE2} Assume $(M,J,\theta )$ has vanishing torsion and
constant Webster scalar curvature. Let $F_t(\Sigma)$ be given by (\ref{fg}).
We have 
\begin{equation}  \label{firstvariationofE2}
\frac{d}{dt}\int_{F_{t}(\Sigma )}dA_{2}=\int_{\Sigma }\mathcal{E}_{2}h\theta
\wedge e^{1}.
\end{equation}%
where 
\begin{eqnarray}
\mathcal{E}_{2} &=&\frac{4}{9}[He_{1}e_{1}(H)+3e_{1}V(H)+e_{1}(H)^{2}+\frac{1%
}{3}H^{4}  \notag \\
&&+3e_{1}(\alpha )^{2}+12\alpha ^{2}e_{1}(\alpha )+12\alpha ^{4}  \notag \\
&&-\alpha He_{1}(H)+2H^{2}e_{1}(\alpha )+5\alpha ^{2}H^{2}  \label{mathcalE}
\\
&&+\frac{3}{2}W(e_{1}(\alpha )+\frac{2}{3}H^{2}+5\alpha ^{2}+\frac{1}{2}W)].
\notag
\end{eqnarray}
\end{theorem}

%TCIMACRO{\TeXButton{Proof}{\proof} }%
%BeginExpansion
\proof
%EndExpansion
From (\ref{e-12-0}), for $f,$ $g,$ and hence $h$ $=$ $f-\alpha g$ having
compact support in nonsingular domain of $\Sigma ,$ we compute 
\begin{eqnarray*}
&&\frac{d}{dt}\int_{F_{t}(\Sigma )}dA_{2} \\
&=&\frac{d}{dt}\int_{\Sigma }F_{t}^{\ast }([\frac{2}{3}e_{1}(\alpha )+\frac{4%
}{3}\alpha ^{2}+\frac{1}{6}W+\frac{2}{27}H^{2}]H)F_{t}^{\ast }(\theta \wedge
e^{1}) \\
&=&\frac{2}{3}\int_{\Sigma }[(e_{1}(\alpha )+2\alpha ^{2}+\frac{1}{3}H^{2}+%
\frac{1}{4}W)\frac{dH}{dt}+H\frac{d}{dt}(e_{1}(\alpha ))+4\alpha H\frac{%
d\alpha }{dt}]\theta \wedge e^{1} \\
&&+\int_{\Sigma }[\frac{2}{3}e_{1}(\alpha )+\frac{4}{3}\alpha ^{2}+\frac{1}{6%
}W+\frac{2}{27}H^{2}]H\frac{d}{dt}[F_{t}^{\ast }(\theta \wedge e^{1})] \\
&:&=I+II,
\end{eqnarray*}%
\noindent where 
\begin{eqnarray*}
I &:&=\frac{2}{3}\int_{\Sigma }[(e_{1}(\alpha )+2\alpha ^{2}+\frac{1}{3}%
H^{2}+\frac{1}{4}W)\frac{dH}{dt}+H\frac{d}{dt}(e_{1}(\alpha ))+4\alpha H%
\frac{d\alpha }{dt}]\theta \wedge e^{1} \\
&=&\int_{\Sigma }\frac{2}{3}(e_{1}(\alpha )+2\alpha ^{2}+\frac{1}{3}H^{2}+%
\frac{1}{4}W) \\
&&[e_{1}e_{1}(h)+2\alpha e_{1}(h)+4h(e_{1}(\alpha )+\alpha ^{2}+\frac{1}{4}%
H^{2}+\frac{1}{2}W)+gV(H)]\theta \wedge e^{1} \\
&&+\int_{\Sigma }\frac{2}{3}H[e_{1}V(h)+ge_{1}V(\alpha )+2fV(\alpha
)+fHe_{1}(\alpha )]\theta \wedge e^{1} \\
&&+\int_{\Sigma }\frac{8}{3}\alpha H[V(h)+gV(\alpha )]\theta \wedge e^{1}
\end{eqnarray*}

\noindent by (\ref{e-18-2}), (\ref{e-18-4}) and (\ref{e-18-3}), and 
\begin{eqnarray*}
II&:&=\int_{\Sigma }[\frac{2}{3}e_{1}(\alpha )+\frac{4}{3}\alpha ^{2}+\frac{1%
}{6}W+\frac{2}{27}H^{2}]H\frac{d}{dt}[F_{t}^{\ast }(\theta \wedge e^{1})] \\
&=&\int_{\Sigma }(\frac{2}{3}He_{1}(\alpha )+\frac{4}{3}\alpha ^{2}H+\frac{2%
}{27}H^{3}+\frac{1}{6}WH)[-fH+V(g)]\theta \wedge e^{1}
\end{eqnarray*}%
\noindent by (\ref{e-18-1}). Let 
\begin{eqnarray*}
H_{1}:&=&\int_{\Sigma }\frac{2}{3}(e_{1}(\alpha )+2\alpha ^{2}+\frac{1}{3}%
H^{2}+\frac{1}{4}W)[e_{1}e_{1}(h)+2\alpha e_{1}(h)]\theta \wedge e^{1}, \\
H_{2}:&=&\int_{\Sigma }\frac{8}{3}(e_{1}(\alpha )+2\alpha ^{2}+\frac{1}{3}%
H^{2}+\frac{1}{4}W)(e_{1}(\alpha )+\alpha ^{2}+\frac{1}{4}H^{2}+\frac{1}{2}%
W)h\theta \wedge e^{1}, \\
H_{3}:&=&\int_{\Sigma }(\frac{2}{3}He_{1}V(h)+\frac{8}{3}\alpha HV(h))\theta
\wedge e^{1}, \\
F :&=&\int_{\Sigma }[\frac{4}{3}HV(\alpha )-H(\frac{4}{3}\alpha ^{2}H+\frac{2%
}{27}H^{3}+\frac{1}{6}WH)]f\theta \wedge e^{1} \\
G_{1} :&=&\int_{\Sigma }\frac{2}{3}(e_{1}(\alpha )+2\alpha ^{2}+\frac{1}{3}%
H^{2}+\frac{1}{4}W)V(H)g\theta \wedge e^{1} \\
&&+\int_{\Sigma }[\frac{2}{3}He_{1}V(\alpha )+\frac{8}{3}\alpha HV(\alpha
)]g\theta \wedge e^{1}, \\
G_{2}:&=&\int_{\Sigma }(\frac{2}{3}He_{1}(\alpha )+\frac{4}{3}\alpha ^{2}H+%
\frac{2}{27}H^{3}+\frac{1}{6}WH)V(g)\theta \wedge e^{1}.
\end{eqnarray*}%
\noindent So we have%
\begin{eqnarray}
&&\frac{d}{dt}\int_{F_{t}(\Sigma )}dA_{2}  \label{e-24-0} \\
&=&I+II=H_{1}+H_{2}+H_{3}+F+G_{1}+G_{2}.  \notag
\end{eqnarray}

Applying the formulas (\ref{e-23}) and (\ref{e-24}) of integration by parts
and the ODE (\ref{odealpha}) for $\alpha $, we compute 
\begin{eqnarray*}
H_{1} &:&=\int_{\Sigma }\frac{2}{3}(e_{1}(\alpha )+2\alpha ^{2}+\frac{1}{3}%
H^{2}+\frac{1}{4}W)[e_{1}e_{1}(h)+2\alpha e_{1}(h)]\theta \wedge e^{1} \\
&=&-\frac{2}{3}\int_{\Sigma }e_{1}(e_{1}(\alpha )+2\alpha ^{2}+\frac{1}{3}%
H^{2}+\frac{1}{4}W)e_{1}(h)\theta \wedge e^{1} \\
&=&-\frac{2}{3}\int_{\Sigma }[\frac{2}{3}He_{1}(H)-2\alpha e_{1}(\alpha
)-4\alpha ^{3}-\alpha H^{2}+V(H)-2W\alpha ]e_{1}(h)\theta \wedge e^{1} \\
&=&\frac{2}{3}\int_{\Sigma }\{e_{1}[\frac{2}{3}He_{1}(H)-2\alpha
e_{1}(\alpha )-4\alpha ^{3}-\alpha H^{2}+V(H)-2W\alpha] \\
&&+2\alpha \lbrack \frac{2}{3}He_{1}(H)-2\alpha e_{1}(\alpha )-4\alpha
^{3}-\alpha H^{2}+V(H)-2W\alpha ]\}h\theta \wedge e^{1} \\
&=&\frac{4}{9}\int_{\Sigma }[He_{1}e_{1}(H)+e_{1}(H)^{2}-3e_{1}(\alpha
)^{2}-\alpha He_{1}(H)-\frac{3}{2}H^{2}e_{1}(\alpha ) \\
&&-6\alpha ^{2}e_{1}(\alpha )+\frac{3}{2}e_{1}V(H) -3We_{1}(\alpha )]h\theta
\wedge e^{1}
\end{eqnarray*}

\noindent and 
\begin{eqnarray*}
H_{3} &:&=\int_{\Sigma }(\frac{2}{3}He_{1}V(h)+\frac{8}{3}\alpha
HV(h))\theta \wedge e^{1} \\
&=&-\int_{\Sigma }\frac{2}{3}[e_{1}(H)+2\alpha H]V(h)\theta \wedge
e^{1}+\int_{\Sigma }\frac{8}{3}\alpha HV(h)\theta \wedge e^{1} \\
&=&-\int_{\Sigma }\frac{2}{3}[e_{1}(H)-2\alpha H]V(h)\theta \wedge e^{1} \\
&=&\int_{\Sigma }\frac{2}{3}(V[e_{1}(H)-2\alpha H]-\alpha H[e_{1}(H)-2\alpha
H])h\theta \wedge e^{1} \\
&=&\int_{\Sigma }[\frac{2}{3}e_{1}V(H)-\frac{4}{3}HV(\alpha )+\frac{4}{3}%
\alpha ^{2}H^{2}]h\theta \wedge e^{1}
\end{eqnarray*}%
\noindent where we have used 
\begin{eqnarray}
\lbrack e_{1},V] &=&[e_{1},T+\alpha e_{2}]=-\omega (T)e_{2}+e_{1}(\alpha
)e_{2}-\alpha \lbrack e_{2},e_{1}]  \label{e-24-1} \\
&=&\alpha h^{-1}e_{1}(h)e_{2}-\alpha He_{1}-\alpha
h^{-1}e_{1}(h)e_{2}-2\alpha ^{2}e_{2}-2\alpha T  \notag \\
&=&-\alpha He_{1}-2\alpha V.  \notag
\end{eqnarray}%
\noindent by (\ref{e-13}), (\ref{e-12-1}), (\ref{omegaT}) and (\ref{omegae2}%
). Therefore we get%
\begin{eqnarray}
&&H_{1}+H_{2}+H_{3}  \label{allH} \\
&=&\frac{4}{9}\int_{\Sigma }[He_{1}e_{1}(H)+3e_{1}V(H)+e_{1}(H)^{2}+\frac{1}{%
2}H^{4}]h\theta \wedge e^{1}  \notag \\
&&+\frac{4}{3}\int_{\Sigma }[e_{1}(\alpha )^{2}+4\alpha^{2}e_{1}(\alpha
)+4\alpha ^{4}]h\theta \wedge e^{1}  \notag \\
&&-\frac{4}{9}\int_{\Sigma }[\alpha He_{1}(H)+3HV(\alpha
)-2H^{2}e_{1}(\alpha )-8\alpha ^{2}H^{2}]h\theta \wedge e^{1}  \notag \\
&&+\frac{2}{3}\int_{\Sigma }W(e_{1}(\alpha )+\frac{11}{12}H^{2}+5\alpha^{2}+%
\frac{W}{2}) h\theta \wedge e^{1}.  \notag
\end{eqnarray}

On the other hand, from (\ref{e-24}) and (\ref{e-24-1}) we compute 
\begin{eqnarray*}
G_{2} &:&=\int_{\Sigma }(\frac{2}{3}He_{1}(\alpha )+\frac{4}{3}\alpha ^{2}H+%
\frac{2}{27}H^{3}+\frac{1}{6}WH)V(g)\theta \wedge e^{1} \\
&=&\int_{\Sigma }-V(\frac{2}{3}He_{1}(\alpha )+\frac{4}{3}\alpha ^{2}H+\frac{%
2}{27}H^{3}+\frac{1}{6}WH)g\theta \wedge e^{1} \\
&&+\int_{\Sigma }\alpha H(\frac{2}{3}He_{1}(\alpha )+\frac{4}{3}\alpha ^{2}H+%
\frac{2}{27}H^{3}+\frac{1}{6}WH)g\theta \wedge e^{1} \\
&=&\int_{\Sigma }-[\frac{2}{3}HVe_{1}(\alpha )+\frac{2}{3}e_{1}(\alpha )V(H)+%
\frac{4}{3}\alpha ^{2}V(H)+\frac{8}{3}\alpha HV(\alpha ) \\
&&+\frac{2}{9}H^{2}V(H)+\frac{1}{6}WV(H)]g\theta \wedge e^{1} \\
&&+\int_{\Sigma }\alpha H(\frac{2}{3}He_{1}(\alpha )+\frac{4}{3}\alpha ^{2}H+%
\frac{2}{27}H^{3}+\frac{1}{6}WH)g\theta \wedge e^{1} \\
&=&\int_{\Sigma }[-\frac{2}{3}He_{1}V(\alpha )-\frac{2}{3}e_{1}(\alpha )V(H)-%
\frac{4}{3}\alpha ^{2}V(H)-4\alpha HV(\alpha ) \\
&&-\frac{2}{9}H^{2}V(H)-\frac{1}{6}WV(H)+\frac{4}{3}\alpha ^{3}H^{2}+\frac{2%
}{27}\alpha H^{4}+\frac{1}{6}W\alpha H^{2}]g\theta \wedge e^{1}.
\end{eqnarray*}%
\noindent So we have 
\begin{equation*}
G_{1}+G_{2}=\int_{\Sigma }(-\frac{4}{3}HV(\alpha )+\frac{4}{3}\alpha
^{2}H^{2}+\frac{2}{27}H^{4}+\frac{1}{6}WH^{2})\alpha g\theta \wedge e^{1},
\end{equation*}%
\begin{equation}
F+G_{1}+G_{2}=\int_{\Sigma }[\frac{4}{3}HV(\alpha )-\frac{4}{3}\alpha
^{2}H^{2}-\frac{2}{27}H^{4}-\frac{1}{6}WH^{2}]h\theta \wedge e^{1}.
\label{FG}
\end{equation}%
(\ref{firstvariationofE2}) then follows from (\ref{e-24-0}), (\ref{allH})
and (\ref{FG}).

%TCIMACRO{\TeXButton{End Proof}{\endproof}}%
%BeginExpansion
\endproof%
%EndExpansion

\bigskip

\begin{theorem}
\label{Ee2} Let $\Sigma $ be a ($C^{\infty }$ smooth) surface in a
3-dimensional pseudohermitian manifold $(M,J,\theta )$ with vanishing
torsion and constant Webster scalar curvature. Then $\Sigma $ is critical,
subject to variations supported in any nonsingular domain, for the energy
functional $E_{2}$ if and only if
\end{theorem}

\begin{equation}
\mathcal{E}_{2}=0.  \label{EulerLagrange2}
\end{equation}

\bigskip

\subsection{First variation of $E_{1}$-An alternative approach}

For the reader's reference, we give an alternative deduction for the first
variation of $E_{1}$ in this subsection. Assume $A_{\overline{1}}^{1}=0$ and 
$W$ is constant.

\bigskip

\begin{theorem}
Assume $(M,J,\theta )$ has vanishing torsion and constant Webster scalar
curvature. Let $F_{t}(\Sigma )$ be given by (\ref{fg}). We have%
\begin{equation}
\frac{d}{dt}\int_{F_{t}(\Sigma )}dA_{1}=\int_{\Sigma }\mathcal{E}_{1}h\theta
\wedge e^{1},  \label{firstvariationofE1}
\end{equation}
where $\mathcal{E}_{1}$ is given by the right hand side of (\ref%
{thesameEpsilon1}).
\end{theorem}

\bigskip

%TCIMACRO{\TeXButton{Proof}{\proof} }%
%BeginExpansion
\proof
%EndExpansion
The proof is similar to that of Theorem \ref{firstvariationE2}. We only
sketch the main steps. We have 
\begin{eqnarray*}
\frac{d}{dt}E_{1}(\Sigma _{t}) &=&\frac{d}{dt}\int_{F_{t}(\Sigma
)}[(e_{1}(\alpha )+\frac{1}{2}\alpha ^{2}+\frac{1}{6}H^{2}+\frac{1}{4}%
W)^{2}]^{3/4}\theta \wedge e^{1} \\
&=&\int_{\Sigma }\frac{3}{2}(H_{cr}^{2})^{-\frac{1}{4}}H_{cr}\frac{d}{dt}%
[e_{1}(\alpha )+\frac{1}{2}\alpha ^{2}+\frac{1}{6}H^{2}+\frac{1}{4}W]\theta
\wedge e^{1} \\
&&+\int_{\Sigma }(H_{cr}^{2})^{\frac{3}{4}}\frac{d}{dt}F_{t}^{\ast }(\theta
\wedge e^{1}).
\end{eqnarray*}%
Using Lemma \ref{evolutionequations}, we get 
\begin{eqnarray*}
&&\frac{d}{dt}E_{1}(\Sigma _{t}) \\
&=&\int_{\Sigma }\frac{3}{2}(H_{cr}^{2})^{-\frac{1}{4}%
}H_{cr}[e_{1}V(h)+ge_{1}V(\alpha )+2fV(\alpha )+fHe_{1}(\alpha )]\theta
\wedge e^{1} \\
&&+\int_{\Sigma }\frac{3}{2}(H_{cr}^{2})^{-\frac{1}{4}}H_{cr}\alpha \lbrack
V(h)+gV(\alpha )]\theta \wedge e^{1}+\int_{\Sigma }\frac{1}{2}(H_{cr}^{2})^{-%
\frac{1}{4}}H_{cr}HV(H)g\theta \wedge e^{1} \\
&&+\int_{\Sigma }\frac{1}{2}(H_{cr}^{2})^{-\frac{1}{4}%
}H_{cr}H[e_{1}e_{1}(h)+2\alpha e_{1}(h)+4h(e_{1}(\alpha )+\alpha ^{2}+\frac{1%
}{4}H^{2}+\frac{1}{2}W)]\theta \wedge e^{1} \\
&&+\int_{\Sigma }(H_{cr}^{2})^{\frac{3}{4}}[-fH+V(g)]\theta \wedge e^{1}.
\end{eqnarray*}%
Let 
\begin{eqnarray*}
H_{1}:= &&\int_{\Sigma }\frac{3}{2}(H_{cr}^{2})^{-\frac{1}{4}%
}H_{cr}e_{1}V(h)\theta \wedge e^{1}+\int_{\Sigma }\frac{3}{2}(H_{cr}^{2})^{-%
\frac{1}{4}}H_{cr}\alpha V(h)\theta \wedge e^{1}, \\
H_{2}:= &&\int_{\Sigma }\frac{1}{2}(H_{cr}^{2})^{-\frac{1}{4}%
}H_{cr}H[e_{1}e_{1}(h)+2\alpha e_{1}(h)]\theta \wedge e^{1}, \\
H_{3}:= &&\int_{\Sigma }2(H_{cr}^{2})^{-\frac{1}{4}}H_{cr}H(e_{1}(\alpha
)+\alpha ^{2}+\frac{1}{4}H^{2}+\frac{1}{2}W)h\theta \wedge e^{1}, \\
F:= &&\int_{\Sigma }\frac{3}{2}(H_{cr}^{2})^{-\frac{1}{4}}H_{cr}[2V(\alpha
)+He_{1}(\alpha )]f\theta \wedge e^{1}+\int_{\Sigma }-H(H_{cr}^{2})^{\frac{3%
}{4}}f\theta \wedge e^{1}, \\
G_{1}:= &&\int_{\Sigma }\frac{3}{2}(H_{cr}^{2})^{-\frac{1}{4}%
}H_{cr}[e_{1}V(\alpha )+\alpha V(\alpha )+\frac{1}{3}HV(H)]g\theta \wedge
e^{1}, \\
G_{2}:= &&\int_{\Sigma }(H_{cr}^{2})^{\frac{3}{4}}V(g)\theta \wedge e^{1}.
\end{eqnarray*}%
Then 
\begin{equation*}
\frac{d}{dt}E_{1}(\Sigma _{t})=H_{1}+H_{2}+H_{3}+F+G_{1}+G_{2}.
\end{equation*}

Using Lemma \ref{integrationbypart}, we obtain 
\begin{eqnarray*}
H_{1} &=&\int_{\Sigma }\frac{3}{2}[V((H_{cr}^{2})^{-\frac{1}{4}}(\frac{1}{2}%
e_{1}(H_{cr})+\alpha H_{cr}))-\alpha H(H_{cr}^{2})^{-\frac{1}{4}}(\frac{1}{2}%
e_{1}(H_{cr})+\alpha H_{cr})]h\theta \wedge e^{1},
\end{eqnarray*}%
\begin{eqnarray*}
H_{2} &=&\int_{\Sigma }\frac{1}{2}(e_{1}e_{1}[(H_{cr}^{2})^{-\frac{1}{4}%
}H_{cr}H]+2\alpha e_{1}[(H_{cr}^{2})^{-\frac{1}{4}}H_{cr}H])h\theta \wedge
e^{1}
\end{eqnarray*}%
and 
\begin{eqnarray*}
G_{2} &=&\int_{\Sigma }[-\frac{3}{2}(H_{cr}^{2})^{-\frac{1}{4}%
}H_{cr}V(H_{cr})+\alpha H(H_{cr}^{2})^{\frac{3}{4}}]g\theta \wedge e^{1}.
\end{eqnarray*}%
Hence, by using (\ref{e-24-1}) we have 
\begin{eqnarray*}
G_{1}+G_{2} &=&\int_{\Sigma }\frac{3}{2}(H_{cr}^{2})^{-\frac{1}{4}%
}H_{cr}[2V(\alpha )+He_{1}(\alpha )](-\alpha g)\theta \wedge e^{1} \\
&&+\int_{\Sigma }(H_{cr}^{2})^{\frac{3}{4}}\alpha Hg\theta \wedge e^{1},
\end{eqnarray*}%
\begin{eqnarray*}
F+G_{1}+G_{2} &=&\int_{\Sigma }\frac{3}{2}(H_{cr}^{2})^{-\frac{1}{4}%
}H_{cr}[2V(\alpha )+He_{1}(\alpha )]h\theta \wedge e^{1} \\
&&-\int_{\Sigma }(H_{cr}^{2})^{\frac{3}{4}}Hh\theta \wedge e^{1}.
\end{eqnarray*}%
Therefore we obtain the first variation of $E_{1}:$ 
\begin{equation*}
\frac{d}{dt}E_{1}(\Sigma _{t})=\int_{\Sigma }\mathcal{E}_{1}h\theta \wedge
e^{1}
\end{equation*}%
\noindent where 
\begin{eqnarray*}
\mathcal{E}_{1} &=&\frac{3}{2}[V((H_{cr}^{2})^{-\frac{1}{4}}(\frac{1}{2}%
e_{1}(H_{cr})+\alpha H_{cr}))-\alpha H(H_{cr}^{2})^{-\frac{1}{4}}(\frac{1}{2}%
e_{1}(H_{cr})+\alpha H_{cr})]  \label{ScriptE1} \\
&&+\frac{1}{2}(e_{1}e_{1}[(H_{cr}^{2})^{-\frac{1}{4}}H_{cr}H]+2\alpha
e_{1}[(H_{cr}^{2})^{-\frac{1}{4}}H_{cr}H])  \notag \\
&&+2(H_{cr}^{2})^{-\frac{1}{4}}H_{cr}H(e_{1}(\alpha )+\alpha ^{2}+\frac{1}{4}%
H^{2}+\frac{1}{2}W)  \notag \\
&&+\frac{3}{2}(H_{cr}^{2})^{-\frac{1}{4}}H_{cr}[2V(\alpha )+He_{1}(\alpha
)]-(H_{cr}^{2})^{\frac{3}{4}}H  \notag
\end{eqnarray*}

\noindent assuming $H_{cr}$ never vanishes on the (nonsingular) support of $%
h.$ We compute 
\begin{eqnarray*}
\mathcal{E}_1&=&\frac{3}{2}V((H_{cr}^2)^{-\frac{1}{4}}(\frac{1}{2}%
e_1(H_{cr})+\alpha H_{cr}))+\frac{1}{2}e_1e_1[(H_{cr}^2)^{-\frac{1}{4}%
}H_{cr}H] +\alpha e_1[(H_{cr}^2)^{-\frac{1}{4}}H_{cr}H] \\
&&+|H_{cr}|^{-\frac{1}{2}}\{-\frac{3}{2}\alpha H(\frac{1}{2}%
e_1(H_{cr})+\alpha H_{cr})+2H_{cr}H(e_1(\alpha)+\alpha^2+\frac{1}{4}H^2+%
\frac{1}{2}W) \\
&&+\frac{3}{2}H_{cr}[2V(\alpha)+He_1(\alpha)]-H_{cr}^2H\},
\end{eqnarray*}
where 
\begin{eqnarray*}
\lefteqn{\frac{3}{2}V((H_{cr}^2)^{-\frac{1}{4}}(\frac{1}{2}%
e_1(H_{cr})+\alpha H_{cr}))} \\
&=&|H_{cr}|^{-\frac{1}{2}}\{\frac{3}{4}Ve_1(H_{cr})+\frac{3}{2}%
V(\alpha)H_{cr}+\frac{3}{4}\alpha V(H_{cr}) \\
&&-\frac{3}{8}|H_{cr}|^{-2}H_{cr}e_1(H_{cr})V(H_{cr})\},
\end{eqnarray*}
\begin{eqnarray*}
\lefteqn{\frac{1}{2}e_1e_1[(H_{cr}^2)^{-\frac{1}{4}}H_{cr}H]} \\
&=&\frac{1}{2}|H_{cr}|^{-\frac{1}{2}}(\frac{3}{2}e_1(H)e_1(H_{cr})+\frac{1}{2%
}He_1e_1(H_{cr})+e_1e_1(H)H_{cr}) \\
&&-\frac{1}{4}|H_{cr}|^{-\frac{1}{2}}[\frac{1}{2}%
|H_{cr}|^{-2}H_{cr}He_1(H_{cr})^2+e_1(H)e_1(H_{cr})].
\end{eqnarray*}
Therefore, 
\begin{eqnarray*}
\mathcal{E}_1&=&|H_{cr}|^{-\frac{1}{2}}\{\frac{3}{4}Ve_1(H_{cr})+\frac{3}{2}%
V(\alpha)H_{cr}+\frac{3}{4}\alpha V(H_{cr}) -\frac{3}{8}%
|H_{cr}|^{-2}H_{cr}e_1(H_{cr})V(H_{cr}) \\
&&+(\frac{1}{2}He_1(H_{cr})+e_1(H)H_{cr})\alpha \\
&&+\frac{1}{2}(\frac{3}{2}e_1(H)e_1(H_{cr})+\frac{1}{2}%
He_1e_1(H_{cr})+e_1e_1(H)H_{cr}) \\
&&-\frac{1}{4}[\frac{1}{2}%
|H_{cr}|^{-2}H_{cr}He_1(H_{cr})^2+e_1(H)e_1(H_{cr})] \\
&&-\frac{3}{2}\alpha H(\frac{1}{2}e_1(H_{cr})+\alpha
H_{cr})+2H_{cr}H(e_1(\alpha)+\alpha^2+\frac{1}{4}H^2+\frac{1}{2}W) \\
&&+\frac{3}{2}H_{cr}[2V(\alpha)+He_1(\alpha)]-H_{cr}^2H\}.
\end{eqnarray*}

Using (\ref{odealpha}), one can obtain 
\begin{eqnarray*}
&&\frac{3}{8}|H_{cr}|^{-2}H_{cr}e_{1}(H_{cr})V(H_{cr})+\frac{1}{8}%
|H_{cr}|^{-2}H_{cr}He_{1}(H_{cr})^{2} \\
&=&\frac{1}{4}sign(H_{cr})\mathfrak{f}e_{1}(H_{cr})-\frac{1}{4}%
e_{1}(H_{cr})(e_{1}(H)-\alpha H).
\end{eqnarray*}%
Therefore, by using (\ref{e-24-1}), we get 
\begin{eqnarray*}
\mathcal{E}_{1} &=&|H_{cr}|^{-\frac{1}{2}}\{-\frac{1}{4}sign(H_{cr})%
\mathfrak{f}e_{1}(H_{cr})+e_{1}[\frac{1}{2}e_{1}(H)H_{cr}+\frac{1}{4}%
He_{1}(H_{cr})+\frac{3}{4}V(H_{cr})] \\
&&+\frac{1}{4}\alpha He_{1}(H_{cr})+\frac{9}{4}\alpha V(H_{cr}) \\
&&+H_{cr}[\frac{9}{2}V(\alpha )+\frac{5}{2}He_{1}(\alpha )+\alpha e_{1}(H)+%
\frac{1}{3}H^{3}+\frac{3}{4}WH]\}.
\end{eqnarray*}%
Note that 
\begin{equation*}
\frac{1}{2}|H_{cr}|\mathfrak{f}\mathfrak{=}\frac{1}{2}e_{1}(H)H_{cr}+\frac{3%
}{4}V(H_{cr})+\frac{1}{4}He_{1}(H_{cr})-\frac{1}{2}\alpha HH_{cr}
\end{equation*}%
Therefore, 
\begin{eqnarray}
\mathcal{E}_{1} &=&|H_{cr}|^{-\frac{1}{2}}\{-\frac{1}{4}sign(H_{cr})%
\mathfrak{f}e_{1}(H_{cr})+\frac{1}{2}e_{1}(|H_{cr}|\mathfrak{f})
\label{Epsilon1} \\
&&+\frac{3}{2}|H_{cr}|\mathfrak{f}\alpha +H_{cr}[\frac{9}{2}V(\alpha
)+3HH_{cr}-\frac{1}{6}H^{3}]\}.  \notag
\end{eqnarray}%
It is easy to see that $\mathcal{E}_{1}$ is just the right hand side of (\ref%
{thesameEpsilon1}).

%TCIMACRO{\TeXButton{End Proof}{\endproof}}%
%BeginExpansion
\endproof%
%EndExpansion

\bigskip

\section{Examples}

\subsection{Examples of critical points of $E_{1}$}

First we consider some minimizers for the energy functional $E_{1}$ with
zero energy in $\mathcal{H}_{1}.$ That is, a surface satisfies the following
equation:%
\begin{equation}
H_{cr}=e_{1}(\alpha )+\frac{1}{2}\alpha ^{2}+\frac{1}{6}H^{2}=0  \label{2-1}
\end{equation}

\noindent by (\ref{e-1-1}) with $W$ $=$ $0$. Recall that we take the contact
form%
\begin{equation*}
\theta =dt+xdy-ydx
\end{equation*}

\noindent for $\mathcal{H}_{1},$ where $x,$ $y,$ $t$ are coordinates of $%
\mathcal{H}_{1}.$

\bigskip

\textbf{Example 1}. Consider the vertical planes defined by $ax+by$ $=$ $c$
in $\mathcal{H}_{1}.$ Observe that $T$ $=$ $\frac{\partial }{\partial t}$
lies in the tangent plane everywhere, so $\alpha $ $\equiv $ $0.$ Let $u$ $=$
$ax+by-c$ be a defining function. Compute%
\begin{eqnarray}
e_{2} &=&\frac{\nabla _{b}u}{|\nabla _{b}u|}  \label{2-2} \\
&=&\frac{a\mathring{e}_{1}+b\mathring{e}_{2}}{\sqrt{a^{2}+b^{2}}}  \notag
\end{eqnarray}

\noindent where the length $|\cdot |$ is with respect to the Levi metric $%
\theta ^{2}+dx^{2}+dy^{2}$ and 
\begin{equation*}
\mathring{e}_{1}:=\frac{\partial }{\partial x}+y\frac{\partial }{\partial t},%
\text{ }\mathring{e}_{2}:=\frac{\partial }{\partial y}-x\frac{\partial }{%
\partial t}
\end{equation*}

\noindent are standard left invariant (parallel with respect to the
pseudohermitian connection $\nabla )$ vector fields lying in contact planes
described by $\ker \theta .$ The subgradient $\nabla _{b}$ acting on a
function $u$ is given by $\nabla _{b}u$ $=$ ($\mathring{e}_{1}u)\mathring{e}%
_{1}+(\mathring{e}_{2}u)\mathring{e}_{2}.$ Recall that the CR structure $J$
is defined to satisfy $J\mathring{e}_{1}=\mathring{e}_{2}$ and $J\mathring{e}%
_{2}=-\mathring{e}_{1}.$ From (\ref{2-2}) we comput%
\begin{eqnarray*}
e_{1} &=&-Je_{2} \\
&=&\frac{b\mathring{e}_{1}-a\mathring{e}_{2}}{\sqrt{a^{2}+b^{2}}}
\end{eqnarray*}

We can then compute $H$ as follows:%
\begin{equation*}
\nabla e_{1}=\frac{b}{\sqrt{a^{2}+b^{2}}}\nabla \mathring{e}_{1}-\frac{a}{%
\sqrt{a^{2}+b^{2}}}\nabla \mathring{e}_{2}=0
\end{equation*}

\noindent since $\nabla \mathring{e}_{1}=\nabla \mathring{e}_{2}=0.$ It
follows that $H=0.$ So together with $\alpha $ $=$ $0$ the vertical planes
satisfy equation (\ref{2-1}).

\bigskip

For surfaces of the form $t^{2}$ $=$ $f(r^{2})$ in $\mathcal{H}_{1},$ where $%
r$ $=$ $\sqrt{x^{2}+y^{2}}$, we have formulas for $e_{1},$ $\alpha $ and $H$
as follows (see Section 3 in \cite{CCHY}; $e_{2}$ is chosen so that it
equals $\nabla _{b}u/|\nabla _{b}u|$ for $u$ $=$ $f-t^{2}$):%
\begin{eqnarray}
e_{1} &=&-Je_{2}=\frac{(yf^{\prime }+tx)\mathring{e}_{1}-(xf^{\prime }-ty)%
\mathring{e}_{2}}{r\sqrt{(f^{\prime })^{2}+f}},  \label{2-3} \\
\alpha &=&\frac{t}{r\sqrt{(f^{\prime })^{2}+f}},  \notag \\
H &=&\frac{r^{2}-f^{\prime }}{r\sqrt{(f^{\prime })^{2}+f}}-\frac{%
(1+2f^{\prime \prime })rf}{[(f^{\prime })^{2}+f]^{3/2}}  \notag
\end{eqnarray}

\noindent (note that the formula for $\alpha $ in Section 3 of \cite{CCHY}
also holds for the case $n=1).$ In terms of polar coordinates, we can
express $e_{1}$ as follows:%
\begin{equation}
e_{1}=\frac{-f^{\prime }\partial _{\theta }+tr\partial _{r}+r^{2}f^{\prime
}\partial _{t}}{r\sqrt{(f^{\prime })^{2}+f}}.  \label{2-3-1}
\end{equation}

\noindent Here we have used formulas: $y\partial _{x}-x\partial
_{y}=-\partial _{\theta },$ $x\partial _{x}+y\partial _{y}$ $=$ $r\partial
_{r}.$

\bigskip

\textbf{Example 2}. Consider shifted Heisenberg spheres defined by $%
(r^{2}+\lambda )^{2}+4t^{2}=\rho _{0}^{4}$ in $\mathcal{H}_{1},$ where $%
r^{2} $ $=$ $x^{2}+y^{2}$ and constants $\rho _{0}>0,$ $\lambda \geq 0$.
Write $t^{2}$ $=$ $\frac{1}{4}[\rho _{0}^{4}-(r^{2}+\lambda )^{2}]$ with 
\begin{equation*}
f(r^{2})=\frac{1}{4}[\rho _{0}^{4}-(r^{2}+\lambda )^{2}]\text{ or }f(s)=%
\frac{1}{4}[\rho _{0}^{4}-(s+\lambda )^{2}]
\end{equation*}

\noindent where $s=r^{2}.$ First we compute%
\begin{eqnarray}
\sqrt{(f^{\prime })^{2}+f} &=&\sqrt{\frac{1}{4}(s+\lambda )^{2}+\frac{1}{4}%
[\rho _{0}^{4}-(s+\lambda )^{2}]}  \label{2-4} \\
&=&\frac{1}{2}\rho _{0}^{2}.  \notag
\end{eqnarray}
\noindent By (\ref{2-3}), (\ref{2-3-1}) and (\ref{2-4}) we obtain, through a
direct computation,%
\begin{eqnarray*}
\alpha &=&\frac{2t}{\rho _{0}^{2}r}, \\
e_{1}(\alpha ) &=&\frac{4}{\rho _{0}^{4}}(-\frac{t^{2}}{r^{2}}-\frac{1}{2}%
(r^{2}+\lambda )), \\
H &=&\frac{1}{\rho _{0}^{2}}(3r+\frac{\lambda }{r}).
\end{eqnarray*}

\noindent We then compute%
\begin{eqnarray*}
H_{cr} &=&e_{1}(\alpha )+\frac{1}{2}\alpha ^{2}+\frac{1}{6}H^{2} \\
&=&...=\frac{-\frac{1}{2}\rho _{0}^{4}+\frac{2}{3}\lambda ^{2}}{\rho
_{0}^{4}r^{2}}=0
\end{eqnarray*}

\noindent if and only if $\lambda $ $=$ $\frac{\sqrt{3}}{2}\rho _{0}^{2}.$
Therefore Shifted Heisenberg spheres defined by $(r^{2}+\frac{\sqrt{3}}{2}%
\rho _{0}^{2})^{2}+4t^{2}=\rho _{0}^{4}$ are closed (compact with no
boundary) minimizers for the energy functional $E_{1}$ with zero energy.

\bigskip

\begin{conjecture}
\label{conj_3-1} The above shifted Heisenberg spheres defined by $(r^{2}+%
\frac{\sqrt{3}}{2}\rho _{0}^{2})^{2}+4t^{2}=\rho _{0}^{4}$ with $\rho _{0}$ $%
>$ $0$ are the only closed minimizers for $E_{1}$ (with zero energy) in $%
\mathcal{H}_{1}.$
\end{conjecture}

\bigskip

\textbf{Example 3}. Consider a family of surfaces defined by $t=cr^{2},$ $%
c\geq 0$ in $\mathcal{H}_{1},$ which are invariant under Heisenberg dilations%
$.$ Write $t^{2}=f(r^{2})$ with $f(s)$ $=$ $c^{2}s^{2}.$ So $\sqrt{%
(f^{\prime })^{2}+f}$ $=$ $c\sqrt{4c^{2}+1}r^{2}.$ A direct computation shows%
\begin{eqnarray}
H &=&-\frac{2c}{\sqrt{4c^{2}+1}}\frac{1}{r},  \label{2-4-1} \\
\alpha &=&\frac{1}{\sqrt{4c^{2}+1}}\frac{1}{r},  \notag \\
e_{1}(\alpha ) &=&-\frac{1}{4c^{2}+1}\frac{1}{r^{2}}  \notag
\end{eqnarray}

\noindent by (\ref{2-3}) and (\ref{2-3-1}). We then have%
\begin{eqnarray}
H_{cr} &=&e_{1}(\alpha )+\frac{1}{2}\alpha ^{2}+\frac{1}{6}H^{2}  \label{2-5}
\\
&=&\frac{\frac{2}{3}c^{2}-\frac{1}{2}}{4c^{2}+1}\frac{1}{r^{2}}  \notag
\end{eqnarray}

\noindent which vanishes if and only if $c$ $=$ $\frac{\sqrt{3}}{2}.$ We
conclude that $t=\frac{\sqrt{3}}{2}r^{2}$ is a minimizer for $E_{1}$ with
zero energy.

We are going to examine that for which $c,$ $c$ $\neq $ $\frac{\sqrt{3}}{2},$
$t=cr^{2}$ is still a solution to the Euler-Lagrange equation (\ref{e-4}).
Recall that 
\begin{eqnarray}
h_{11} &=&H=-\frac{2c}{\sqrt{4c^{2}+1}}\frac{1}{r},\text{ }  \label{2-6} \\
h_{10} &=&e_{1}(\alpha )+\frac{1}{2}\alpha ^{2}=-\frac{1}{2}\frac{1}{4c^{2}+1%
}\frac{1}{r^{2}}  \notag
\end{eqnarray}

\noindent from (\ref{2-4-1}). On the other hand, observe that $2\phi
_{1r}^{1}$ $=$ $\phi $ $=$ $-2\alpha e^{1}$ (by (\ref{1-9-1}); we have
replaced the notations $\tilde{\phi}_{1r}^{1},$ $\tilde{\phi}$ by $\phi
_{1r}^{1},$ $\phi ,$ resp. here) and $\phi _{r}^{1}$ $=$ $0.$ From (\ref{e-2}%
) we compute%
\begin{eqnarray}
h_{111} &=&\frac{2c}{4c^{2}+1}\frac{1}{r^{2}}+\alpha \frac{-2c}{\sqrt{%
4c^{2}+1}}\frac{1}{r}=0,  \label{2-7} \\
h_{110} &=&\frac{4c^{2}}{(4c^{2}+1)^{3/2}}\frac{1}{r^{3}},  \notag \\
h_{100} &=&...=\frac{2c}{(4c^{2}+1)^{2}}\frac{1}{r^{4}}.  \notag
\end{eqnarray}

\noindent We can now compute, by (\ref{2-5}), (\ref{2-6}) and (\ref{2-7}),%
\begin{eqnarray*}
\mathfrak{f} &\mathfrak{:}&=|H_{cr}|^{-1}\{h_{10}h_{111}+\frac{1}{3}%
h_{11}^{2}h_{111} \\
&&+h_{11}h_{110}+\frac{3}{2}h_{100}\} \\
&=&...=\frac{6c(3-8c^{2})}{|4c^{2}-3|(4c^{2}+1)}\frac{1}{r^{2}}
\end{eqnarray*}

\noindent (cf. (\ref{e-3})). We also need to know%
\begin{eqnarray*}
h_{00} &=&(T+\alpha e_{2})(\alpha ) \\
&=&...=\frac{-2c}{(4c^{2}+1)^{3/2}}\frac{1}{r^{3}}.
\end{eqnarray*}

\noindent It follows that%
\begin{eqnarray*}
&&9h_{00}+6h_{11}h_{10}+\frac{2}{3}h_{11}^{3} \\
&=&...=-\frac{(36+16c^{2})c}{3(4c^{2}+1)^{3/2}}\frac{1}{r^{3}}.
\end{eqnarray*}

\noindent It turns out that the first two terms of $\sigma ,$ involving $%
e^{1}$ cancel (see (\ref{e-4})). So we end up to conclude that the
Euler-Lagrange equation (\ref{e-4}) holds if and only if%
\begin{equation*}
\frac{|4c^{2}-3|^{1/2}(36+16c^{2})c}{3\sqrt{6}(4c^{2}+1)^{2}}=0.
\end{equation*}

\noindent For $c$ $\neq $ $\frac{\sqrt{3}}{2},$ we get $c$ $=$ $0.$
Therefore the surface defined by $t=0$ is a solution to the Euler-Lagrange
equation of $E_{1}$ (of higher energy level).

We remark that under the (inverse) Cayley transform, $\{t=0\}$ is mapped to
a great 2-sphere of $S^{3}$ (A great 2-sphere is the intersection of $S^{3}$
and a $3$-plane in $R^{4}$ passing through the origin)$.$ All great
2-spheres of $S^{3}$ are CR equivalent. One of them corresponds to the
distance sphere $\rho $ $=$ $1$ ($\rho $ $=$ ($r^{4}+4t^{2})^{1/4})$ in $%
\mathcal{H}_{1}$, so it (and hence $\rho $ $=$ $\rho _{0}$ $>$ $0$ due to
dilations in $\mathcal{H}_{1}$ being $CR$ transformations) is a critical
point of higher energy level with respect to $E_{1}.$ In fact, we can verify
that (\ref{e-10}) and (\ref{e-11}) vanish for distance spheres (or
Heisenberg spheres) $\rho $ $=$ $\rho _{0}$ $>$ $0$ by a direct computation.

\bigskip

\textbf{Example 4}. Define a closed surface $\Sigma _{c}$ in $S^{3}$ $%
\subset $ $C^{2}$ by $\rho _{1}=c$ where ($z_{1},z_{2})$ $\in $ $C^{2}$ such
that $S^{3}$ is described by%
\begin{eqnarray*}
z_{1} &=&\rho _{1}e^{i\varphi _{1}},z_{2}=\rho _{2}e^{i\varphi _{2}}, \\
\rho _{1}^{2}+\rho _{2}^{2} &=&1.
\end{eqnarray*}%
\noindent The contact form $\hat{\Theta}$ on $S^{3}$ reads%
\begin{eqnarray*}
\hat{\Theta} &=&\func{Im}\partial (|z_{1}|^{2}+|z_{2}|^{2}) \\
&=&x_{1}dy_{1}-y_{1}dx_{1}+x_{2}dy_{2}-y_{2}dx_{2} \\
&=&\rho _{1}^{2}d\varphi _{1}+\rho _{2}^{2}d\varphi _{2}.
\end{eqnarray*}%
\noindent It is not hard to see that the Reeb vector field 
\begin{equation*}
T=\frac{\partial }{\partial \varphi _{1}}+\frac{\partial }{\partial \varphi
_{2}},
\end{equation*}%
\noindent and hence $T$ $\in $ $T\Sigma _{c}.$ So we get

\begin{equation}
\alpha \equiv 0  \label{torusalpha}
\end{equation}
\noindent on $\Sigma _{c}.$ On the other hand, we can express the CR
structure $\hat{J}$ in terms of polar coordinates by%
\begin{equation*}
\hat{J}(\rho _{j}\frac{\partial }{\partial \rho _{j}})=\frac{\partial }{%
\partial \varphi _{j}},\hat{J}\frac{\partial }{\partial \varphi _{j}}=-\rho
_{j}\frac{\partial }{\partial \rho _{j}}
\end{equation*}%
\noindent for $j=1,2$ (note that $\rho _{j}\frac{\partial }{\partial \rho
_{j}}$ $=$ $2\func{Re}(z_{j}\frac{\partial }{\partial z_{j}}),$ $\frac{%
\partial }{\partial \varphi _{j}}$ $=$ $2\func{Re}(iz_{j}\frac{\partial }{%
\partial z_{j}})).$ For $\Sigma _{c}$ we can easily find%
\begin{eqnarray*}
e_{1} &=&-\frac{\rho _{2}}{\rho _{1}}\frac{\partial }{\partial \varphi _{1}}+%
\frac{\rho _{1}}{\rho _{2}}\frac{\partial }{\partial \varphi _{2}}, \\
e_{2} &=&\hat{J}e_{1}=\rho _{2}\frac{\partial }{\partial \rho _{1}}-\rho _{1}%
\frac{\partial }{\partial \rho _{2}}.
\end{eqnarray*}%
\noindent It follows that 
\begin{eqnarray*}
e^{1} &=&-\rho _{1}\rho _{2}(d\varphi _{1}-d\varphi _{2}), \\
e^{2} &=&\rho _{2}d\rho _{1}-\rho _{1}d\rho _{2}.
\end{eqnarray*}%
\noindent Hence the $p$-area form and the volume form read%
\begin{eqnarray}
\hat{\Theta}\wedge e^{1} &=&\rho _{1}\rho _{2}d\varphi _{1}\wedge d\varphi
_{2},  \label{2-8} \\
\hat{\Theta}\wedge e^{1}\wedge e^{2} &=&\rho _{1}d\rho _{1}\wedge d\varphi
_{1}\wedge d\varphi _{2}.  \notag
\end{eqnarray}%
\noindent Here we have used $\rho _{1}^{2}+\rho _{2}^{2}=1.$ We compute $H$
as follows:%
\begin{eqnarray}
d(\hat{\Theta}\wedge e^{1}) &=&(\rho _{2}^{2}-\rho _{1}^{2})\rho
_{2}^{-1}d\rho _{1}\wedge d\varphi _{1}\wedge d\varphi _{2}  \label{2-9} \\
&=&-H\hat{\Theta}\wedge e^{1}\wedge e^{2}  \notag
\end{eqnarray}%
\noindent So comparing (\ref{2-8}) with (\ref{2-9}), we get%
\begin{equation}
H=-(\rho _{2}^{2}-\rho _{1}^{2})\rho _{2}^{-1}\rho _{1}^{-1}  \label{torusH}
\end{equation}%
\noindent for $\rho _{1}\rho _{2}\neq 0.$ Consider the Clifford torus $%
\Sigma _{\frac{\sqrt{2}}{2}}$ in $S^{3}$ with $c$ $=$ $\rho _{1}$ $=$ $\rho
_{2}$ $=$ $\frac{\sqrt{2}}{2}.$ We get $H=0.$ From (\ref{1-13}) we compute
the first energy of $\Sigma _{\frac{\sqrt{2}}{2}}:$ 
\begin{eqnarray*}
&&E_{1}(\Sigma _{\frac{\sqrt{2}}{2}}) \\
&=&\int_{\Sigma _{\frac{\sqrt{2}}{2}}}|e_{1}(\alpha )+\frac{1}{2}\alpha ^{2}-%
\func{Im}A_{11}+\frac{1}{4}W+\frac{1}{6}H^{2}|^{3/2}\hat{\Theta}\wedge e^{1}
\\
&=&\int_{\Sigma _{\frac{\sqrt{2}}{2}}}|\frac{1}{4}W|^{3/2}\rho _{1}\rho
_{2}d\varphi _{1}\wedge d\varphi _{2}\text{ (since }\alpha =H=0,A_{11}=0) \\
&=&(\frac{1}{2})^{3/2}(\frac{\sqrt{2}}{2})^{2}(2\pi )^{2}\text{ (since }%
W=2,\rho _{1}=\rho _{2}=\frac{\sqrt{2}}{2}) \\
&=&\frac{\pi ^{2}}{\sqrt{2}}.
\end{eqnarray*}

Next we claim that $\Sigma _{\frac{\sqrt{2}}{2}}$ satisfies the
Euler-Lagrange equation for $E_{1}.$ Observe by (\ref{1-8}) that $h_{11}$ $=$
$H$ $=$ $0,$ $h_{10}$ $=$ $e_{1}(\alpha )$ $+$ $\frac{1}{2}\alpha ^{2}$ $-$ $%
\func{Im}A_{11}$ $+$ $\frac{1}{4}W$ $=$ $\frac{1}{4}W$ = $\frac{1}{2}$ since 
$\alpha $ $=$ $0,$ $A_{11}$ $=$ $0$ and $W$ $=$ $2.$ Also $\phi _{r}^{1}$
(or $\tilde{\phi}_{r}^{1})$ $=$ $-\frac{W}{4}\alpha \hat{\Theta}$ $=$ $0,$ $%
\phi $ (or $\tilde{\phi})$ $=$ -2$\alpha e^{1}$ $=$ $0.$ It follows that $%
h_{111}$ $=$ $h_{100}$ $=$ $h_{110}$ $=$ $0.$ So $\mathfrak{f}$ in (\ref{e-3}%
) vanishes (which can also be proved by the formula (\ref{e-10}) due to $H$ $%
=$ $0$ and $\alpha $ $=$ $0)$. From (\ref{1-12}) we have%
\begin{eqnarray*}
h_{00} &=&(T+\alpha e_{2})(\alpha ) \\
&&+\func{Im}[\frac{1}{6}W^{,1}+\frac{2i}{3}(A^{11})_{,1}]-\alpha (\func{Re}%
A_{\bar{1}}^{1}) \\
&=&0
\end{eqnarray*}

\noindent since $\alpha $ $=$ $0,$ $A_{11}$ $=$ $0$ and $W$ $=$ $2.$
Altogether we have shown that (\ref{e-4}) holds, i.e., the Clifford torus $%
\Sigma _{\frac{\sqrt{2}}{2}}$ satisfies the Euler-Lagrange equation for $%
E_{1}.$

\bigskip

\begin{conjecture}
\label{conj_3-2} The Clifford torus $\Sigma _{\frac{\sqrt{2}}{2}}$ $\subset $
$S^{3}$ is a unique minimizer among all surfaces of torus type (genus=1) for 
$E_{1}$ up to $CR$ automorphisms of $S^{3}.$
\end{conjecture}

\bigskip

We remark that one can see that if $T^{2}\subset S^{3}$ is a torus without
singular points, then it cannot have $E_{1}=0$. Otherwise 
\begin{equation*}
e_{1}(\alpha )=-\frac{1}{2}\alpha ^{2}-\frac{1}{6}H^{2}-\frac{1}{2}\leq -%
\frac{1}{2},
\end{equation*}%
so that any characteristic curve of $e_{1}$ is open and along the curve $%
\alpha $ goes to infinity, which contradicts to $\alpha $ being bounded on a
nonsingular compact surface. We compute the value of $E_{1}$ 
\begin{equation*}
\int (\frac{1}{6}H^{2}+\frac{W}{4})^{\frac{3}{2}}\rho _{1}\rho _{2}d\varphi
_{1}d\varphi _{2},
\end{equation*}%
which amounts to compare the value of 
\begin{equation*}
(H^{2}+3)(\rho _{1}\rho _{2})^{\frac{2}{3}}=x^{-\frac{4}{3}}-x^{\frac{2}{3}%
},\quad x=\rho _{1}\rho _{2}\in (0,1/2],
\end{equation*}%
which is decreasing in $x$, hence we have the minimal value (when $\rho
_{1}=\rho _{2}=\frac{1}{\sqrt{2}},H=0$) 
\begin{equation*}
\min_{\Sigma _{\rho _{1},\rho _{2}}}\int (\frac{1}{6}H^{2}+\frac{W}{4})^{%
\frac{3}{2}}\rho _{1}\rho _{2}d\varphi _{1}d\varphi _{2}=\frac{\sqrt{2}}{2}%
\pi ^{2}.
\end{equation*}

\bigskip

\subsection{Examples of critical points of $E_{2}$}

For critical points of $E_{2},$ we have the following examples.

\bigskip

\textbf{Example 1}. Consider the vertical planes defined by $ax+by$ $=$ $c$
in $\mathcal{H}_{1}.$ As we have seen, the vertical planes have $\alpha =0$
and $H=0$. So the vertical planes satisfy the Euler-Lagrange equation (\ref%
{EulerLagrange2}) (see (\ref{mathcalE}) for the definition of $\mathcal{E}%
_{2}$).

\bigskip

\textbf{Example 2}. Consider a surface foliated by a linear combination of $%
\mathring{e}_{1}$ and $\mathring{e}_{2}$ where 
\begin{equation*}
\mathring{e}_{1}:=\frac{\partial }{\partial x}+y\frac{\partial }{\partial t},%
\text{ }\mathring{e}_{2}:=\frac{\partial }{\partial y}-x\frac{\partial }{%
\partial t}.
\end{equation*}%
If $\Sigma $ is foliated by a distribution $\mathbb{R}X$, where $X$ is a
unit vector field 
\begin{equation*}
X=a\mathring{e}_{1}+b\mathring{e}_{2},\quad a,b\in \mathbb{R},\quad
a^{2}+b^{2}=1,
\end{equation*}%
then we can take along $\Sigma $ 
\begin{equation*}
e_{1}=X,\quad e_{2}=-b\mathring{e}_{1}+a\mathring{e}_{2}.
\end{equation*}%
Note that 
\begin{equation*}
\nabla _{\mathring{e}_{i}}\mathring{e}_{j}=0,
\end{equation*}%
so we take $e_{2}=-b\mathring{e}_{1}+a\mathring{e}_{2}$ near $\Sigma $. Then
near $\Sigma $, $e_{1}=-Je_{2}=a\mathring{e}_{1}+b\mathring{e}_{2}$ and 
\begin{equation*}
H=\langle \nabla _{e_{1}}e_{1},e_{2}\rangle =0.
\end{equation*}%
Note that $\nabla _{T}e_{1}=0$, hence on $\Sigma $ 
\begin{equation*}
\omega (T)=e_{1}(\alpha )+2\alpha ^{2}=0.
\end{equation*}%
Hence $\Sigma $ satisfies the Euler-Lagrange equation (\ref{EulerLagrange2}).

In particular 
\begin{equation*}
\Sigma _{+}:=\{t-xy=c\}\text{ and }\Sigma _{-}:=\{t+xy=c\}
\end{equation*}%
are surfaces foliated by the distribution $\mathbb{R}\mathring{e}_{1}$ and $%
\mathbb{R}\mathring{e}_{2}$ respectively.

\bigskip

\textbf{Example 3}. Consider a surface in $\mathcal{H}_{1}$ defined by $%
t=cr^{2},c>0$. As we have seen, 
\begin{eqnarray}
H &=&-\frac{2c}{\sqrt{4c^{2}+1}}\frac{1}{r},  \label{4-2-0} \\
\alpha &=&\frac{1}{\sqrt{4c^{2}+1}}\frac{1}{r}.  \notag
\end{eqnarray}
Moreover, we have%
\begin{eqnarray}
e_{1}(\alpha ) &=&-\alpha ^{2},\quad e_{1}(H)=2c\alpha ^{2},  \label{4-2-1}
\\
V(H) &=&4c^{2}\alpha ^{3},\quad e_{1}V(H)=-12c^{2}\alpha ^{4},  \notag \\
e_{1}e_{1}(H) &=&-4c\alpha ^{3}.  \notag
\end{eqnarray}
Substituting (\ref{4-2-1}), $H$ $=$ $-2c\alpha $ (from (\ref{4-2-0})) and $W$
$=$ $0$ into (\ref{mathcalE}), we find that when $c^{2}=\frac{3}{4},$ (\ref%
{EulerLagrange2}) holds. So we have the following critical point of $E_{2}:$ 
\begin{equation*}
\Sigma =\{t=\frac{\sqrt{3}}{2}(x^{2}+y^{2})\}.
\end{equation*}

\bigskip

Remark 1. $\{t=0\}\subset \mathcal{H}_{1}$ and the distance sphere $\{\frac{1%
}{4}|z|^{4}+t^{2}=1\}\subset \mathcal{H}_{1}$ are not critical for $E_{2}$.
Equivalently, the great 2-spheres of $S^{3}$ are not critical for $E_{2}$.

Remark 2. $E_{2}$ is unbounded from below and above in general. Consider $%
\Sigma _{\rho _{1}}$ in $S^{3}$ $\subset $ $C^{2}$ (see Example 4 in the
last subsection) as follows:%
\begin{equation*}
\Sigma _{\rho _{1}}=\{|z_{1}|^{2}=x_{1}^{2}+y_{1}^{2}=\rho _{1}^{2}\}.
\end{equation*}%
From (\ref{torusalpha}), (\ref{torusH}) and (\ref{2-8}), we learn 
\begin{eqnarray}
\alpha &=&0,  \label{rmk2} \\
H &=&\frac{\rho _{1}^{2}-\rho _{2}^{2}}{\rho _{1}\rho _{2}}=\frac{\rho _{1}}{%
\rho _{2}}-\frac{\rho _{2}}{\rho _{1}},  \notag \\
\hat{\Theta}\wedge e^{1} &=&\rho _{1}\rho _{2}d\varphi _{1}\wedge d\varphi
_{2}  \notag
\end{eqnarray}
\noindent on $\Sigma _{\rho _{1}}.$ By (\ref{rmk2}) and $W$ $=$ $2$ on $%
S^{3},$ we can then compute 
\begin{eqnarray*}
E_{2}(\Sigma _{\rho _{1}}) &=&\int_{\Sigma }(\frac{1}{6}WH+\frac{2}{27}H^{3})%
\hat{\Theta} \wedge e^{1} \\
&=&\int_{\Sigma }(\frac{1}{3}+\frac{2}{27}H^{2})(\rho _{1}^{2}-\rho
_{2}^{2})d\varphi _{1}\wedge d\varphi _{2} \\
&=&(\frac{1}{3}+\frac{2}{27}H^{2})(\rho _{1}^{2}-\rho _{2}^{2})(2\pi )^{2}
\end{eqnarray*}%
which tends to $+\infty $ ($-\infty ,$ resp.) as $\rho _{1}\rightarrow 1$ ($%
\rho _{1}\rightarrow 0$, resp.). On the other hand, from (\ref{mathcalE}) we
compute%
\begin{equation*}
\mathcal{E}_{2}(\Sigma _{\rho _{1}})=\frac{4}{9}(\frac{1}{3}%
H^{4}+2H^{2}+3)>0.
\end{equation*}
\noindent for all $0<\rho _{1}<1.$ So all $\Sigma _{\rho _{1}},0<\rho _{1}<1$%
, are not critical points of $E_{2}$.

\bigskip

\section{Singular CR Yamabe problem: formal solutions}

In this section, we consider formal solutions to the singular CR Yamabe
problem on a 3-dimensional CR manifold with boundary. In particular, we get
the expansion of the formal solution in a specific defining function for the
boundary. In the next section, we will prove that $E_{2}$ can be interpreted
as the coefficient of the log term in the volume renormalization for a
formal solution to the singular CR Yamabe problem.

\bigskip

\subsection{A defining function for the boundary and a frame near the
boundary}

Let $(M,J,\theta )$ be a $3$-dimensional pseudohermitian manifold with
boundary $\Sigma $ $=$ $\partial M$. For other notations and some basic
formulas, we refer the reader to the Appendix. Let $x$ be a regular (or
nonsingular) point of $\Sigma $, that is 
\begin{equation*}
T_{x}\Sigma \neq \xi _{x}:=ker\theta _{x},
\end{equation*}%
i.e. 
\begin{equation*}
T\neq \pm \nu ,
\end{equation*}%
where $\nu $ is the unit inner normal of $\Sigma $ with respect to $%
g_{\theta }$. Let $e_{1}$ be the unit vector in $T_{x}\Sigma \cap \xi _{x}$
and $e_{2}=Je_{1}$ such that 
\begin{equation}
\langle e_{2},\nu \rangle >0.  \label{orientatione2}
\end{equation}%
There exists $\alpha (x)$ $\in $ $R$ such that 
\begin{equation*}
V:=T+\alpha e_{2}\in T_{x}\Sigma .
\end{equation*}%
Note that 
\begin{equation*}
\nu =\frac{e_{2}-\alpha T}{\sqrt{1+\alpha ^{2}}},\quad \langle e_{2},\nu
\rangle =\frac{1}{\sqrt{1+\alpha ^{2}}},\quad \langle T,\nu \rangle =\frac{%
-\alpha }{\sqrt{1+\alpha ^{2}}}.
\end{equation*}

Let $x$ be a regular point of $\Sigma $ and $e_{2}\in \xi _{x}$ be chosen as
above. We now define a geodesic $\gamma (\rho )$, with $\dot{\gamma}=\frac{d%
}{d\rho }\gamma (\rho )\in \xi =\ker \theta $, by 
\begin{equation}
\left\{ 
\begin{array}{ll}
\nabla _{\dot{\gamma}}\dot{\gamma}=0, &  \\ 
\gamma (0)=x,\quad \dot{\gamma}(0)=e_{2} & 
\end{array}%
\right.  \label{geodesic}
\end{equation}%
Note that $\dot{\gamma}$ is a unit vector. Denote $\exp _{x}(\rho
e_{2})=\gamma (\rho )$, and $\rho $ is the defining function such that 
\begin{equation*}
\rho (\exp _{x}\rho e_{2})=\rho .
\end{equation*}%
We define the level surface 
\begin{equation*}
\Sigma _{\rho }=\{\exp _{x}(\rho e_{2}),\quad x\in \Sigma \}.
\end{equation*}%
We can extend $\{e_{1},e_{2}\}$ from the boundary $\Sigma $ to an
orthonormal frame $\{e_{1},e_{2}\}$ of $\xi $ into the interior of $M$ as
follows: at $y=\gamma (\rho )=\exp _{x}(\rho e_{2})$ let 
\begin{equation}
e_{2}(y):=\dot{\gamma}(\rho ),\quad e_{1}(y):=-Je_{2}(y).
\label{geodesicframe}
\end{equation}%
The region near $\Sigma $ can be endowed with the natural coordinate charts $%
(x^{1},x^{2},\rho )$ and $e_{2}=\partial _{\rho }$. Throughout this section, 
$\{e_{1},e_{2}\}$ is given by (\ref{geodesicframe}) unless specified
otherwise. Note that $e_{1}$ in this frame may not be tangent to the level
sets of $\rho $ other than $\rho $ $=$ $0.$

\bigskip

\subsection{Formal solutions to the singular CR Yamabe problem}

Let 
\begin{equation*}
\widetilde{\theta }=u^{-2}\theta ,
\end{equation*}%
where $u>0$ in the interior of $M$. We consider the following singular CR
Yamabe problem: 
\begin{equation}
\left\{ 
\begin{array}{ll}
\widetilde{R}=-8,\text{on}\quad M &  \\ 
u=0,\quad \text{on}\,\ \Sigma . & 
\end{array}%
\right.  \label{CRYamabe}
\end{equation}%
Here $\widetilde{R}=2\widetilde{W}$ and $\widetilde{W}$ is the Webster
scalar curvature of $\widetilde{\theta}$. Note that

\begin{eqnarray}
\widetilde{R}&=&R_{\widetilde{\theta}}=(-4\triangle_Pu^{-1}+Ru^{-1})u^{3} 
\notag \\
&=&4u\triangle_Pu-8|d_Pu|^2+Ru^2.  \label{transformationlawofWS}
\end{eqnarray}
where 
\begin{eqnarray}  \label{sublaplacian}
\triangle_Pu&=&e_1e_1(u)+e_2e_2(u)-\omega(e_1)e_2(u)+\omega(e_2)e_1(u) 
\notag \\
&=&e_1e_1(u)+e_2e_2(u)-\omega(e_1)e_2(u)
\end{eqnarray}
and 
\begin{equation}
|d_Pu|^2=e_1(u)^2+e_2(u)^2.
\end{equation}

There is a formal solution of the form 
\begin{equation}
u(x,\rho )=c(x)\rho +v(x)\rho ^{2}+w(x)\rho ^{3}+z(x)\rho ^{4}+l(x)\rho
^{5}\log \rho +h(x)\rho ^{5}+O(\rho ^{6})  \label{formalsolution}
\end{equation}%
to the singular CR Yamabe problem (\ref{CRYamabe}). $u$ satisfies 
\begin{equation}
\left\{ 
\begin{array}{ll}
|d_{P}u|^{2}-\frac{1}{2}u\triangle _{P}u-\frac{1}{8}Ru^{2}=1+O(\rho
^{5}\log\rho), &  \\ 
u=0,\,\ \mbox{on}\,\ \Sigma . & 
\end{array}%
\right.  \label{formalCRYamabe}
\end{equation}%
In order to compute the renormalized volume of the solution $\widetilde{%
\theta }$, we need to identify the coefficients $c(x),v(x),w(x)$ and $z(x)$
in (\ref{formalsolution}). $h(x)$ is not locally and formally determined and 
$l(x)\neq 0$ is the obstruction to the smoothness of the solution of (\ref%
{CRYamabe}). Actually, we can let $h(x)\equiv 0$ in the formal solution.

\bigskip

\subsection{The expansion of a formal solution to the singular CR Yamabe
problem}

Let 
\begin{equation*}
u(x,\rho)=u[\exp_x(\rho e_2)], \quad x\in \Sigma.
\end{equation*}
By (\ref{formalCRYamabe}), at $\rho=0$ 
\begin{equation*}
|d_Pu|^2=e_2(u)^2=u_\rho^2=1,
\end{equation*}
and according to (\ref{orientatione2}), we have 
\begin{equation*}
u(x,\rho)=\rho+\cdots
\end{equation*}

\begin{proposition}
We have 
\begin{equation}  \label{v}
v(x)=-\frac{1}{6}H(x),
\end{equation}
where $H$ is the p-mean curvature of $\Sigma$, defined by 
\begin{equation*}
H\equiv\omega(e_1)=\langle \nabla_{e_1}e_1,e_2\rangle.
\end{equation*}
\end{proposition}

%TCIMACRO{\TeXButton{Proof}{\proof} }%
%BeginExpansion
\proof
%EndExpansion
Let 
\begin{equation*}
u(x,\rho )=\rho +v(x)\rho ^{2}+\cdots
\end{equation*}%
where $v(x)$ will be determined by 
\begin{equation*}
\frac{d}{d\rho }|_{\rho =0}[|d_{P}u|^{2}-\frac{1}{2}u\triangle _{P}u-\frac{1%
}{8}Ru^{2}]=0.
\end{equation*}%
We have 
\begin{eqnarray*}
\frac{d}{d\rho }(|d_{P}u|^{2}) &=&\frac{d}{d\rho }[e_{1}(u)^{2}+e_{2}(u)^{2}]
\\
&=&2e_{1}(u)e_{2}e_{1}(u)+2e_{2}(u)e_{2}e_{2}(u),
\end{eqnarray*}%
hence it follows from the fact that at $\Sigma $ 
\begin{equation*}
e_{1}(u)=0,\quad e_{2}(u)=1
\end{equation*}%
that 
\begin{equation*}
\frac{d}{d\rho }|_{\rho =0}(|d_{P}u|^{2})=2u_{\rho \rho }.
\end{equation*}

Recall (\ref{sublaplacian}) 
\begin{eqnarray*}
\triangle_Pu&=&e_1e_1(u)+e_2e_2(u)-\omega(e_1)e_2(u),
\end{eqnarray*}
hence 
\begin{equation}
(\triangle_Pu)|_{\rho=0}=u_{\rho\rho}-H,  \label{laplacianu}
\end{equation}
so we have 
\begin{eqnarray*}
\frac{d}{d\rho}|_{\rho=0}[u\triangle_Pu]=(\triangle_Pu)|_{\rho=0}=u_{\rho%
\rho}-H.
\end{eqnarray*}

It then follows from 
\begin{equation*}
\frac{d}{d\rho }|_{\rho =0}[|d_{P}u|^{2}-\frac{1}{2}u\triangle _{P}u-\frac{1%
}{8}Ru^{2}]=0
\end{equation*}%
that 
\begin{equation*}
2u_{\rho \rho }-\frac{1}{2}(u_{\rho \rho }-H)=0,
\end{equation*}%
hence 
\begin{equation*}
u_{\rho \rho }=-\frac{1}{3}H,
\end{equation*}%
\begin{equation*}
u(x,\rho )=\rho -\frac{1}{6}H(x)\rho ^{2}+\cdots
\end{equation*}

%TCIMACRO{\TeXButton{End Proof}{\endproof}}%
%BeginExpansion
\endproof%
%EndExpansion

\bigskip

Let 
\begin{equation*}
u(x,\rho )=\rho -\frac{1}{6}H(x)\rho ^{2}+w(x)\rho ^{3}+\cdots
\end{equation*}%
$w(x)$ will be determined by 
\begin{equation*}
\frac{d^{2}}{d\rho ^{2}}|_{\rho =0}[|d_{P}u|^{2}-\frac{1}{2}u\triangle _{P}u-%
\frac{1}{8}Ru^{2}]=0.
\end{equation*}%
Note that for $\{e_{1},e_{2}\}$ given by (\ref{geodesicframe}), 
\begin{equation*}
\nabla _{e_{2}}e_{2}=0,\quad i.e.\quad \omega (e_{2})=0.
\end{equation*}%
One can easily check the following

\bigskip

\begin{lemma}
We have 
\begin{eqnarray}  \label{commutators}
&&\nabla_{e_2}e_1=0, \quad \nabla_{e_2}e_2=0,  \notag \\
&&\nabla_{e_1}e_1=\omega(e_1)e_2, \quad \nabla_{e_1}e_2=-\omega(e_1)e_1, 
\notag \\
&& [e_2,e_1]=\omega(e_1)e_1+2T, \quad T_\nabla(e_1,e_2)=2T,  \notag \\
&&[T, e_1]=-a_1e_1+[a_2+\omega(T)]e_2, \\
&& [e_2,T]=(-a_2+\omega(T))e_1-a_1e_2,  \notag \\
&& e_2\omega(e_1)=R+\omega(e_1)^2+2\omega(T).  \notag
\end{eqnarray}
\end{lemma}

\begin{lemma}
At $\Sigma $, we have 
\begin{eqnarray}
e_{1}e_{2}(u) &=&0,\quad e_{2}e_{1}(u)=-2\alpha ,\quad \triangle _{P}u=8v, 
\notag \\
\omega (T) &=&e_{1}(\alpha )+2\alpha ^{2}-a_{2},  \notag \\
e_{1}T(u) &=&-e_{1}(\alpha ),\quad Te_{1}(u)=2\alpha ^{2},  \label{lemmaofwz}
\\
Te_{2}(u) &=&-2\alpha v,\quad e_{2}T(u)=-2\alpha v-a_{1},  \notag \\
e_{2}e_{1}e_{1}(u) &=&-2e_{1}(\alpha )+4\alpha ^{2}.  \notag
\end{eqnarray}
\end{lemma}

\bigskip

%TCIMACRO{\TeXButton{Proof}{\proof} }%
%BeginExpansion
\proof
%EndExpansion
By (\ref{commutators}), 
\begin{equation*}
e_{2}e_{1}(u)=e_{1}e_{2}(u)+\omega _{1}(e_{1})e_{1}(u)+2T(u).
\end{equation*}%
In particular, on $\Sigma $ 
\begin{equation*}
e_{2}e_{1}(u)=2T(u)=2(V-\alpha e_{2})(u)=-2\alpha .
\end{equation*}%
We have at $\Sigma $, 
\begin{equation*}
Te_{1}(u)=e_{1}Tu+[T,e_{1}]u=-e_{1}(\alpha )+\omega (T)+a_{2}.
\end{equation*}%
On the other hand, 
\begin{equation*}
Te_{1}(u)=(V-\alpha e_{2})e_{1}u=2\alpha ^{2}.
\end{equation*}%
Hence 
\begin{equation*}
\omega (T)=e_{1}(\alpha )+2\alpha ^{2}-a_{2}.
\end{equation*}%
We have at $\Sigma $, 
\begin{eqnarray*}
Te_{2}(u) &=&(V-\alpha e_{2})e_{2}(u)=-2\alpha v, \\
e_{2}T(u) &=&Te_{2}(u)+[e_{2},T](u)=-2\alpha v-a_{1}.
\end{eqnarray*}%
Finally, 
\begin{eqnarray*}
e_{2}e_{1}e_{1}(u) &=&e_{1}e_{2}e_{1}(u)+[e_{2},e_{1}]e_{1}(u) \\
&=&e_{1}e_{1}e_{2}(u)+e_{1}[e_{2},e_{1}](u)+[e_{2},e_{1}]e_{1}(u).
\end{eqnarray*}%
Recall that 
\begin{equation*}
\lbrack e_{2},e_{1}]=\omega (e_{1})e_{1}+2T,
\end{equation*}%
and hence at $\Sigma $, 
\begin{equation*}
e_{2}e_{1}e_{1}(u)=2e_{1}T(u)+2Te_{1}(u)=-2e_{1}(\alpha )+4\alpha ^{2}.
\end{equation*}

%TCIMACRO{\TeXButton{End Proof}{\endproof}}%
%BeginExpansion
\endproof%
%EndExpansion

\bigskip

\begin{proposition}
We have 
\begin{equation}
w=-\frac{2}{3}[e_{1}(\alpha )+2\alpha ^{2}+\frac{1}{6}H^{2}+\frac{3}{16}R-%
\frac{1}{2}a_{2}].  \label{w}
\end{equation}
\end{proposition}

\bigskip

%TCIMACRO{\TeXButton{Proof}{\proof} }%
%BeginExpansion
\proof
%EndExpansion
$w(x)$ will be determined by 
\begin{equation*}
\frac{d^{2}}{d\rho ^{2}}|_{\rho =0}[|d_{P}u|^{2}-\frac{1}{2}u\triangle _{P}u-%
\frac{1}{8}Ru^{2}]=0.
\end{equation*}%
We compute 
\begin{eqnarray*}
\frac{d^{2}}{d\rho ^{2}}(|d_{P}u|^{2}) &=&\frac{d^{2}}{d\rho ^{2}}%
[e_{1}(u)^{2}+e_{2}(u)^{2}] \\
&=&2e_{1}(u)e_{2}e_{2}e_{1}(u)+2(e_{2}e_{1}(u))^{2}+2(u_{\rho \rho
})^{2}+2u_{\rho }u_{\rho \rho \rho }.
\end{eqnarray*}%
Then 
\begin{equation*}
\frac{d^{2}}{d\rho ^{2}}|_{\rho
=0}(|d_{P}u|^{2})=2(e_{2}e_{1}(u))^{2}+2(u_{\rho \rho })^{2}+2u_{\rho \rho
\rho }=8\alpha ^{2}+8v^{2}+12w.
\end{equation*}%
We also have 
\begin{equation*}
\frac{d^{2}}{d\rho ^{2}}(u\triangle _{P}u)=u_{\rho \rho }\triangle
_{P}u+2u_{\rho }e_{2}(\triangle _{P}u)+ue_{2}e_{2}(\triangle _{P}u),
\end{equation*}%
and hence 
\begin{equation*}
\frac{d^{2}}{d\rho ^{2}}|_{\rho =0}(u\triangle
_{P}u)=16v^{2}+2e_{2}(\triangle _{P}u).
\end{equation*}%
Now it follows from 
\begin{equation*}
\triangle _{P}u=e_{1}e_{1}(u)+e_{2}e_{2}(u)-\omega (e_{1})e_{2}(u)
\end{equation*}%
that 
\begin{equation*}
e_{2}(\triangle _{P}u)=e_{2}e_{1}e_{1}(u)+u_{\rho \rho \rho }-e_{2}(\omega
(e_{1}))u_{\rho }-\omega (e_{1})u_{\rho \rho }.
\end{equation*}%
Therefore by (\ref{lemmaofwz}) we get 
\begin{eqnarray*}
\frac{d}{d\rho }|_{\rho =0}(\triangle _{P}u) &=&e_{2}e_{1}e_{1}(u)+u_{\rho
\rho \rho }-e_{2}(\omega (e_{1}))-\omega (e_{1})u_{\rho \rho } \\
&=&-2e_{1}(\alpha )+4\alpha ^{2}+6w-e_{2}(\omega (e_{1}))-2Hv.
\end{eqnarray*}%
By (\ref{Websterscalar}), 
\begin{equation*}
-R=-e_{2}\omega (e_{1})+\omega (e_{1})^{2}+2\omega (T).
\end{equation*}%
Therefore, using (\ref{lemmaofwz}) we get 
\begin{equation}
\frac{d}{d\rho }|_{\rho =0}(\triangle _{P}u)=6w-4e_{1}(\alpha
)-2Hv-R-H^{2}+2a_{2}.  \label{e2laplacianu}
\end{equation}

Now 
\begin{eqnarray*}
0 &=&\frac{d^{2}}{d\rho ^{2}}|_{\rho =0}[|d_{P}u|^{2}-\frac{1}{2}u\triangle
_{P}u-\frac{1}{8}Ru^{2}] \\
&=&8\alpha ^{2}+8v^{2}+12w \\
&&-[8v^{2}+6w-4e_{1}(\alpha )-2Hv-R-H^{2}+2a_{2}] \\
&&-\frac{1}{4}R.
\end{eqnarray*}%
Hence 
\begin{eqnarray*}
0 &=&6w+4e_{1}(\alpha )+8\alpha ^{2}+\frac{2}{3}H^{2}+\frac{3}{4}R-2a_{2}, \\
w &=&-\frac{2}{3}[e_{1}(\alpha )+2\alpha ^{2}+\frac{1}{6}H^{2}+\frac{3}{16}R-%
\frac{1}{2}a_{2}].
\end{eqnarray*}

%TCIMACRO{\TeXButton{End Proof}{\endproof}}%
%BeginExpansion
\endproof%
%EndExpansion

\bigskip

We now start to compute the coefficient $z(x)$.

\bigskip

\begin{lemma}
At $\Sigma $, we have 
\begin{equation}
e_{2}^{2}e_{1}(u):=e_{2}e_{2}e_{1}(u)=2e_{1}(v)+4\alpha v-2a_{1}.
\label{221u}
\end{equation}
\end{lemma}

\bigskip

%TCIMACRO{\TeXButton{Proof}{\proof} }%
%BeginExpansion
\proof
%EndExpansion
(\ref{221u}) follows from (\ref{lemmaofwz}) and the following, at $x\in
\Sigma $, 
\begin{eqnarray*}
\lefteqn{e_{2}e_{2}e_{1}(u)=e_{2}e_{1}e_{2}(u)+e_{2}[e_{2},e_{1}](u)} \\
&=&e_{1}e_{2}e_{2}(u)+[e_{2},e_{1}]e_{2}(u)+e_{2}[e_{2},e_{1}](u) \\
&=&2e_{1}(v)+(He_{1}+2T)e_{2}(u)+e_{2}(\omega (e_{1})e_{1}+2T)(u) \\
&=&2e_{1}(v)+2Te_{2}(u)+He_{2}e_{1}(u)+2e_{2}T(u).
\end{eqnarray*}

%TCIMACRO{\TeXButton{End Proof}{\endproof}}%
%BeginExpansion
\endproof%
%EndExpansion

\bigskip

We have 
\begin{equation*}
\frac{d^{3}}{d\rho ^{3}}%
(|d_{P}u|^{2})=2e_{1}(u)e_{2}^{3}e_{1}(u)+6e_{2}e_{1}(u)e_{2}^{2}e_{1}(u)+2e_{2}(u)e_{2}^{4}(u)+6e_{2}^{2}(u)e_{2}^{3}(u).
\end{equation*}%
Hence by (\ref{221u}), 
\begin{eqnarray}
\lefteqn{\frac{d^{3}}{d\rho ^{3}}|_{\rho =0}(|d_{P}u|^{2})=-12\alpha
e_{2}^{2}e_{1}(u)+48z+72vw}  \notag \\
&=&-12\alpha (2e_{1}(v)+4\alpha v-2a_{1})+48z+72vw.  \label{thirdofgradient}
\end{eqnarray}%
Now we compute 
\begin{equation*}
\frac{d^{3}}{d\rho ^{3}}(u\triangle _{P}u)=u_{\rho \rho \rho }\triangle
_{P}u+3u_{\rho \rho }e_{2}\triangle _{P}u+3u_{\rho }e_{2}^{2}\triangle
_{P}u+ue_{2}^{3}\triangle _{P}u.
\end{equation*}%
Using (\ref{laplacianu}) and (\ref{e2laplacianu}), we find at $\Sigma $, 
\begin{eqnarray}
\frac{d^{3}}{d\rho ^{3}}|_{\rho =0}(u\triangle _{P}u)
&=&48vw+6v(6w-4e_{1}(\alpha )-24v^{2}-R+2a_{2})  \notag \\
&&+3e_{2}^{2}\triangle _{P}u.  \label{222ulaplacianu0}
\end{eqnarray}

\bigskip

\begin{lemma}
On $\Sigma$, we have 
\begin{eqnarray*}
\frac{d^3}{d\rho^3}|_{\rho=0}(u\triangle_Pu)&=&72z+6e_1e_1(v)-24V(%
\alpha)-3e_2(R)-6Rv \\
&&-12\alpha e_1(v)+16ve_1(\alpha)-48v^3-208\alpha^2v \\
&&+36a_1\alpha+28a_2v-12e_1(a_1)+6e_2(a_2).
\end{eqnarray*}
\end{lemma}

\bigskip

%TCIMACRO{\TeXButton{Proof}{\proof} }%
%BeginExpansion
\proof
%EndExpansion
The term we need to handle is $e_{2}^{2}\triangle _{P}u$. We compute 
\begin{eqnarray}
e_{2}^{2}\triangle _{P}u &=&e_{2}e_{2}(e_{1}e_{1}u+e_{2}e_{2}u-\omega
(e_{1})e_{2}(u))  \notag \\
&=&u_{\rho \rho \rho \rho }-Hu_{\rho \rho \rho }-2e_{2}(\omega
(e_{1}))u_{\rho \rho }-e_{2}e_{2}\omega (e_{1})+e_{2}e_{2}e_{1}e_{1}(u).
\label{z1}
\end{eqnarray}%
Recall (\ref{Websterscalar}) 
\begin{equation*}
-R=-e_{2}\omega (e_{1})+\omega (e_{1})^{2}+2\omega (T).
\end{equation*}%
Then 
\begin{eqnarray*}
e_{2}e_{2}(\omega (e_{1})) &=&e_{2}(R)+2\omega (e_{1})e_{2}(\omega
(e_{1}))+2e_{2}(\omega (T)) \\
&=&e_{2}(R)+2H(R+H^{2}+2\omega (T))+2e_{2}(\omega (T)).
\end{eqnarray*}%
Using (\ref{domega}) and (\ref{D12}), we have 
\begin{equation*}
d\omega (e_{2},T)=D_{2}=e_{1}(a_{1})-e_{2}(a_{2})+2a_{2}\omega (e_{1}).
\end{equation*}%
On the other hand, using (\ref{commutators}) we get 
\begin{eqnarray*}
d\omega (e_{2},T) &=&e_{2}(\omega (T))-\omega ([e_{2},T]) \\
&=&e_{2}(\omega (T))-\omega (e_{1})\omega (T)+a_{2}\omega (e_{1}),
\end{eqnarray*}%
and hence from (\ref{lemmaofwz}) we get 
\begin{equation}
e_{2}(\omega (T))=\omega (e_{1})(e_{1}(\alpha )+2\alpha
^{2})+e_{1}(a_{1})-e_{2}(a_{2}),  \label{e2omegaT}
\end{equation}%
\begin{equation*}
e_{2}e_{2}(\omega (e_{1}))=e_{2}(R)+2H(R+H^{2}+3e_{1}(\alpha )+6\alpha
^{2})-4Ha_{2}+2e_{1}(a_{1})-2e_{2}(a_{2}),
\end{equation*}%
\begin{eqnarray}
e_{2}^{2}\triangle _{P}u &=&24z-6Hw-4v(R+H^{2}+2e_{1}(\alpha )+4\alpha
^{2}-2a_{2})  \notag \\
&&-[e_{2}(R)+2H(R+H^{2}+3e_{1}(\alpha )+6\alpha
^{2})-4Ha_{2}+2e_{1}(a_{1})-2e_{2}(a_{2})]  \label{z2} \\
&&+e_{2}e_{2}e_{1}e_{1}(u).  \notag
\end{eqnarray}

Now we compute 
\begin{eqnarray*}
\lefteqn{%
e_{2}e_{2}e_{1}e_{1}(u)=e_{2}e_{1}e_{2}e_{1}(u)+e_{2}[e_{2},e_{1}]e_{1}(u)}
\\
&=&e_{1}e_{2}e_{2}e_{1}(u)+[e_{2},e_{1}]e_{2}e_{1}(u)+e_{2}[e_{2},e_{1}]e_{1}(u)
\\
&=&e_{1}e_{2}e_{1}e_{2}(u)+e_{1}e_{2}[e_{2},e_{1}](u)+[e_{2},e_{1}]e_{2}e_{1}(u)+e_{2}[e_{2},e_{1}]e_{1}(u)
\\
&=&e_{1}e_{1}e_{2}e_{2}(u)+e_{1}[e_{2},e_{1}]e_{2}(u)+e_{1}e_{2}[e_{2},e_{1}](u)
\\
&&+2[e_{2},e_{1}]e_{2}e_{1}(u)+[e_{2},[e_{2},e_{1}]]e_{1}(u) \\
:= &&I+II+III+IV+V.
\end{eqnarray*}%
We have 
\begin{eqnarray*}
I &=&e_{1}e_{1}e_{2}e_{2}(u)=2e_{1}e_{1}(v), \\
II &=&e_{1}[e_{2},e_{1}]e_{2}(u)=e_{1}[(\omega (e_{1})e_{1}+2T)e_{2}(u)] \\
&=&2e_{1}Te_{2}(u)=2e_{1}(V-\alpha e_{2})e_{2}(u) \\
&=&-4e_{1}(\alpha v), \\
III &=&e_{1}e_{2}[e_{2},e_{1}](u)=e_{1}e_{2}(\omega (e_{1})e_{1}+2T)u \\
&=&e_{1}e_{2}(\omega (e_{1})e_{1}(u))+2e_{1}e_{2}T(u) \\
&=&e_{1}(\omega (e_{1}))e_{2}e_{1}(u)+He_{1}e_{2}e_{1}(u)+2e_{1}(-2\alpha
v-a_{1}) \\
&=&-2\alpha e_{1}(H)-2He_{1}(\alpha )-4e_{1}(\alpha v)-2e_{1}(a_{1}) \\
&=&8e_{1}(\alpha v)-2e_{1}(a_{1}).
\end{eqnarray*}%
Using (\ref{221u}), we get 
\begin{eqnarray*}
Te_{2}e_{1}(u) &=&(V-\alpha e_{2})e_{2}e_{1}(u) \\
&=&V(-2\alpha )-\alpha (2e_{1}(v)+4\alpha v-2a_{1}) \\
&=&-2V(\alpha )-2\alpha e_{1}(v)-4\alpha ^{2}v+2a_{1}\alpha .
\end{eqnarray*}%
Then 
\begin{eqnarray*}
\frac{IV}{2} &=&[e_{2},e_{1}]e_{2}e_{1}(u)=(\omega
(e_{1})e_{1}+2T)e_{2}e_{1}(u) \\
&=&-2He_{1}(\alpha )+2Te_{2}e_{1}(u), \\
&=&-4V(\alpha )-4\alpha e_{1}(v)+12ve_{1}(\alpha )-8\alpha
^{2}v+4a_{1}\alpha .
\end{eqnarray*}%
Finally, 
\begin{eqnarray*}
V &=&[e_{2},[e_{2},e_{1}]]e_{1}(u) \\
&=&[e_{2},\omega (e_{1})e_{1}+2T]e_{1}(u) \\
&=&H[e_{2},e_{1}]e_{1}(u)-2a_{1}e_{2}e_{1}(u) \\
&=&2HTe_{1}(u)+4a_{1}\alpha \\
&=&-24\alpha ^{2}v+4a_{1}\alpha .
\end{eqnarray*}%
So we find 
\begin{equation*}
e_{2}e_{2}e_{1}e_{1}(u)=2e_{1}e_{1}(v)-8V(\alpha )-4\alpha
e_{1}(v)+28ve_{1}(\alpha )-40\alpha ^{2}v+12a_{1}\alpha -2e_{1}(a_{1}),
\end{equation*}%
and it follows from (\ref{z2}) that 
\begin{eqnarray*}
e_{2}^{2}\triangle _{P}u &=&24z-6Hw-4v(R+H^{2}+2e_{1}(\alpha )+4\alpha
^{2}-2a_{2}) \\
&&-[e_{2}(R)+2H(R+H^{2}+3e_{1}(\alpha )+6\alpha
^{2})-4Ha_{2}+2e_{1}(a_{1})-2e_{2}(a_{2})] \\
&&+2e_{1}e_{1}(v)-8V(\alpha )-4\alpha e_{1}(v)+28ve_{1}(\alpha )-40\alpha
^{2}v+12a_{1}\alpha -2e_{1}(a_{1}) \\
&=&24z-6Hw-4v(R+H^{2}+2e_{1}(\alpha )+4\alpha ^{2}) \\
&&-[e_{2}(R)+2H(R+H^{2}+3e_{1}(\alpha )+6\alpha ^{2})] \\
&&+2e_{1}e_{1}(v)-8V(\alpha )-4\alpha e_{1}(v)+28ve_{1}(\alpha )-40\alpha
^{2}v \\
&&+12a_{1}\alpha -16a_{2}v-4e_{1}(a_{1})+2e_{2}(a_{2}).
\end{eqnarray*}%
Hence by (\ref{222ulaplacianu0}) and (\ref{w}), we have 
\begin{eqnarray}
\frac{d^{3}}{d\rho ^{3}}|_{\rho =0}(u\triangle _{P}u)
&=&48vw+6v(6w-4e_{1}(\alpha )-24v^{2}-R+2a_{2})+3e_{2}^{2}\triangle _{P}u 
\notag \\
&=&72z+6e_{1}e_{1}(v)-24V(\alpha )-3e_{2}(R)-6Rv  \notag \\
&&-12\alpha e_{1}(v)+16ve_{1}(\alpha )-48v^{3}-208\alpha ^{2}v  \label{z3} \\
&&+36a_{1}\alpha +28a_{2}v-12e_{1}(a_{1})+6e_{2}(a_{2}).  \notag
\end{eqnarray}

%TCIMACRO{\TeXButton{End Proof}{\endproof}}%
%BeginExpansion
\endproof%
%EndExpansion

\bigskip

\begin{proposition}
\begin{eqnarray}
12z &=&3e_{1}e_{1}(v)-12V(\alpha )-\frac{3}{4}e_{2}(R)+\frac{15}{2}Rv
\label{z} \\
&&+18\alpha e_{1}(v)+56ve_{1}(\alpha )+264v^{3}+40\alpha ^{2}v  \notag \\
&&-6a_{1}\alpha -10a_{2}v-6e_{1}(a_{1})+3e_{2}(a_{2}).  \notag
\end{eqnarray}
\end{proposition}

\bigskip

%TCIMACRO{\TeXButton{Proof}{\proof} }%
%BeginExpansion
\proof
%EndExpansion
By (\ref{thirdofgradient}) and (\ref{z3}), we have%
\begin{eqnarray*}
0 &=&\frac{d^{3}}{d\rho ^{3}}|_{\rho =0}[|d_{P}u|^{2}-\frac{1}{2}u\triangle
_{P}u-\frac{1}{8}Ru^{2}] \\
&=&-12\alpha (2e_{1}(v)+4\alpha v-2a_{1})+48z+72vw \\
&&-\frac{1}{2}[72z+6e_{1}e_{1}(v)-24V(\alpha )-3e_{2}(R)-6Rv \\
&&-12\alpha e_{1}(v)+16ve_{1}(\alpha )-48v^{3}-208\alpha ^{2}v \\
&&+36a_{1}\alpha +28a_{2}v-12e_{1}(a_{1})+6e_{2}(a_{2})] \\
&&-\frac{3}{4}(e_{2}(R)+2Rv).
\end{eqnarray*}%
Then (\ref{z}) follows from (\ref{v}) and (\ref{w}).

%TCIMACRO{\TeXButton{End Proof}{\endproof}}%
%BeginExpansion
\endproof%
%EndExpansion

\bigskip

\section{Volume renormalization and $E_{2}$}

In this section we will find that $E_{2}$ can be interpreted as the
coefficient of the log term in the volume renormalization for a formal
solution to the singular CR Yamabe problem. We also prove that the
coefficient of the first log term in the expansion of a formal solution is
given by the variational derivative of $E_{2}$. This is an analogue of the
result in singular Yamabe problem (\cite{GrahamYamabe}). We follow Graham's
ideas (\cite{GrahamYamabe}) in the conformal case.

\bigskip

\subsection{Expansions of tangent vector fields of $\Sigma_\protect\rho$}

In order to carry out the volume renormalization for a formal solution $%
\widetilde{\theta }=u^{-2}\theta $, we need also to expand the volume form $%
\theta \wedge d\theta $ in $\rho $. More precisely, for $\widetilde{\theta }%
=u^{-2}\theta $ we have $\widetilde{\theta }\wedge d\widetilde{\theta }%
=u^{-4}\theta \wedge d\theta $ and 
\begin{equation*}
Vol(\{\rho >\epsilon \})=\int_{\{\rho >\epsilon \}}u^{-4}\langle e_{2},\nu
\rangle d\mu _{\Sigma _{\rho }}d\rho ,
\end{equation*}%
where $\nu $ is the unit inner normal of $\Sigma _{\rho }$. In this
subsection we compute the expansion of $\langle e_{2},\nu \rangle d\mu
_{\Sigma _{\rho }}$ in $\rho $.

Consider a family of geodesics 
\begin{equation*}
\gamma _{s}(\rho )=\gamma (s,\rho )=\exp _{x(s)}(\rho e_{2}),\quad x(s)\in
\Sigma ,
\end{equation*}%
and let 
\begin{equation*}
U=\frac{\partial }{\partial s}\gamma (s,\rho ),\quad e_{2}=\frac{\partial }{%
\partial \rho }\gamma (s,\rho ).
\end{equation*}%
By definition $U\in T_{y}\Sigma _{\rho }$. Note that $[U,e_{2}]=d\gamma
([\partial _{s},\partial _{\rho }])=0$, and if along $\gamma _{0}(\rho )$, 
\begin{equation*}
U(\rho )=f(\rho )e_{1}+g(\rho )e_{2}+h(\rho )T,
\end{equation*}%
then 
\begin{eqnarray*}
0 &=&[e_{2},fe_{1}+ge_{2}+hT] \\
&=&e_{2}(f)e_{1}+e_{2}(g)e_{2}+e_{2}(h)T+f\omega (e_{1})e_{1}+2fT+h\omega
(T)e_{1}-h(a_{2}e_{1}+a_{1}e_{2}),
\end{eqnarray*}%
that is 
\begin{equation}
\left\{ 
\begin{array}{lll}
e_{2}(f)+f\omega (e_{1})+h\omega (T)-ha_{2}=0, &  &  \\ 
e_{2}(g)-ha_{1}=0, &  &  \\ 
\label{odeofvectorfield}e_{2}(h)+2f=0. &  & 
\end{array}%
\right.
\end{equation}%
Differentiating (\ref{odeofvectorfield}) with respect to $e_{2}$ and using (%
\ref{odeofvectorfield}) and (\ref{e2omegaT}), we get 
\begin{eqnarray*}
0 &=&e_{2}e_{2}(f)+e_{2}(f)\omega (e_{1})+fe_{2}(\omega (e_{1}))-2f\omega
(T)+he_{2}(\omega (T))+2fa_{2}-he_{2}(a_{2}) \\
&=&e_{2}e_{2}(f)-\omega (e_{1})(f\omega (e_{1})+h\omega (T)-ha_{2}) \\
&&+f(R+\omega (e_{1})^{2}+2\omega (T))-2f\omega (T)+he_{2}(\omega
(T))+2fa_{2}-he_{2}(a_{2}) \\
&=&e_{2}e_{2}(f)+fR+h[e_{2}(\omega (T))-\omega (e_{1})\omega (T)+a_{2}\omega
(e_{1})]+2fa_{2}-he_{2}(a_{2}) \\
&=&e_{2}e_{2}(f)+fR+2fa_{2}+h[e_{1}(a_{1})-2e_{2}(a_{2})+2a_{2}\omega
(e_{1})].
\end{eqnarray*}%
Therefore we have the Jacobi equations for $U$: 
\begin{equation}
\left\{ 
\begin{array}{lll}
e_{2}e_{2}(f)+fR+2fa_{2}+h[e_{1}(a_{1})-2e_{2}(a_{2})+2a_{2}\omega
(e_{1})]=0, &  &  \\ 
e_{2}e_{2}(g)+2fa_{1}-he_{2}(a_{1})=0, &  &  \\ 
e_{2}e_{2}(h)+2e_{2}(f)=0. &  & 
\end{array}%
\right.  \label{Jacobi}
\end{equation}

Let 
\begin{equation*}
U_1(0)=e_1,
\end{equation*}
then 
\begin{equation*}
U_1(\rho)=f(\rho)e_1+g(\rho)e_2+h(\rho)T,
\end{equation*}
satisfies 
\begin{equation*}
f(0)=1, \quad g(0)=0, \quad h(0)=0.
\end{equation*}
It follows from (\ref{odeofvectorfield}) that 
\begin{equation*}
f^{\prime }(0)=-H, \quad g^{\prime }(0)=0, \quad h^{\prime }(0)=-2.
\end{equation*}
Hence by (\ref{Jacobi}) 
\begin{equation*}
f^{\prime \prime }(0)=-R-2a_2,\quad g^{\prime \prime }(0)=-2a_1, \quad
h^{\prime \prime }(0)=2H.
\end{equation*}
Then by differentiating (\ref{Jacobi}), we get 
\begin{eqnarray*}
f^{\prime \prime \prime }(0)&=&HR-e_2(R)+6Ha_2+2e_1(a_1)-6e_2(a_2), \\
g^{\prime \prime \prime }(0)&=&2Ha_1-4e_2(a_1), \\
h^{\prime \prime \prime }(0)&=&2R+4a_2.
\end{eqnarray*}
%So we find
%\begin{eqnarray*}
%U_1(\rho)&=&(1-H\rho-\frac{1}{2}(R+2a_2)\rho^2+\frac{1}{6}[HR-e_2(R)+6Ha_2+2e_1(a_1)-6e_2(a_2)]\rho^3)e_1
%\\&&+(-a_1\rho^2+\frac{1}{3}(Ha_1-2e_2(a_1))\rho^3)e_2
%\\&&+(-2\rho+H\rho^2+\frac{1}{3}(W+2a_2)\rho^3)T+o(\rho^3).
%\end{eqnarray*}

Let 
\begin{equation*}
U_{2}(0)=V=T+\alpha e_{2},
\end{equation*}%
and for $U_{2}$ we have 
\begin{equation*}
f(0)=0,\quad g(0)=\alpha ,\quad h(0)=1,
\end{equation*}%
it follows from (\ref{odeofvectorfield}) that 
\begin{equation*}
f^{\prime }(0)=-\omega (T)+a_{2}=-e_{1}(\alpha )-2\alpha ^{2}+2a_{2},\quad
g^{\prime }(0)=a_{1},\quad h^{\prime }(0)=0.
\end{equation*}%
From (\ref{Jacobi}) we also have 
\begin{eqnarray*}
f^{\prime \prime }(0) &=&-[e_{1}(a_{1})-2e_{2}(a_{2})+2Ha_{2}], \\
g^{\prime \prime }(0) &=&e_{2}(a_{1}), \\
h^{\prime \prime }(0) &=&2e_{1}(\alpha )+4\alpha ^{2}-4a_{2},
\end{eqnarray*}%
and hence 
\begin{eqnarray*}
f^{\prime \prime \prime }(0) &=&(R+2a_{2})(e_{1}(\alpha )+2\alpha
^{2}-2a_{2})-e_{2}[e_{1}(a_{1})-2e_{2}(a_{2})+2\omega (e_{1})a_{2}], \\
g^{\prime \prime \prime }(0) &=&2a_{1}[e_{1}(\alpha )+2\alpha
^{2}-2a_{2}]+e_{2}e_{2}(a_{1}), \\
h^{\prime \prime \prime }(0) &=&2[e_{1}(a_{1})-2e_{2}(a_{2})+2Ha_{2}].
\end{eqnarray*}%
%
%
%
%
%
%
%
%
%
%then
%\begin{eqnarray*}
%U_2(\rho)&=&[-(\omega(T)-a_2)\rho-\frac{1}{2}[e_1(a_1)-2e_2(a_2)+2\omega(e_1)a_2]\rho^2
%\\&&+\frac{1}{6}((W+2a_2)(e_1(\alpha)+2\alpha^2-2a_2)-e_2[e_1(a_1)-2e_2(a_2)+2\omega(e_1)a_2])\rho^3+\cdots]e_1
%\\&&+[\alpha+a_1\rho+\frac{1}{2}e_2(a_1)\rho^2+\frac{1}{6}(2\omega(T)a_1-2a_1a_2+e_2e_2(a_1))\rho^3+\cdots]e_2
%\\&&+[1+(\omega(T)-a_2)\rho^2+\frac{1}{3}[e_1(a_1)-2e_2(a_2)+2\omega(e_1)a_2]\rho^3+\cdots]T.
%\end{eqnarray*}

\bigskip

\begin{lemma}
There holds 
\begin{eqnarray}
\lefteqn{\langle e_{2},\nu \rangle d\mu _{\Sigma _{\rho }}=[1-H\rho
-(e_{1}(\alpha )+2\alpha ^{2}+\frac{1}{2}R-a_{2})\rho ^{2}}  \notag \\
&&+\frac{1}{6}(RH-e_{2}(R)-2Ha_{2}-2e_{1}(a_{1})+2e_{2}(a_{2}))\rho
^{3}+o(\rho ^{3})]\frac{d\mu _{\Sigma }}{\sqrt{1+\alpha ^{2}}}.
\label{expansiondmu}
\end{eqnarray}
\end{lemma}

\bigskip

%TCIMACRO{\TeXButton{Proof}{\proof} }%
%BeginExpansion
\proof
%EndExpansion
For $y=\gamma (x,\rho )$, we write 
\begin{equation}
U_{1}(y)=(d\gamma )_{y}(e_{1}):=f_{1}e_{1}+g_{1}e_{2}+h_{1}T+o(\rho ^{3}),
\label{U1}
\end{equation}%
and 
\begin{equation}
U_{2}=(d\gamma )_{y}(V):=f_{2}e_{1}+g_{2}e_{2}+h_{2}T+o(\rho ^{3}).
\label{U2}
\end{equation}%
We have 
\begin{equation*}
d\mu _{\Sigma _{\rho }}=\frac{\sqrt{|U_{1}|^{2}|U_{2}|^{2}-\langle
U_{1},U_{2}\rangle ^{2}}}{\sqrt{1+\alpha ^{2}}}d\mu _{\Sigma },
\end{equation*}%
and the unit inner normal to $\Sigma _{\rho }$ is given by 
\begin{eqnarray*}
\lefteqn{\nu (x,\rho )=-\frac{U_{1}\wedge U_{2}}{|U_{1}\wedge U_{2}|}} \\
&=&-\frac{%
(g_{1}h_{2}-g_{2}h_{1})e_{1}+(h_{1}f_{2}-h_{2}f_{1})e_{2}+(f_{1}g_{2}-f_{2}g_{1})T+o(\rho ^{3})%
}{\sqrt{|U_{1}|^{2}|U_{2}|^{2}-\langle U_{1},U_{2}\rangle ^{2}}}.
\end{eqnarray*}%
Then 
\begin{equation*}
\langle e_{2},\nu \rangle =\frac{f_{1}h_{2}-f_{2}h_{1}+o(\rho ^{3})}{\sqrt{%
|U_{1}|^{2}|U_{2}|^{2}-\langle U_{1},U_{2}\rangle ^{2}}},
\end{equation*}%
and 
\begin{equation*}
\langle e_{2},\nu \rangle d\mu _{\Sigma _{\rho
}}=(f_{1}h_{2}-f_{2}h_{1}+o(\rho ^{3}))\frac{d\mu _{\Sigma }}{\sqrt{1+\alpha
^{2}}}.
\end{equation*}%
(\ref{expansiondmu}) then follows from 
\begin{eqnarray*}
\lefteqn{f_{1}h_{2}-f_{2}h_{1}=(1-H\rho -\frac{1}{2}(R+2a_{2})\rho ^{2}+%
\frac{1}{6}[HR-e_{2}(R)+6Ha_{2}+2e_{1}(a_{1})-6e_{2}(a_{2})]\rho ^{3})} \\
&&\cdot (1+(e_{1}(\alpha )+2\alpha ^{2}-2a_{2})\rho ^{2}+\frac{1}{3}%
[e_{1}(a_{1})-2e_{2}(a_{2})+2Ha_{2}]\rho ^{3}) \\
&&+[(e_{1}(\alpha )+2\alpha ^{2}-2a_{2})\rho +\frac{1}{2}%
[e_{1}(a_{1})-2e_{2}(a_{2})+2Ha_{2}]\rho ^{2}](-2\rho +H\rho ^{2})+o(\rho
^{3}) \\
&=&1-H\rho -(e_{1}(\alpha )+2\alpha ^{2}+\frac{1}{2}R-a_{2})\rho ^{2} \\
&&+\frac{1}{6}(RH-e_{2}(R)-2Ha_{2}-2e_{1}(a_{1})+2e_{2}(a_{2}))\rho
^{3}+o(\rho ^{3}).
\end{eqnarray*}

%TCIMACRO{\TeXButton{End Proof}{\endproof}}%
%BeginExpansion
\endproof%
%EndExpansion

\bigskip

\subsection{Volume renormalizatioin}

In this subsection we carry out the volume renormalization for a formal
solution $\widetilde{\theta }=u^{-2}\theta $ to the singular CR Yamabe
problem. We are going to obtain the coefficient of $\log \frac{1}{\epsilon }$
in the volume renormalizatioin.

\bigskip

\begin{proposition}
For $\epsilon >0$ sufficiently small, we have 
\begin{equation}
Vol(\{\rho >\epsilon \})=c_{0}\epsilon ^{-3}+c_{1}\epsilon
^{-2}+c_{2}\epsilon ^{-1}+L\log \frac{1}{\epsilon }+V_{0}+o(1),
\label{volumerenormalization}
\end{equation}%
where 
\begin{eqnarray}
c_{0} &=&\frac{1}{3}\int_{\Sigma }\theta \wedge e^{1},  \label{c0c1c2L} \\
c_{1} &=&-\frac{1}{6}\int_{\Sigma }H\theta \wedge e^{1},  \notag \\
c_{2} &=&\frac{1}{3}\int_{\Sigma }(5e_{1}(\alpha )+10\alpha ^{2}+\frac{1}{6}%
H^{2}-a_{2})\theta \wedge e^{1},  \notag \\
L &=&\int_{\Sigma }(\frac{1}{6}e_{1}e_{1}(H)+4V(\alpha )+\alpha
e_{1}(H)+2He_{1}(\alpha )+\frac{4}{27}H^{3}+\frac{1}{3}WH  \notag \\
&&+\frac{1}{6}e_{2}(W)+2a_{1}\alpha +\frac{5}{3}e_{1}(a_{1})-\frac{2}{3}%
e_{2}(a_{2}))\theta \wedge e^{1}.  \notag
\end{eqnarray}
\end{proposition}

\bigskip

%TCIMACRO{\TeXButton{Proof}{\proof} }%
%BeginExpansion
\proof
%EndExpansion
For $\widetilde{\theta }=u^{-2}\theta $, we have $\widetilde{\theta }\wedge d%
\widetilde{\theta }=u^{-4}\theta \wedge d\theta $ and 
\begin{equation}
Vol(\{\rho >\epsilon \})=\int_{\{\rho >\epsilon \}}u^{-4}\langle e_{2},\nu
\rangle d\mu _{\Sigma _{\rho }}d\rho .  \label{volumeintegral}
\end{equation}
It follows from 
\begin{equation*}
u(x,\rho )=\rho +v(x)\rho ^{2}+w(x)\rho ^{3}+z(x)\rho ^{4}+o(\rho ^{4})
\end{equation*}%
that 
\begin{equation*}
u^{-4}=\rho ^{-4}(1-4v\rho -(4w-10v^{2})\rho ^{2}-(4z-20vw+20v^{3})\rho
^{3}+o(\rho ^{3})).
\end{equation*}%
Then by (\ref{expansiondmu}), we have 
\begin{eqnarray*}
\lefteqn{u^{-4}\langle e_{2},\nu \rangle d\mu _{\Sigma _{\rho }}=\rho
^{-4}(1-4v\rho -(4w-10v^{2})\rho ^{2}-(4z-20vw+20v^{3})\rho ^{3}+o(\rho
^{3}))} \\
&&\cdot \lbrack 1-H\rho -(e_{1}(\alpha )+2\alpha ^{2}+\frac{1}{2}%
R-a_{2})\rho ^{2} \\
&&+\frac{1}{6}(HR-e_{2}(R)-2Ha_{2}-2e_{1}(a_{1})+2e_{2}(a_{2}))\rho
^{3}+o(\rho ^{3})]\frac{d\mu _{\Sigma }}{\sqrt{1+\alpha ^{2}}} \\
&=&\rho ^{-4}[1+2v\rho -(4w+14v^{2}+e_{1}(\alpha )+2\alpha ^{2}+\frac{1}{2}%
R-a_{2})\rho ^{2} \\
&&+(-4z-4vw+40v^{3}+4ve_{1}(\alpha )+8\alpha ^{2}v+Rv-\frac{1}{6}e_{2}(R) \\
&&-2a_{2}v-\frac{1}{3}e_{1}(a_{1})+\frac{1}{3}e_{2}(a_{2}))\rho ^{3}+o(\rho
^{3})]\frac{d\mu _{\Sigma }}{\sqrt{1+\alpha ^{2}}}.
\end{eqnarray*}%
Substituting the formulas of $v,w$ and $z$ in Subsection 5.3 into the above
expression, we get 
\begin{equation}
u^{-4}\langle e_{2},\nu \rangle d\mu _{\Sigma _{\rho }}=\rho
^{-4}[1+v^{(1)}\rho +v^{(2)}\rho ^{2}+v^{(3)}\rho ^{3}+o(\rho ^{3})]\frac{%
d\mu _{\Sigma }}{\sqrt{1+\alpha ^{2}}},  \label{volumeelementexpansion}
\end{equation}
\noindent where 
\begin{eqnarray*}
v^{(1)} &=&-\frac{1}{3}H, \\
v^{(2)} &=&\frac{1}{3}(5e_{1}(\alpha )+10\alpha ^{2}+\frac{1}{6}H^{2}-a_{2})
\\
v^{(3)} &=&\frac{1}{6}e_{1}e_{1}(H)+4V(\alpha )+\alpha
e_{1}(H)+2He_{1}(\alpha )+\frac{4}{27}H^{3}+\frac{1}{3}WH \\
&&+\frac{1}{6}e_{2}(W)+2a_{1}\alpha +\frac{5}{3}e_{1}(a_{1})-\frac{2}{3}%
e_{2}(a_{2}),
\end{eqnarray*}%
where in the last line we have replaced $R$ by $2W$. Then (\ref%
{volumerenormalization}) follows from (\ref{volumeintegral}), (\ref%
{volumeelementexpansion}) and the fact that $\frac{d\mu _{\Sigma }}{\sqrt{%
1+\alpha ^{2}}}=\theta \wedge e^{1}$. In particular, 
\begin{equation}
L(\Sigma )=\int_{\Sigma }v^{(3)}\theta \wedge e^{1}.  \label{Lv3}
\end{equation}

%TCIMACRO{\TeXButton{End Proof}{\endproof}}%
%BeginExpansion
\endproof%
%EndExpansion

\bigskip

\subsection{Conformal invariance of $L(\Sigma)$}

In this subsection we consider the transformation law of $l(x)$ and $%
L(\Sigma)$ under a conformal change of $\theta$. We will show that $l(x)$ is
conformally covariant and $L(\Sigma)$ is a conformal invariant.

Let $\widehat{\theta }=U^{2}\theta $, where $U$ is a smooth positive
function defined on $M$, then at $\Sigma $, 
\begin{equation*}
\widehat{e}_{1}=U^{-1}e_{1},\quad \widehat{e}_{2}=U^{-1}e_{2}.
\end{equation*}%
The defining function $\widehat{\rho }$ with respect to $\widehat{\theta }$
satisfies $\widehat{e}_{2}(\widehat{\rho })|_{\Sigma }=1$, so that 
\begin{equation*}
\widehat{\rho }=U\rho +O(\rho ^{2}).
\end{equation*}%
Note that $\widehat{u}:=uU$ is a formal solution to 
\begin{equation*}
R_{\widehat{u}^{-2}\widehat{\theta }}=-8+O(\widehat{\rho }^{5}\log\widehat{%
\rho }),
\end{equation*}

\noindent and hence 
\begin{eqnarray*}
\widehat{u} &=&uU=(\rho +v\rho ^{2}+w\rho ^{3}+z\rho ^{4}+h\rho ^{5}+l\rho
^{5}\log \rho +o(\rho ^{5}))U \\
&=&\widehat{\rho }+\widehat{v}\widehat{\rho }^{2}+\widehat{w}\widehat{\rho }%
^{3}+\widehat{z}\widehat{\rho }^{4}+\widehat{h}\widehat{\rho }^{5}+\widehat{l%
}\widehat{\rho }^{5}\log \widehat{\rho }+o(\widehat{\rho }^{5}).
\end{eqnarray*}%
Note that $\widehat{\rho }$ is smooth up to the boundary in $\rho $, so we
find 
\begin{equation*}
\widehat{l}=(U|_{\Sigma })^{-4}l.
\end{equation*}

\bigskip

\begin{proposition}
$L(\Sigma)$ is invariant under conformal transformations of $\theta$.
\end{proposition}

\bigskip

%TCIMACRO{\TeXButton{Proof}{\proof} }%
%BeginExpansion
\proof
%EndExpansion
For $\widehat{\theta }=U^{2}\theta $, we have 
\begin{equation*}
\widehat{\rho }(x,\rho )=e^{\Upsilon (x,\rho )}\rho
\end{equation*}%
for some smooth function $\Upsilon $. We can solve the relation $\widehat{%
\rho }=e^{\Upsilon }\rho $ to get 
\begin{equation*}
\rho (x,\widehat{\rho })=b(x,\widehat{\rho })\widehat{\rho },
\end{equation*}%
where $b(x,\widehat{\rho })$ is a smooth function which is bounded from
below by a positive constant. Let 
\begin{equation*}
\widehat{\epsilon }(x,\epsilon )=\epsilon b(x,\epsilon ),
\end{equation*}%
so that $\widehat{\rho }>\epsilon $ is equivalent to $\rho >\widehat{%
\epsilon }(x,\epsilon )$. It follows from (\ref{volumeintegral}) and (\ref%
{volumeelementexpansion}) that 
\begin{eqnarray*}
\lefteqn{Vol(\{\rho >\epsilon \})-Vol(\{\widehat{\rho }>\epsilon
\})=\int_{\{\rho >\epsilon \}}d\mu _{u^{-2}\theta }-\int_{\{\widehat{\rho }%
>\epsilon \}}d\mu _{u^{-2}\theta }} \\
&=&\int_{\Sigma }\int_{\epsilon }^{\widehat{\epsilon }(x,\epsilon
)}[\sum_{k=0}^{3}\rho ^{-4+k}v^{(k)}+o(\rho ^{-1})]d\rho \theta \wedge e^{1}
\\
&=&\sum_{k=0}^{2}\epsilon ^{k-3}\int_{\Sigma }\frac{v^{(k)}}{k-3}%
(b(x,\epsilon )^{k-3}-1)\theta \wedge e^{1}+\int_{\Sigma }v^{(3)}\log
b(x,\epsilon )\theta \wedge e^{1}+o(1).
\end{eqnarray*}%
As $\epsilon \rightarrow 0$, the above contains no $\log \frac{1}{\epsilon }$
term. Hence $L(\Sigma ,\theta )=L(\Sigma ,\widehat{\theta })$.

%TCIMACRO{\TeXButton{End Proof}{\endproof}}%
%BeginExpansion
\endproof%
%EndExpansion

\bigskip

\subsection{Comparison between $L$ and $E_2$}

We are going to show that the difference $dA_{2}-\frac{1}{2}v^{(3)}\theta
\wedge e^{1}$ is an exact form. In particular, if $\Sigma $ is a closed,
nonsingular surface, $E_{2}(\Sigma )=\frac{1}{2}L(\Sigma )$. Recall that 
\begin{eqnarray*}
E_{2}(\Sigma ) &=&\int_{\Sigma }dA_{2} \\
&=&\int_{\Sigma }\{V(\alpha )+\frac{2}{3}[e_{1}(\alpha )+\frac{1}{2}\alpha
^{2}-\func{Im}A_{11}+\frac{1}{4}W]H+\frac{2}{27}H^{3} \\
&&+\func{Im}[\frac{1}{6}W^{,1}+\frac{2i}{3}(A^{11})_{,1}]-\alpha (\func{Re}%
A_{\bar{1}}^{1})\}\theta \wedge e^{1}.
\end{eqnarray*}%
Recall that 
\begin{eqnarray*}
A_{1}^{\overline{1}} &=&A_{11}=a_{1}+ia_{2}, \\
\tau (e_{1}) &=&a_{1}e_{1}-a_{2}e_{2},\quad \tau
(e_{2})=-a_{2}e_{1}-a_{1}e_{2}.
\end{eqnarray*}

\noindent Write

\begin{equation*}
\tau =\tau _{i}^{j}e^{i}\otimes e_{j},\quad \tau _{i}^{j}=\langle \tau
(e_{i}),e_{j}\rangle .
\end{equation*}%
Then 
\begin{eqnarray*}
e_{1}(a_{1}) &=&e_{1}\langle \tau (e_{1}),e_{1}\rangle \\
&=&\langle (\nabla _{e_{1}}\tau )(e_{1})+\tau (\nabla
_{e_{1}}e_{1}),e_{1}\rangle +\langle \tau (e_{1}),\nabla _{e_{1}}e_{1}\rangle
\\
&=&\nabla _{e_{1}}\tau _{1}^{1}-2Ha_{2},
\end{eqnarray*}%
\begin{eqnarray*}
-e_{2}(a_{2}) &=&e_{2}\langle \tau (e_{1}),e_{2}\rangle \\
&=&\langle (\nabla _{e_{2}}\tau )(e_{1})+\tau (\nabla
_{e_{2}}e_{1}),e_{2}\rangle +\langle \tau (e_{1}),\nabla _{e_{2}}e_{2}\rangle
\\
&=&\nabla _{e_{2}}\tau _{1}^{2}-2\omega (e_{2})a_{1}=\nabla _{e_{2}}\tau
_{1}^{2}.
\end{eqnarray*}%
So we have 
\begin{equation*}
A_{\overline{1}}^{1}=A^{11}=a_{1}-ia_{2}=\tau _{1}^{1}+i\tau _{1}^{2},
\end{equation*}%
and 
\begin{eqnarray*}
Re[(A^{11})_{,1}] &=&Re(\nabla _{\frac{1}{2}(e_{1}-ie_{2})}A^{11})=\frac{1}{2%
}(\nabla _{e_{1}}\tau _{1}^{1}+\nabla _{e_{2}}\tau _{1}^{2}) \\
&=&\frac{1}{2}[e_{1}(a_{1})+2Ha_{2}-e_{2}(a_{2})].
\end{eqnarray*}%
So in the the frame given by (\ref{geodesicframe}), we have 
\begin{equation*}
E_{2}(\Sigma )=\int_{\Sigma }dA_{2},
\end{equation*}%
where 
\begin{eqnarray}
dA_{2} &=&\{V(\alpha )+\frac{2}{3}[e_{1}(\alpha )+\frac{1}{2}\alpha ^{2}+%
\frac{1}{4}W]H+\frac{2}{27}H^{3}  \notag \\
&&+\frac{1}{12}e_{2}(W)+\frac{1}{3}[e_{1}(a_{1})-e_{2}(a_{2})]-a_{1}\alpha
\}\theta \wedge e^{1}.  \label{E2measure}
\end{eqnarray}%
On the other hand, $L(\Sigma )=\int_{\Sigma }v^{(3)}\theta \wedge e^{1}$ and 
\begin{eqnarray}
v^{(3)}\theta \wedge e^{1} &=&(\frac{1}{6}e_{1}e_{1}(H)+4V(\alpha )+\alpha
e_{1}(H)+2He_{1}(\alpha )+\frac{4}{27}H^{3}+\frac{1}{3}WH  \notag \\
&&+\frac{1}{6}e_{2}(W)+2a_{1}\alpha +\frac{5}{3}e_{1}(a_{1})-\frac{2}{3}%
e_{2}(a_{2}))\theta \wedge e^{1}.  \label{Lmeasure}
\end{eqnarray}

Note that we have the orthonormal frame of $T\Sigma $: 
\begin{equation*}
e_{1},\quad \widetilde{V}=\frac{V}{\sqrt{1+\alpha ^{2}}}=\frac{T+\alpha e_{2}%
}{\sqrt{1+\alpha ^{2}}},
\end{equation*}%
and the coframe 
\begin{equation*}
e^{1},\quad \sigma =\frac{1}{\sqrt{1+\alpha ^{2}}}(\theta +\alpha e^{2}).
\end{equation*}%
Then on $\Sigma $, we have%
\begin{eqnarray}
d\mu _{\Sigma } &=&\sigma \wedge e^{1}=\sqrt{1+\alpha ^{2}}\theta \wedge
e^{1},  \notag \\
\sigma \wedge \theta &=&e^{2}\wedge \theta =0,  \label{related2form} \\
e^{2}\wedge e^{1} &=&\alpha \theta \wedge e^{1}=\frac{\alpha }{\sqrt{%
1+\alpha ^{2}}}d\mu _{\Sigma }.  \notag
\end{eqnarray}

\bigskip

\begin{lemma}
On $\Sigma$, we have 
\begin{eqnarray}
d(e_1(H)\theta)&=&-[e_1e_1(H)+2\alpha e_1(H)]\theta\wedge e^1,  \notag \\
d(\alpha e^1)&=&[V(\alpha)-\alpha^2H+a_1\alpha]\theta\wedge e^1,  \notag \\
d(\alpha H\theta)&=&-[\alpha e_1(H)+He_1(\alpha)+2\alpha^2 H]\theta\wedge
e^1,  \label{exactforms} \\
d(a_1\theta)&=&-[e_1(a_1)+2a_1\alpha]\theta\wedge e^1.  \notag
\end{eqnarray}
\end{lemma}

\bigskip

%TCIMACRO{\TeXButton{Proof}{\proof} }%
%BeginExpansion
\proof
%EndExpansion
Using (\ref{dtheta}), (\ref{de}) and (\ref{related2form}), we compute 
\begin{eqnarray*}
d(e_{1}(H)\theta ) &=&-e_{1}e_{1}(H)\theta \wedge e^{1}+e_{1}(H)d\theta \\
&=&-e_{1}e_{1}(H)\theta \wedge e^{1}+e_{1}(H)2e^{1}\wedge e^{2} \\
&=&-[e_{1}e_{1}(H)+2\alpha e_{1}(H)]\theta \wedge e^{1}, \\
d(\alpha e^{1}) &=&\widetilde{V}(\alpha )\sigma \wedge e^{1}+\alpha de^{1} \\
&=&V(\alpha )\theta \wedge e^{1}+\alpha (-e^{2}\wedge \omega +\theta \wedge
(a_{1}e^{1}-a_{2}e^{2})) \\
&=&[V(\alpha )-\alpha ^{2}H+a_{1}\alpha ]\theta \wedge e^{1}, \\
d(\alpha H\theta ) &=&-e_{1}(\alpha H)\theta \wedge e^{1}+2\alpha
He^{1}\wedge e^{2} \\
&=&-[\alpha e_{1}(H)+He_{1}(\alpha )+2\alpha ^{2}H]\theta \wedge e^{1}, \\
d(a_{1}\theta ) &=&-e_{1}(a_{1})\theta \wedge e^{1}+2a_{1}e^{1}\wedge e^{2}
\\
&=&-[e_{1}(a_{1})+2a_{1}\alpha ]\theta \wedge e^{1}.
\end{eqnarray*}

%TCIMACRO{\TeXButton{End Proof}{\endproof}}%
%BeginExpansion
\endproof%
%EndExpansion

\bigskip

\begin{proposition}
We have 
\begin{eqnarray}
&&dA_{2}-\frac{1}{2}v^{(3)}\theta \wedge e^{1}  \label{dA2andL} \\
&=&\frac{1}{12}d(e_{1}(H)\theta )-d(\alpha e^{1})+\frac{1}{3}d(\alpha
H\theta)+\frac{1}{2}d(a_{1}\theta ).  \notag
\end{eqnarray}%
In particular, if $\Sigma $ is a closed surface without singular points,
then 
\begin{equation}
E_{2}(\Sigma )=\frac{1}{2}L(\Sigma ).  \label{E2L}
\end{equation}
\end{proposition}

\bigskip

%TCIMACRO{\TeXButton{Proof}{\proof} }%
%BeginExpansion
\proof
%EndExpansion
It follows from (\ref{E2measure}) and $(\ref{Lmeasure})$ that 
\begin{eqnarray*}
dA_{2}-\frac{1}{2}v^{(3)}\theta \wedge e^{1} &=&-\frac{1}{12}%
e_{1}e_{1}(H)-V(\alpha )-\frac{1}{2}\alpha e_{1}(H)-\frac{1}{3}He_{1}(\alpha
) \\
&&+\frac{1}{3}\alpha ^{2}H-\frac{1}{2}e_{1}(a_{1})-2a_{1}\alpha .
\end{eqnarray*}%
(\ref{dA2andL}) then follows from (\ref{exactforms}).

%TCIMACRO{\TeXButton{End Proof}{\endproof}}%
%BeginExpansion
\endproof%
%EndExpansion

\bigskip

\subsection{Smoothness of solutions to the singular Yamabe problem and
critical surfaces of $L(\Sigma)$}

In the expression of a formal solution to the singular CR Yamabe problem, $%
l\neq 0$ is an obstruction to the smoothness of the solution up to the
boundary. In this subsection we show that $l(x)$ is a multiple of $\mathcal{E%
}_2$, see (\ref{mathcalE}). It follows that the solution to the singular CR
Yamabe problem (\ref{CRYamabe}) can be smooth up to the boundary only if $%
\Sigma $ is a critical point of $L(\Sigma )$ or $E_{2}(\Sigma )$.

Let $F_{t}:\Sigma \rightarrow M^{3},t\in (-\delta ,\delta )$ be a variation
of $\Sigma $, such that 
\begin{equation*}
\partial _{t}F_{t}=X=fe_{2}^{t}+gT\in \Gamma (TM|_{\Sigma _{t}})
\end{equation*}%
where $\Sigma _{t}=F_{t}(\Sigma )$, $e_{2}^{t}=Je_{1}^{t}$ and $e_{1}^{t}$
is a unit vector in $\Gamma (T\Sigma _{t}\cap \xi )$. We assume $f,g$ are
supported in a domain of $\Sigma $ away form the singular set of $\Sigma $.
Let $\rho _{t}$ be the defining function for $\Sigma _{t}$, which is defined
by $\rho _{t}(\exp _{\Sigma _{t}}(ae_{2}^{t}))=a$. Here $\gamma (a):=\exp
_{F_{t}(x)}(ae_{2}^{t})$ is the geodesic initiating from $F_{t}(x)$ with $%
\dot{\gamma}(0)=e_{2}^{t}$ such that $\widetilde{e}_{2}^{t,a}:=\dot{\gamma}%
(a)$ satisfies $\nabla _{\widetilde{e}_{2}^{t,a}}\widetilde{e}_{2}^{t,a}=0$.
Let $\Sigma _{t,a}=\exp _{\Sigma _{t}}(ae_{2}^{t})$. For any fixed $t$ and $%
\Sigma _{t}$, we have a formal solution $u_{t}$ with respect to $\Sigma _{t}$%
, and the expansion 
\begin{equation*}
u_{t}^{-4}\langle \widetilde{e}_{2}^{t,\rho _{t}},\nu \rangle d\mu _{\Sigma
_{t,\rho _{t}}}=\rho _{t}^{-4}(1+v_{t}^{(1)}\rho _{t}+v_{t}^{(2)}\rho
_{t}^{2}+v_{t}^{(3)}\rho _{t}^{3}+o(\rho _{t}^{3}))
\end{equation*}%
where $\nu $ is the unit inner normal of $\Sigma _{t,\rho _{t}}$. Let 
\begin{equation*}
L(\Sigma _{t})=\int_{\Sigma _{t}}v_{t}^{(3)}\theta \wedge (e^{1})^{t}.
\end{equation*}

\bigskip

\begin{theorem}
We have 
\begin{equation}
\frac{d}{dt}|_{t=0}L(\Sigma _{t})=10\int_{\Sigma }(f-\alpha g)l\theta \wedge
e^{1},  \label{log}
\end{equation}%
where $l$ is the coefficient of the log term of a formal solution $u$, given
by (\ref{formalsolution}). By (\ref{dA2andL}), (\ref{firstvariationofE2})
and (\ref{mathcalE}), we have%
\begin{equation}
l=\frac{1}{5}\mathcal{E}_{2}.  \label{obstructionlandmathcalE}
\end{equation}
\end{theorem}

\bigskip

%TCIMACRO{\TeXButton{Proof}{\proof} }%
%BeginExpansion
\proof
%EndExpansion
We use an argument which is analogous to \cite{GrahamYamabe}. We have 
\begin{equation*}
u_{t}=\rho _{t}+v_{t}\rho _{t}^{2}+w_{t}\rho _{t}^{3}+z_{t}\rho
_{t}^{4}+h_{t}\rho _{t}^{5}+l_{t}\rho _{t}^{5}\log \rho _{t}+O(\rho _{t}^{6})
\end{equation*}%
such that 
\begin{equation*}
R_{u_{t}^{-2}\theta }=-8+O(\rho _{t}^{5}\log\rho_t).
\end{equation*}%
Differentiating the volume expansion 
\begin{equation*}
Vol_{u_{t}^{-2}\theta }(\{\rho _{t}>\epsilon
\})=\sum_{k=0}^{2}c_{k}^{t}\epsilon ^{-3+k}+L_{t}\log \frac{1}{\epsilon }%
+V_{t}+o(1)
\end{equation*}%
gives 
\begin{equation}
\frac{d}{dt}|_{t=0}Vol_{u_{t}^{-2}\theta }(\{\rho _{t}>\epsilon
\})=\sum_{k=0}^{2}\dot{c}_{k}\epsilon ^{-3+k}+\dot{L}\log \frac{1}{\epsilon }%
+O(1).  \label{dotL1}
\end{equation}%
On the other hand, 
\begin{eqnarray}
\lefteqn{\frac{d}{dt}|_{t=0}Vol_{u_{t}^{-2}\theta }(\{\rho _{t}>\epsilon \})=%
\frac{d}{dt}|_{t=0}\int_{\{\rho _{t}>\epsilon \}}d\mu _{u_{t}^{-2}\theta }} 
\notag \\
&=&\frac{d}{dt}|_{t=0}\int_{\{\rho _{t}>\epsilon \}}d\mu _{u^{-2}\theta
}+\int_{\{\rho >\epsilon \}}\frac{d}{dt}|_{t=0}(u_{t}^{-4})d\mu _{\theta } 
\notag \\
&=&\frac{d}{dt}|_{t=0}\int_{\{\rho _{t}>\epsilon \}}d\mu _{u^{-2}\theta
}-4\int_{\{\rho >\epsilon \}}u^{-1}\dot{u}d\mu _{u^{-2}\theta }.
\label{dotL2}
\end{eqnarray}%
We then compute the coefficient of $\log \frac{1}{\epsilon }$ for the
expansion of (\ref{dotL2}).

We first deal with the term 
\begin{equation*}
\frac{d}{dt}|_{t=0}\int_{\{\rho_t>\epsilon\}}d\mu_{u^{-2}\theta}.
\end{equation*}
One can express $\{\rho_t>\epsilon\}$ as $\{\rho=\rho_0>\psi(x,t,\epsilon),
x\in \Sigma\}$ for a smooth function $\psi(x,t,\epsilon)$. Here $%
\psi(x,t,\epsilon)$ is the $\rho$-coordinate of the point $%
\exp_{F_t(x)}(\epsilon e_2^t)$. Let 
\begin{equation*}
\dot{\psi}(x,\epsilon)=\frac{\partial}{\partial t}\psi(x,0,\epsilon).
\end{equation*}
We have 
\begin{equation}  \label{dotpsi}
\dot{\psi}(x,0)=f(x)-\alpha(x)g(x).
\end{equation}
Actually, 
\begin{equation*}
\dot{\psi}(x,0)=\frac{\partial}{\partial t}\psi(x,0,0)=\frac{d}{dt}%
\rho(F_tx)=(fe_2+gT)\rho,
\end{equation*}
note that 
\begin{equation*}
e_1(\rho)=(T+\alpha e_2)\rho=0,\quad e_2(\rho)=1,
\end{equation*}
hence $T(\rho)=-\alpha$ and (\ref{dotpsi}) follows. For $\epsilon_0>0$ small
and fixed and $\epsilon\ll \epsilon_0$, we have 
\begin{eqnarray*}
\frac{d}{dt}|_{t=0}\int_{\{\rho_t>\epsilon\}}d\mu_{u^{-2}\theta} &=&\frac{d}{%
dt}|_{t=0}\int_\Sigma\int_{\psi(x,t,\epsilon)}^{\epsilon_0}u(x,\rho)^{-4}%
\langle \widetilde{e}_2^{0,\rho},\nu\rangle \frac{d\mu_{\Sigma_{0,\rho}}}{%
d\mu_\Sigma}d\rho d\mu_\Sigma \\
&=&-\int_\Sigma \dot{\psi}(x,\epsilon)u(x,\epsilon)^{-4}\langle \widetilde{e}%
_2^{0,\epsilon},\nu\rangle \frac{d\mu_{\Sigma_{0,\epsilon}}}{d\mu_\Sigma}%
d\mu_\Sigma,
\end{eqnarray*}
where $\nu$ is the unit inner normal of $\Sigma_{0,\epsilon}$. $\dot{\psi}%
(x,\epsilon)\langle \widetilde{e}_2^{0,\epsilon},\nu\rangle \frac{%
d\mu_{\Sigma_{0,\epsilon}}}{d\mu_\Sigma}$ is smooth in $\epsilon$ and it
follows from (\ref{dotpsi}) that 
\begin{equation*}
[\dot{\psi}(x,\epsilon)\langle \widetilde{e}_2^{0,\epsilon},\nu\rangle \frac{%
d\mu_{\Sigma_{0,\epsilon}}}{d\mu_\Sigma}]|_{\epsilon=0} =\frac{f-\alpha g}{%
\sqrt{1+\alpha(x)^2}}.
\end{equation*}
Note that 
\begin{equation*}
u(x,\epsilon)=\epsilon(1+v\epsilon+\cdots+h\epsilon^{4}+l\epsilon^{4}\log
\epsilon+O(\epsilon^5)),
\end{equation*}
therefore the coefficient of $\log\frac{1}{\epsilon}$ in $\frac{d}{dt}%
|_{t=0}\int_{\{\rho_t>\epsilon\}}d\mu_{u^{-2}\theta}$ is 
\begin{equation}  \label{firsttermlog}
-4\int_\Sigma (f-\alpha g)l\theta\wedge e^1.
\end{equation}

We then consider the second term in (\ref{dotL2}), i.e. $-4\int_{\{\rho
>\epsilon \}}u^{-1}\dot{u}d\mu _{u^{-2}\theta }$. For $\Sigma _{t}=\partial
M_{t}$ we have the formal solution to the singular CR Yamabe problem 
\begin{equation*}
u_{t}=\rho _{t}+v_{t}\rho _{t}^{2}+\cdots +h_{t}\rho _{t}^{5}+l_{t}\rho
_{t}^{5}\log \rho _{t}+O(\rho _{t}^{6}),
\end{equation*}%
which satisfies 
\begin{equation}
R_{u_{t}^{-2}\theta }=-8+O(\rho _{t}^{5}\log \rho _{t}).
\label{usingequation0}
\end{equation}%
Let 
\begin{equation*}
\widetilde{\theta }_{t}=u_{t}^{-2}\theta =(\frac{u_{t}}{u})^{-2}\widetilde{%
\theta }_{0}:=e^{2w_{t}}\widetilde{\theta }_{0},
\end{equation*}%
where 
\begin{equation*}
w_{t}=-\log \frac{u_{t}}{u},\quad u_{t}=ue^{-w_{t}}.
\end{equation*}%
Then by (\ref{transformationlawofWS}), 
\begin{equation*}
R_{u_{t}^{-2}\theta }=R_{e^{2w_{t}}\widetilde{\theta }_{0}}=e^{-2w_{t}}[R_{%
\widetilde{\theta }_{0}}-4\triangle _{P}^{\widetilde{\theta }%
_{0}}w_{t}-4|d_{P}w_{t}|_{\widetilde{\theta }_{0}}^{2}].
\end{equation*}%
Differentiating it in $t$ gives 
\begin{equation}
\frac{d}{dt}|_{t=0}R_{u_{t}^{-2}\theta }=-2R_{\widetilde{\theta }_{0}}\dot{w}%
-4\triangle _{P}^{\widetilde{\theta }_{0}}\dot{w}.  \label{usingequation1}
\end{equation}%
It follows from (\ref{usingequation0}) that 
\begin{equation}
\frac{d}{dt}|_{t=0}R_{u_{t}^{-2}\theta }=O(\rho ^{4}\log \rho ).
\label{usingequation2}
\end{equation}%
Note that $u=O(\rho ),\dot{\rho}=O(1)$, so we have 
\begin{eqnarray}
R_{u^{-2}\theta }\dot{w} &=&-8\dot{w}+O(\dot{w}\rho ^{5}\log \rho )
\label{usingequation3} \\
&=&-8\dot{w}+O(\frac{\dot{u}}{u}\rho ^{5}\log \rho )=-8\dot{w}+O(\rho
^{4}\log \rho ).  \notag
\end{eqnarray}%
Therefore it follows from (\ref{usingequation1}), (\ref{usingequation2}) and
(\ref{usingequation3}) that 
\begin{equation*}
4\dot{w}=\triangle _{P}^{\widetilde{\theta }_{0}}\dot{w}+O(\rho ^{4}\log
\rho ).
\end{equation*}%
Hence 
\begin{equation*}
-4\int_{\{\rho >\epsilon \}}u^{-1}\dot{u}d\mu _{u^{-2}\theta }=4\int_{\{\rho
>\epsilon \}}\dot{w}d\mu _{u^{-2}\theta }=\int_{\{\rho >\epsilon
\}}\triangle _{b}^{\widetilde{\theta }_{0}}\dot{w}d\mu _{u^{-2}\theta }+O(1).
\end{equation*}%
Let $\widetilde{\nu }$ be the outward unit normal with respect to $%
\widetilde{\theta }_{0}=u^{-2}\theta $ and $\widetilde{e}_{2},\widetilde{e}%
^{1},\widetilde{\alpha }$ are relative to $(M,\Sigma _{0,\epsilon },%
\widetilde{\theta }_{0})$. Actually, $\widetilde{e}_{2}=ue_{2}$ and $%
e_{2}=Je_{1}$ with $e_{1}\in T\Sigma _{0,\epsilon }\cap \xi $. We also have $%
\widetilde{e}^{1}=u^{-1}e^{1}+u^{-2}e_{2}(u)\theta $. Then 
\begin{eqnarray*}
\int_{\{\rho >\epsilon \}}\triangle _{b}^{\widetilde{\theta }_{0}}\dot{w}%
d\mu _{u^{-2}\theta } &=&-\int_{\{\rho =\epsilon \}}\widetilde{e}_{2}(\dot{w}%
)\widetilde{\theta }_{0}\wedge \widetilde{e}^{1} \\
&=&-\int_{\{\rho =\epsilon \}}ue_{2}(-\frac{\dot{u}}{u})u^{-3}\theta \wedge
e^{1} \\
&=&\int_{\{\rho =\epsilon \}}(u^{-3}e_{2}(\dot{u})-u^{-4}e_{2}(u)\dot{u}%
)\theta \wedge e^{1}.
\end{eqnarray*}%
Therefore, 
\begin{equation}
-4\int_{\{\rho >\epsilon \}}u^{-1}\dot{u}d\mu _{u^{-2}\theta }=\int_{\{\rho
=\epsilon \}}(u^{-3}e_{2}(\dot{u})-u^{-4}e_{2}(u)\dot{u})\theta \wedge
e^{1}+O(1),  \label{secondterm1}
\end{equation}%
here the function $u$ is $u_{0}$, i.e. the solution to the singular CR
Yamabe problem relative to $\Sigma =\Sigma _{0}$, that is 
\begin{equation*}
u=\rho +v\rho ^{2}+w\rho ^{3}+z\rho ^{4}+h\rho ^{5}+l\rho ^{5}\log \rho
+O(\rho ^{6}).
\end{equation*}

Note that the defining function relative to $\Sigma _{t}=F_{t}(\Sigma )$
satisfies 
\begin{equation*}
\rho _{t}(F_{t}x)=0,\quad x\in \Sigma _{0}
\end{equation*}%
hence by applying $\partial _{t}|_{t=0}$ to it, we find 
\begin{equation*}
\dot{\rho}(x)+(fe_{2}+gT)\rho =\dot{\rho}+(f-\alpha g)=0,\quad x\in \Sigma
_{0}
\end{equation*}%
hence 
\begin{equation*}
\dot{\rho}(x,\rho =0)=-(f-\alpha g),\quad x\in \Sigma _{0}.
\end{equation*}%
Note that $\dot{\rho}(x,\rho )$ is smooth in $\rho $. It then follows from 
\begin{equation*}
u_{t}=\rho _{t}+v_{t}\rho _{t}^{2}+w_{t}\rho _{t}^{3}+z_{t}\rho
_{t}^{4}+h_{t}\rho _{t}^{5}+l_{t}\rho _{t}^{5}\log \rho _{t}+O(\rho _{t}^{6})
\end{equation*}%
that 
\begin{equation*}
\dot{u}(x,\rho )=-k(x,\rho )-5(f-\alpha g)l\rho ^{4}\log \rho +O(\rho
^{5}\log \rho ),
\end{equation*}%
where $k(x,\rho )$ is a smooth function satisfying $k(x,0)=f-\alpha g$. Note
that 
\begin{equation*}
e_{2}=\partial _{\rho }+O(\rho )\mathcal{T}+O(\rho )\partial _{\rho },
\end{equation*}%
here $\mathcal{T}$ is a tangent vector of $\Sigma _{0,\epsilon }$, hence 
\begin{eqnarray*}
e_{2}(u) &=&1+O(\rho )+5l\rho ^{4}\log \rho +O(\rho ^{5}\log \rho ), \\
e_{2}(\dot{u}) &=&s(x,\rho )-20(f-\alpha g)l\rho ^{3}\log \rho +O(\rho
^{4}\log \rho )
\end{eqnarray*}%
for a smooth function $s(x,\rho )$. We have 
\begin{eqnarray*}
u(x,\rho )^{-3} &=&\rho ^{-3}(\lambda (x,\rho )+O(\rho ^{4}\log \rho )), \\
u(x,\rho )^{-4} &=&\rho ^{-4}(\mu (x,\rho )-4l\rho ^{4}\log \rho +O(\rho
^{5})),
\end{eqnarray*}%
where $\lambda $ and $\mu $ are smooth functions satisfying $\lambda
(x,0)=\mu (x,0)=1$. Therefore, 
\begin{eqnarray*}
&&\int_{\{\rho =\epsilon \}}u^{-3}e_{2}(\dot{u})\theta \wedge e^{1} \\
&=&\int_{\{\rho =\epsilon \}}\epsilon ^{-3}(\lambda (x,\epsilon )+O(\epsilon
^{4}\log \epsilon ))[s(x,\epsilon )-20(f-\alpha g)l\epsilon ^{3}\log
\epsilon +O(\epsilon ^{4}\log \epsilon )]\theta \wedge e^{1},
\end{eqnarray*}%
\begin{eqnarray*}
&&\int_{\rho =\epsilon }-u^{-4}e_{2}(u)\dot{u}\theta \wedge e^{1} \\
&=&-\int_{\rho =\epsilon }\epsilon ^{-4}(\mu (x,\epsilon )-4l\epsilon
^{4}\log \epsilon +O(\epsilon ^{5}))[1+O(\epsilon )+5l\epsilon ^{4}\log
\epsilon +O(\epsilon ^{5}\log \epsilon )] \\
&&\cdot \lbrack -k(x,\epsilon )-5(f-\alpha g)l\epsilon ^{4}\log \epsilon
+O(\epsilon ^{5}\log \epsilon )]\theta \wedge e^{1}
\end{eqnarray*}%
and the coefficients of $\log \frac{1}{\epsilon }$ in $\int_{\{\rho
=\epsilon \}}u^{-3}e_{2}(\dot{u})\theta \wedge e^{1}$ and $\int_{\rho
=\epsilon }-u^{-4}e_{2}(u)\dot{u}\theta \wedge e^{1}$ are 
\begin{eqnarray}
&&20\int_{\Sigma }(f-\alpha g)l\theta \wedge e^{1},  \label{secondtermlog1}
\\
&&-6\int_{\Sigma }(f-\alpha g)l\theta \wedge e^{1},  \label{secondtermlog2}
\end{eqnarray}%
respectively. It follows from (\ref{secondterm1}), (\ref{secondtermlog1})
and (\ref{secondtermlog2}) that the coefficient of $\log \frac{1}{\epsilon }$
in $-4\int_{\{\rho >\epsilon \}}u^{-1}\dot{u}d\mu _{u^{-2}\theta }$ is $%
14\int_{\Sigma }(f-\alpha g)l\theta \wedge e^{1}$. Then (\ref{log}) follows
from (\ref{dotL1}), (\ref{dotL2}) and (\ref{firsttermlog}).

%TCIMACRO{\TeXButton{End Proof}{\endproof}}%
%BeginExpansion
\endproof%
%EndExpansion

\bigskip

We would like to mention a simple example in the Heisenberg space $\mathcal{H%
}_{1}$. Let 
\begin{equation*}
\Sigma =\{r=\sqrt{x^{2}+y^{2}}=1,\quad t\in \mathbb{R}\},
\end{equation*}%
and $M=\{(x,y,t)\in \mathcal{H}|r\leq 1\}$. For $\Sigma $, 
\begin{equation*}
\alpha =0,\quad H=1.
\end{equation*}%
One can find a formal solution to the singular CR Yamabe problem. It turns
out that $\rho =1-r$, and 
\begin{equation*}
u=\rho -\frac{1}{6}\rho ^{2}-\frac{1}{9}\rho ^{3}-\frac{11}{108}\rho
^{4}+h\rho ^{5}+\frac{4}{135}\rho ^{5}\log \rho +O(\rho ^{6}).
\end{equation*}%
That is $l=\frac{4}{135}$. On the other hand, by (\ref{mathcalE}), $\mathcal{%
E}_{2}=\frac{4}{27}$.

\bigskip

\section{Appendix}

We collect some basic formulas in this appendix. The reader is referred to 
\cite{CHMY} for fundamental notions and notations. Let $(M,J,\theta )$ be a $%
3$-dimensional pseudohermitian manifold with boundary $\Sigma $ $=$ $%
\partial M$. We take an orthonormal (co)frame 
\begin{equation*}
\{e_{1},e_{2}=Je_{1},T\},\quad \{e^{1},e^{2},\theta \}
\end{equation*}%
with respect to the (adapted) metric 
\begin{equation*}
g_{\theta }=\theta \otimes \theta +\frac{1}{2}d\theta (\cdot ,J\cdot ).
\end{equation*}%
Here we use the convention 
\begin{equation*}
d\theta (X,Y)=X\theta (Y)-Y\theta (X)-\theta ([X,Y]).
\end{equation*}%
We have%
\begin{equation}
d\theta =2e^{1}\wedge e^{2}.  \label{dtheta}
\end{equation}%
Let $\nabla $ be the pseudohermitian (Tanaka-Webster) connection on $M$.
Recall (see the Appendix in \cite{CHMY}) that there is a real 1-form $\omega 
$ such that%
\begin{equation}
\nabla e_{1}=\omega \otimes e_{2},\quad \nabla e_{2}=-\omega \otimes e_{1},
\label{omega}
\end{equation}

\begin{eqnarray}
\nabla _{e_{1}}e_{1} &=&\omega (e_{1})e_{2},\quad \nabla
_{e_{2}}e_{2}=-\omega (e_{2})e_{1},  \label{e-12-1} \\
\nabla _{e_{2}}e_{1} &=&\omega (e_{2})e_{2},\quad \nabla
_{e_{1}}e_{2}=-\omega (e_{1})e_{1},  \notag \\
\lbrack e_{2},e_{1}] &=&\omega (e_{2})e_{2}+\omega (e_{1})e_{1}+2T.  \notag
\end{eqnarray}%
Let 
\begin{equation*}
T_{\nabla }(X,Y):=\nabla _{X}Y-\nabla _{Y}X-[X,Y]
\end{equation*}%
be the torsion of $\nabla $. Then we have 
\begin{equation*}
T_{\nabla }(e_{1},e_{2})=-\theta ([e_{1},e_{2}])T=2T.
\end{equation*}

\noindent Let 
\begin{equation*}
\tau (X):=T_{\nabla }(T,X).
\end{equation*}%
We have 
\begin{equation*}
\tau (Z_{1})=A_{1}^{\overline{1}}Z_{\overline{1}},\quad \tau (Z_{\overline{1}%
})=\overline{A_{1}^{\overline{1}}Z_{\overline{1}}}=A_{\overline{1}}^{1}Z_{1}
\end{equation*}%
where 
\begin{equation*}
Z_{1}=\frac{1}{2}(e_{1}-ie_{2}),\quad Z_{\overline{1}}=\frac{1}{2}%
(e_{1}+ie_{2}).
\end{equation*}

\noindent Write

\begin{equation*}
A_{11}=A_{1}^{\overline{1}}=a_{1}+ia_{2},\quad A_{\overline{11}}=A_{%
\overline{1}}^{1}=\overline{A_{11}}=a_{1}-ia_{2},
\end{equation*}
\noindent for some real valued functions $a_{1},a_{2}$, so that%
\begin{equation}
\tau (e_{1})=a_{1}e_{1}-a_{2}e_{2},\quad \tau (e_{2})=-a_{2}e_{1}-a_{1}e_{2}.
\label{tau}
\end{equation}%
Then we have%
\begin{eqnarray}
de^{1} &=&-e^{2}\wedge \omega +\theta \wedge (a_{1}e^{1}-a_{2}e^{2}),
\label{de} \\
de^{2} &=&e^{1}\wedge \omega -\theta \wedge (a_{2}e^{1}+a_{1}e^{2}),  \notag
\end{eqnarray}%
and%
\begin{equation}
d\omega =-Re^{1}\wedge e^{2}+(D_{1}e^{1}+D_{2}e^{2})\wedge \theta ,
\label{domega}
\end{equation}%
where%
\begin{eqnarray}
D_{1} &=&e_{2}(a_{1})+e_{1}(a_{2})-2a_{1}\omega (e_{1})+2a_{2}\omega (e_{2}),
\label{D12} \\
D_{2} &=&e_{1}(a_{1})-e_{2}(a_{2})+2a_{1}\omega (e_{2})+2a_{2}\omega (e_{1}),
\notag
\end{eqnarray}%
and%
\begin{eqnarray}
-R &=&-R^{\nabla }(e_{1},e_{2},e_{1},e_{2})=\langle \nabla _{e_{1}}\nabla
_{e_{2}}e_{1}-\nabla _{e_{2}}\nabla _{e_{1}}e_{1}-\nabla _{\lbrack
e_{1},e_{2}]}e_{1},e_{2}\rangle  \label{Websterscalar} \\
&=&e_{1}\omega (e_{2})-e_{2}\omega (e_{1})+\omega (e_{2})^{2}+\omega
(e_{1})^{2}+2\omega (T)  \notag \\
&=&-2W.  \notag
\end{eqnarray}

\noindent In the torsion free case, we also have%
\begin{eqnarray}
\lbrack T,e_{1}] &=&\omega (T)e_{2},  \label{e-13} \\
\lbrack e_{2},T] &=&\omega (T)e_{1},  \notag \\
d\omega &=&-2We^{1}\wedge e^{2}.  \notag
\end{eqnarray}

\end{document}